\documentclass[journal]{IEEEtran}

\usepackage{comment}
\usepackage{cite}
\usepackage{amsmath,amssymb,amsfonts,mathrsfs}
\usepackage{algorithm,algorithmic,setspace}
\usepackage{graphicx}
\usepackage{textcomp}
\usepackage{xcolor}
\usepackage{comment}
\usepackage{mathtools}
\usepackage{dsfont}
\usepackage{subcaption}
\usepackage{algorithm}
\usepackage{algorithmic}

\usepackage{makecell}
 
\usepackage{epstopdf} 
\usepackage{graphics} 
\usepackage{epstopdf} 
\usepackage{url}
\usepackage{color,soul}
\usepackage{amsthm}
\usepackage[colorlinks=true,linkcolor=blue, citecolor=green, urlcolor=magenta]{hyperref}

\theoremstyle{definition}

\newtheorem{theorem}{Theorem}
\newtheorem{assumption}{Assumption}

\newtheorem{lemma}{Lemma}
\newtheorem{remark}{Remark}

\newtheorem{fact}{Fact}

\AtEndEnvironment{remark}{\hfill\ensuremath{\square}}

\makeatletter
\@ifundefined{proof}{
}{%
  \let\oldproof\proof
  \renewcommand{\proof}[1]{\oldproof{#1\hfill\ensuremath{\blacksquare}}}
}
\makeatother

\newcommand{\Lim}{\displaystyle\lim}
\allowdisplaybreaks

\title{Constrained Stabilization on the $n$-Sphere with Conic and Star-shaped Constraints
} 
\author{Mayur Sawant and Abdelhamid Tayebi  
	\thanks{This work was supported by the Natural Sciences and Engineering Research Council of Canada (NSERC), under the grant RGPIN-2026-06037.}
	\thanks{M. Sawant and A. Tayebi are with the Department of Electrical and Computer Engineering, Lakehead University, Thunder Bay, ON P7B 5E1, Canada. (e-mail: {\tt\small msawant, atayebi@lakeheadu.ca}).}%
}%
\begin{document}

\maketitle

\begin{abstract}
   The problem of constrained stabilization on the $n$-sphere under star-shaped constraints is considered. We propose a control strategy that allows one to almost globally steer the state to a desired location while avoiding star-shaped constraints on the $n$-sphere. Depending on the state's proximity to the unsafe regions, the state is either guided towards the target location along the geodesic connecting the target to the state or steered towards the antipode of a predefined point lying in the interior of the nearest unsafe region.
    We prove that the target location is almost globally asymptotically stable under the proposed continuous, time-invariant feedback control law. Nontrivial simulation results on the $2$-sphere and the $3$-sphere demonstrate the effectiveness of the theoretical results.
\end{abstract}

\section{Introduction}
Various mechanical systems have states that evolve on the $n$-sphere. Examples include spin-axis stabilization of rigid body systems \cite{bullo1995control}, two-axis gimbal systems \cite{osborne2008global}, thrust-vector control for quad-rotor aircraft \cite{hua2015control}, and the spherical robot \cite{muralidharan2015geometric}. In many practical scenarios, the attitude stabilization problem can also be recast as a stabilization on the $3$-sphere.

The stabilization problem on the $n$-sphere (without constraints) has been dealt with in the literature using differential geometry and hybrid dynamical systems tools, see for instance \cite{bullo1995control, casau2019hybrid, casau2019robust}.
In \cite{lee2014feedback}, a logarithmic barrier function is used to design a quaternion-based feedback controller for attitude control of a rigid body spacecraft in the presence of multiple attitude-constrained zones, characterized by quadratic inequalities.
Another logarithmic barrier function based approach for attitude stabilization on the special orthogonal group $\mathrm{SO}(3)$ under conic constraints is proposed in \cite{kulumani2017constrained}. In \cite{nicotra2019spacecraft}, the authors proposed an explicit reference governor approach for spacecraft attitude control under actuator saturation and conic constraints. 
In \cite{danielson2021spacecraft}, an invariant set motion planner is proposed to plan a sequence of reference quaternion waypoints that safely guides the spacecraft attitude to a desired orientation while avoiding unsafe regions---modeled as conic constraints.
In \cite{su2024practical}, the problem of spacecraft attitude reorientation under conic constraints and physical limitations is addressed by designing a virtual angular velocity, relying on control barrier functions to ensure constraint satisfaction. A prescribed performance controller is then designed for the angular velocity tracking while taking into account the control input saturation.
In \cite{berkane2021constrained}, the authors addressed the stabilization problem on the $n$-sphere under conic constraints by leveraging the stereographic projection to transform the problem into a classical navigation problem in $\mathbb{R}^n$ with spherical obstacles, enabling the use of existing navigation function-based obstacle avoidance methods.
Reference \cite{liu2023adaptive} investigates the problem of reduced attitude control for a rigid spacecraft under elliptical pointing constraints and parameter uncertainties. Employing a diffeomorphic projection and elliptical stereographic mapping, the problem is reformulated as an obstacle avoidance problem in a two-dimensional Euclidean space.

Although these approaches guarantee constrained stabilization on the spherical manifold, in most cases, the characterization of unsafe sets is limited to conic constraints. 
Since the $n$-sphere is a bounded manifold, a more flexible characterization of the unsafe region can result in a larger safe region for stabilization purposes.


In this paper, we design a continuous feedback control law for almost\footnote{An equilibrium point is almost globally asymptotically stable if it is stable and attractive from all initial conditions except a set of zero Lebesgue measure.} global asymptotic stabilization on the $n$-sphere while avoiding star-shaped constraints.
Note that geodesically strongly convex constraints \cite[Chap. IV, Def. 5.1]{sakai1996riemannian},
such as conic constraints, each contained in an open hemisphere on the $n$-sphere,
form a subset of the star-shaped constraints.
Inspired by the obstacle avoidance strategy in \cite{kumar2022navigation}, where the state is steered radially away from the center of an ellipsoidal obstacle in the Euclidean space $\mathbb{R}^n$, the proposed feedback controller steers the state, depending on its proximity to unsafe regions, towards the antipode of a predefined point in the interior of the nearest star-shaped set on the $n$-sphere.



The main contributions of the proposed work are as follows:
\begin{enumerate}
    \item \textit{Safety and almost global asymptotic stability:} The proposed control strategy ensures safety and guarantees almost global asymptotic stabilization of the desired location on the $n$-sphere under star-shaped constraints. To the best of the authors' knowledge, this is the first work in the literature achieving such strong stability results for the constrained stabilization problem on the $n$-sphere with star-shaped constraints.
    \item \textit{Arbitrarily-shaped star-shaped constraints on the $n$-sphere:} The proposed feedback controller is able to handle star-shaped constraints on the $n$-sphere. 
    Note that geodesically strongly convex constraints \cite[Chap. IV, Def. 5.1]{sakai1996riemannian}, 
    such as conic constraints, each contained in an open hemisphere of the $n$-sphere,
    form a subset of the star-shaped constraints.
    \item \textit{Minimal constraint information required:} The proposed feedback controller does not require complete knowledge of the constraint set. It only requires (i) at least one point in the interior of each constraint such that the geodesic connecting any point of the set to it lies entirely within the set, and (ii) a means of measuring proximity to the set in terms of spherical separation, as defined later in Section~\ref{notations}.
\end{enumerate}

The rest of the paper is organized as follows. 
Section \ref{notations} introduces the notations and mathematical preliminaries used throughout the paper, and Section \ref{section:n-sphere-problem-formulation} specifies the problem statement.
In Section \ref{section:conic_constraint_avoidance}, we present a feedback control design for stabilization on the $n$-sphere under conic constraints.
This controller is then modified to address the problem of stabilization on the $n$-sphere under the star-shaped constraints in Section \ref{section:feedback_control_design}.
In Section \ref{section:stability_analysis}, we analyze the safety and stability properties of the resulting closed-loop system.
In Section \ref{section:application}, the proposed controllers are applied to the problem of constrained (partial and full) attitude stabilization, and their effectiveness is demonstrated through non-trivial simulation studies. Finally, concluding remarks are provided in Section \ref{section:conclusion}.

\section{Notations and preliminaries}\label{notations}
The sets of real numbers and natural numbers are represented by $\mathbb{R}$ and $\mathbb{N}$, respectively.
Bold lowercase letters are used to represent vector quantities.
The Euclidean norm of any vector $\mathbf{x}\in\mathbb{R}^n$ is given by $\|\mathbf{x}\| = \sqrt{\mathbf{x}^\top\mathbf{x}}$.
The identity matrix and the zero matrix of dimension $n\in\mathbb{N}$ are denoted by $\mathbf{I}_n$ and $\mathbf{0}_n$, respectively.
Given $\mathcal{A}\subset\mathbb{R}^n$ and $\mathcal{B}\subset\mathbb{R}^n$, the relative complement of $\mathcal{B}$ in $\mathcal{A}$ is given by $\mathcal{A}\setminus\mathcal{B} = \{\mathbf{a}\in\mathcal{A}\mid\mathbf{a}\notin\mathcal{B}\}$.
Given $\mathcal{A}\subset\mathbb{R}^n$, the cardinality of $\mathcal{A}$ is denoted by $\mathrm{card}(\mathcal{A})$.
For a twice differentiable scalar mapping $f:\mathbb{R}\to\mathbb{R}$, we denote its first and second derivatives by $f'(x) = \frac{df(x)}{dx}$ and $f''(x) = \frac{d^2f(x)}{dx^2}$, respectively.

We also define the following subsets of $\mathbb{R}^n$:\\
{\bf Line segment:} Given any two points $\mathbf{a}, \mathbf{b}\in\mathbb{R}^n$, the line segment $\mathcal{L}_s(\mathbf{a}, \mathbf{b})$ joining $\mathbf{a}$ and $\mathbf{b}$ is defined as
\begin{equation}\label{expression_for_line_segment}
    \mathcal{L}_s(\mathbf{a}, \mathbf{b}) = \{\mathbf{x}\in\mathbb{R}^n\mid\mathbf{x} = (1-\lambda)\mathbf{a} + \lambda\mathbf{b}, \lambda\in[0, 1]\}.
\end{equation}
{\bf Convex cone:} Given $\mathbf{a}\in\mathbb{R}^n\setminus\{\mathbf{0}\}$ and $\mathbf{b}\in\mathbb{R}^n\setminus\{\mathbf{0}\}$, a convex cone $\mathcal{C}(\mathbf{a}, \mathbf{b})$ with its vertex at the origin is defined as
\[
\mathcal{C}(\mathbf{a}, \mathbf{b}) = \{\mathbf{x}\in\mathbb{R}^n\mid \mathbf{x} = \lambda_1\mathbf{a} + \lambda_2\mathbf{b}, \lambda_1\geq 0, \lambda_2\geq 0\}.
\]
In the present work, we consider the motion on the unit $n$-sphere which is an $n$-dimensional manifold embedded in the Euclidean space $\mathbb{R}^{n+1}$ and defined as $\mathbb{S}^n:=\{\mathbf{x}\in\mathbb{R}^{n+1}\mid\|\mathbf{x}\| = 1\}$.
Given a set $\mathcal{A}\subset\mathbb{S}^n$, the symbols $\overline{\mathcal{A}}, \mathcal{A}^{\circ}$, and $\partial\mathcal{A}$ represent the closure, interior, and the boundary of $\mathcal{A}$ on $\mathbb{S}^n$, where $\partial\mathcal{A} = \overline{\mathcal{A}}\setminus\mathcal{A}^{\circ}$.

In the following, we will provide the definitions of some concepts that will be used throughout the paper.\\
{\bf Tangent space:} The tangent space to $\mathbb{S}^n$ at $\mathbf{x}\in\mathbb{S}^n$ is given by $\mathbf{T}_{\mathbf{x}}(\mathbb{S}^n) = \{\mathbf{a}\in\mathbb{R}^{n+1}\mid\mathbf{a}^\top\mathbf{x} = 0\}$, which represents all vectors in $\mathbb{R}^{n+1}$ that are perpendicular to $\mathbf{x}$.
Given $\mathbf{x}\in\mathbb{S}^n$ and $\mathbf{a}\in\mathbb{R}^{n+1}$, the orthogonal projection operator $\mathbf{P}(\mathbf{x})$, which is given by 
\begin{equation}\label{orthogonal_projection_operator_formula}
    \mathbf{P}(\mathbf{x}) = \mathbf{I}_{n+1} - \mathbf{x}\mathbf{x}^\top,
\end{equation}
projects $\mathbf{a}$ onto the tangent space $\mathbf{T}_{\mathbf{x}}(\mathbb{S}^n)$, \textit{i.e.}, $\mathbf{P}(\mathbf{x})\mathbf{a}\in\mathbf{T}_{\mathbf{x}}(\mathbb{S}^n)$.\\

\noindent{\bf Spherical separation:} {Given a set $\mathcal{A}\subset\mathbb{S}^n$ and $\mathbf{x}\in\mathbb{S}^n$, the separation between $\mathbf{x}$ and $\mathcal{A}$ is defined as
\begin{equation}\label{definition:spherical_distance_function}
    d_s(\mathbf{x}, \mathcal{A}) := \underset{\mathbf{a}\in\mathcal{A}}{\inf\;}(1- \mathbf{x}^\top\mathbf{a}).
\end{equation}}
Similarly, the separation between any two unit vectors $\mathbf{x}, \mathbf{a}\in\mathbb{S}^n$ is defined as $d_s(\mathbf{x}, \mathbf{a}) = 1 - \mathbf{x}^\top\mathbf{a}$.
We define $\mathcal{P}(\mathbf{x}, \mathcal{A})$ as the set containing the points $\mathbf{a}$ in $\mathcal{A}$ that satisfy $d_s(\mathbf{x}, \mathbf{a}) = d_s(\mathbf{x}, \mathcal{A})$ \textit{i.e.,}
\begin{equation}\label{projection_set_definition}
    \mathcal{P}(\mathbf{x}, \mathcal{A}) = \{\mathbf{a}\in\mathcal{A}\mid d_s(\mathbf{x}, \mathbf{a}) = d_s(\mathbf{x}, \mathcal{A})\}.
\end{equation}
If $\mathrm{card}(\mathcal{P}(\mathbf{x}, \mathcal{A})) = 1$, then the unique element in $\mathcal{P}(\mathbf{x}, \mathcal{A})$ is represented by $\Pi(\mathbf{x}, \mathcal{A})$.

Given a set $\mathcal{A}\subset\mathbb{S}^n$, the dilation of $\mathcal{A}$ by $p > 0$ on $\mathbb{S}^n$ is defined as
\begin{equation}\label{dilation_on_sphere}
    \mathcal{D}_p(\mathcal{A}) = \{\mathbf{x}\in\mathbb{S}^n\mid d_s(\mathbf{x}, \mathcal{A}) \leq p\}.
\end{equation}
Furthermore, the $p$-neighborhood of $\mathcal{A}$ on $\mathbb{S}^n$ is given by $\mathcal{N}_{p}(\mathcal{A}) = \mathcal{D}_{p}(\mathcal{A})\setminus\mathcal{A}^{\circ}$.\\
\noindent{\bf Geodesic:} For any two distinct points $\mathbf{a}, \mathbf{b}\in\mathbb{S}^n$ with $\mathbf{a} \ne -\mathbf{b}$, the unique geodesic connecting $\mathbf{a}$ and $\mathbf{b}$ is given by 
\begin{equation}\label{geodesic_expression}
    \mathcal{G}(\mathbf{a}, \mathbf{b}) = \left\{ \mathbf{x} \in \mathbb{S}^n \mid \mathbf{x} = g(\lambda; \mathbf{a}, \mathbf{b}), \lambda\in[0, 1]
    \right\},
\end{equation}
where, motivated by \cite[Section 3.3]{shoemake1985animating}, the mapping $g:[0, 1]\to\mathbb{S}^n$ is defined as
\begin{equation}\label{geodesic_parameterization}
    g(\lambda; \mathbf{a}, \mathbf{b}) = \frac{\sin((1-\lambda)\theta) \mathbf{a} + \sin(\lambda\theta) \mathbf{b}}{\sin\theta},
\end{equation}
where $\theta = \arccos(\mathbf{a}^\top\mathbf{b})$.
Since $\mathbf{P}(g(\lambda; \mathbf{a}, \mathbf{b}))\frac{d^2g(\lambda;\mathbf{a}, \mathbf{b})}{d\lambda^2} = \mathbf{0}$ for all $\lambda \in [0, 1]$, using  \cite[Chap. 3, Def. 2.1]{do1992riemannian}, one can confirm that $\mathcal{G}(\mathbf{a}, \mathbf{b})$ is a geodesic and is the curve on $\mathbb{S}^n$ with the smallest path length, connecting $\mathbf{a}$ and $\mathbf{b}$.
If $\mathbf{a} = \mathbf{b}$, then the geodesic $\mathcal{G}(\mathbf{a}, \mathbf{b})$ is trivially the point itself, \textit{i.e.,} $\mathcal{G}(\mathbf{a}, \mathbf{a}) = \{\mathbf{a}\}$.

\noindent {\bf Star-shaped sets on $\mathbb{S}^n$:} A set $\mathcal{A} \subset \mathbb{S}^n$ is a star-shaped set on $\mathbb{S}^n$ if there exists $\mathbf{g} \in \mathcal{A}$ with $-\mathbf{g} \notin \mathcal{A}$ such that $\mathcal{G}(\mathbf{g}, \mathbf{x}) \subset \mathcal{A}$ for all $\mathbf{x} \in \mathcal{A}$.

Given a star-shaped set $\mathcal{A}$ on $\mathbb{S}^n$, let $\sigma(\mathcal{A})$ be the set of all points $\mathbf{g}$ in $\mathcal{A}$ such that $-\mathbf{g}\notin\mathcal{A}$ and $\mathcal{G}(\mathbf{g}, \mathbf{x})\subset\mathcal{A}$ for all $\mathbf{x}\in\mathcal{A}$, defined as follows:
\begin{equation}\label{sigma_set}
    \sigma(\mathcal{A}) = \{\mathbf{g}\in\mathcal{A}\mid-\mathbf{g}\notin\mathcal{A},~\forall \mathbf{x}\in\mathcal{A},~\mathcal{G}(\mathbf{g}, \mathbf{x})\subset\mathcal{A}\}.
\end{equation}

Notice that for every point $\mathbf{g}\in\sigma(\mathcal{A})\cap\mathcal{A}^{\circ}$, the geodesic $\mathcal{G}(\mathbf{x}, -\mathbf{g})$ connecting any point $\mathbf{x}$ on the boundary of $\mathcal{A}$ on $\mathbb{S}^n$ to $-\mathbf{g}$ does not intersect with the interior of $\mathcal{A}$ on $\mathbb{S}^n$, as stated in the next lemma.
\begin{lemma}\label{lemma:reverse_geodesic_always_stay_outside}
Let $\mathcal{A}$ be a star-shaped set on $\mathbb{S}^n$. Then, for every $\mathbf{g}\in\sigma(\mathcal{A})\cap\mathcal{A}^{\circ}$ and for all $\mathbf{x}\in\partial\mathcal{A}$, one has
\[\mathcal{G}(\mathbf{x}, -\mathbf{g})\cap\mathcal{A}^{\circ} = \emptyset.\]
\end{lemma}
\proof{See Appendix \ref{proof:lemma:reverse_geodesic_always_stay_outside}.}

\begin{figure}[ht]
\centering
\subfloat[][]{\includegraphics[width =0.47\linewidth]{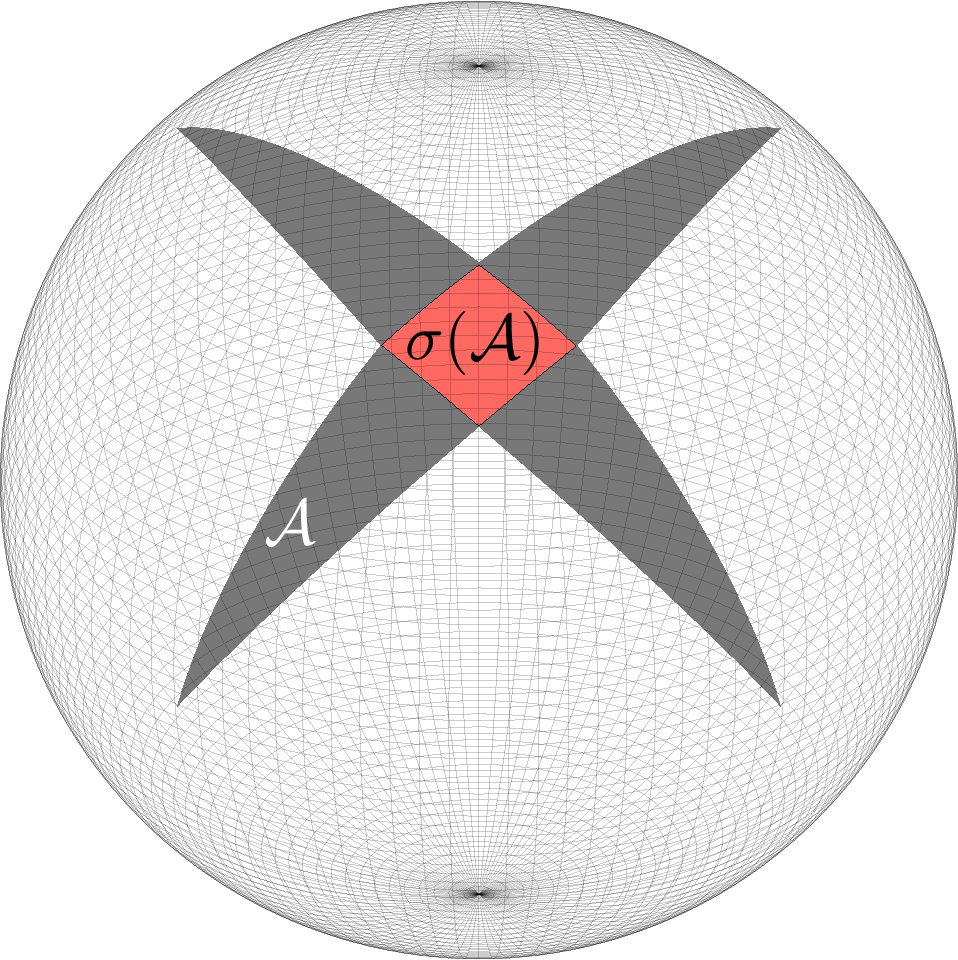}\label{g-star-set-example}}\hspace{0.2cm}
\subfloat[][]{\includegraphics[width =0.47\linewidth]{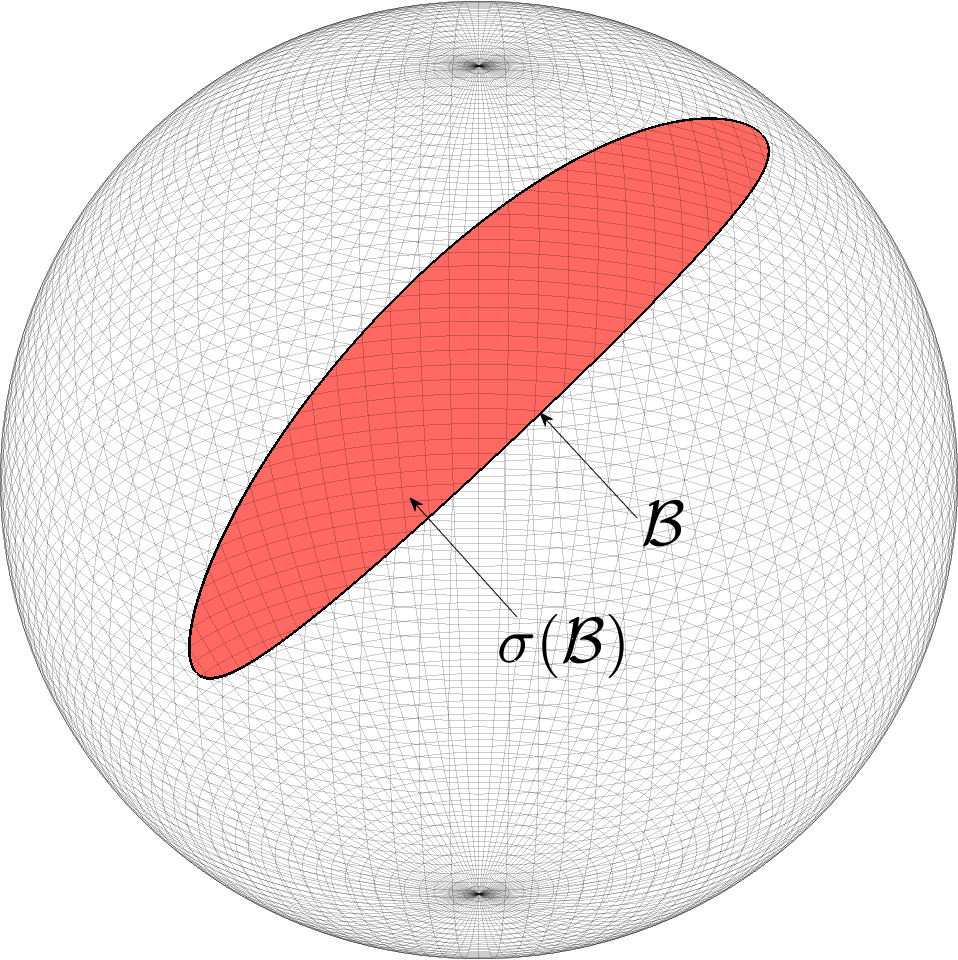}\label{gs-convex-set-example}}
\caption{Illustration of (a) a star-shaped set and (b) a gs-convex set on $\mathbb{S}^n$.}\label{illustration_of_g-star_gs-convex}
\end{figure}

\begin{remark}\label{remark:gs-convex-is-g-star}
    Every geodesically strongly convex (gs-convex) set $\mathcal{A} \subset \mathbb{S}^n$ is a star-shaped set on $\mathbb{S}^n$.
    A set $\mathcal{A} \subset \mathbb{S}^n$ is said to be gs-convex if, for any two points $\mathbf{a}, \mathbf{b} \in \mathcal{A}$, the unique geodesic $\mathcal{G}(\mathbf{a}, \mathbf{b})$ connecting $\mathbf{a}$ and $\mathbf{b}$ lies entirely in $\mathcal{A}$, that is, $\mathcal{G}(\mathbf{a}, \mathbf{b}) \subset \mathcal{A}$ for all $\mathbf{a}, \mathbf{b} \in \mathcal{A}$ \cite[Chap. IV, Def. 5.1]{sakai1996riemannian}. Consequently, if $\mathcal{B}$ is gs-convex, then it is a star-shaped set on $\mathbb{S}^n$ and $\sigma(\mathcal{B}) = \mathcal{B}$, as illustrated in Fig.~\ref{gs-convex-set-example}, where $\sigma(\mathcal{B})$ is defined as in \eqref{sigma_set}.
\end{remark}

\section{Problem formulation}\label{section:n-sphere-problem-formulation}
We consider the problem of constrained stabilization on $\mathbb{S}^n$ for the system
\begin{equation}\label{system_dynamics}
    \dot{\mathbf{x}} = \mathbf{P}(\mathbf{x})\mathbf{u},
\end{equation}
where $\mathbf{x}\in\mathbb{S}^n$ is the state vector, $\mathbf{u}\in\mathbb{R}^{n+1}$ is the control input, and $n\geq 2$.
The orthogonal projection operator $\mathbf{P}(\mathbf{x})$, defined in \eqref{orthogonal_projection_operator_formula}, projects $\mathbf{u}$ onto the tangent space to $\mathbb{S}^n$ at $\mathbf{x}$.
In other words, $\mathbf{P}(\mathbf{x})$ ensures that $\dot{\mathbf{x}}\in\mathbf{T}_{\mathbf{x}}(\mathbb{S}^n)$ for all $\mathbf{x}\in\mathbb{S}^n$, implying that $\mathbf{x}^\top\dot{\mathbf{x}} = 0$.
Consequently, if $\mathbf{x}(0)\in\mathbb{S}^n$, then $\mathbf{x}(t)\in\mathbb{S}^n$ for all future times.

The objective is to stabilize $\mathbf{x}$ at the desired point $\mathbf{x}_d\in\mathbb{S}^n$, while avoiding the interior of an unsafe region $\mathcal{U}\subset\mathbb{S}^n$.
The set $\mathcal{U}$, defined as the union of $m$ closed sets $\mathcal{U}_i$ on $\mathbb{S}^n$, where $i\in\{1, \ldots, m\}=:\mathbb{I}$ and $m\in\mathbb{N}$, is given by
\begin{equation}\label{combined-unsafe-set}
    \mathcal{U} = \bigcup_{i\in\mathbb{I}}\mathcal{U}_i.
\end{equation}

For safe stabilization, the condition $\mathbf{x}(t)\notin\mathcal{U}^{\circ}$ should hold for all $t\geq 0$. 
Defining the set 
\begin{equation}\label{eroded_free_space_definition}
    \mathcal{M}_p = \{\mathbf{x}\in\mathbb{S}^n
    \setminus\mathcal{U}^{\circ}\mid d_s(\mathbf{x}, \mathcal{U}) \geq p\},
\end{equation}
for $p\geq 0$, safe stabilization is, therefore, ensured if and only if $\mathbf{x}(t)\in\mathcal{M}_0$ for all $t\geq 0$.

In Section \ref{section:conic_constraint_avoidance}, the unsafe regions $\mathcal{U}_i$ represent conic constraints, whereas Section \ref{section:feedback_control_design} considers them to be star-shaped on $\mathbb{S}^n$.
To ensure the feasibility of the problem, we assume that the sets $\mathcal{U}_i$, where $i\in\mathbb{I}$, do not overlap with each other, as stated in the following assumption:
\begin{assumption}\label{assumption:non-overlapping-constraints}
    For any two distinct sets $\mathcal{U}_i$ and $\mathcal{U}_j$, the cosine of the angle between any pair of unit vectors $\mathbf{a}\in\mathcal{U}_i$ and $\mathbf{b}\in\mathcal{U}_j$ is always less than or equal to $\cos(2\delta)$, where $\delta \in \big(0, \frac{\pi}{2}\big]$. 
    In other words, for $i, j\in\mathbb{I}, i\ne j$, 
    \begin{equation*}
    \underset{\mathbf{a}\in\mathcal{U}_i, \mathbf{b}\in\mathcal{U}_j}{\max}\mathbf{a}^\top\mathbf{b} \leq \cos(2\delta).
    \end{equation*}
\end{assumption}
According to Assumption \ref{assumption:non-overlapping-constraints}, the smallest positive angle between any pair of unit vectors $\mathbf{a}\in\mathcal{U}_i$ and $\mathbf{b}\in\mathcal{U}_j$ is always greater than or equal to $2\delta$, where $i, j\in\mathbb{I}, i\ne j$ and $\delta\in\big(0, \frac{\pi}{2}\big]$.

The task is to design $\mathbf{u}$ in \eqref{system_dynamics} such that the following objectives are satisfied:
\begin{enumerate}
    \item The set $\mathcal{M}_0$, defined according to  \eqref{eroded_free_space_definition}, is forward invariant. 
    That is, if $\mathbf{x}(0)\in\mathcal{M}_0$, then $\mathbf{x}(t)\in\mathcal{M}_0$ for all $t\geq 0$.
    \item The target location $\mathbf{x}_d\in\mathcal{M}_0^{\circ}$ is almost globally asymptotically stable\footnote{The equilibrium $\mathbf{x}_d\in\mathcal{M}_0^{\circ}$ is stable and attractive from all initial conditions in $\mathcal{M}_0$ except a set of zero Lebesgue measure.} over $\mathcal{M}_0$.
\end{enumerate}

\section{Constrained stabilization under conic constraints}
\label{section:conic_constraint_avoidance}
\begin{figure}[ht]
    \centering
    \includegraphics[width=0.6\linewidth]{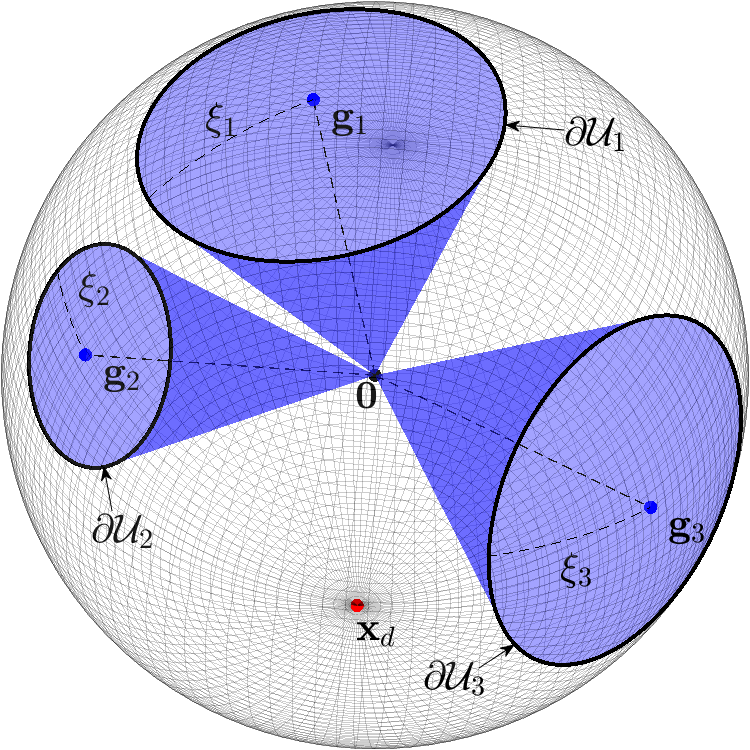}
    \caption{Conic constraints \eqref{definition:conic_constraints}.}
    \label{fig:conic_constraint_depiction}
\end{figure}

For each $i\in\mathbb{I}$, a conic constraint $\mathcal{U}_i$ on $\mathbb{S}^n$ is defined as 
\begin{equation}\label{definition:conic_constraints}
    \mathcal{U}_i = \{\mathbf{x}\in\mathbb{S}^n\mid \mathbf{x}^\top\mathbf{g}_i \geq \cos(\xi_i)\},
\end{equation}
where $\mathbf{g}_i\in\mathbb{S}^n\setminus\{\mathbf{x}_d\}$ and $\xi_i \in (0, \pi)$.
The constant unit vectors $\mathbf{g}_i$ and the scalar parameters $\xi_i$ are set such that the unsafe regions $\mathcal{U}_i$ satisfy Assumption \ref{assumption:non-overlapping-constraints}, as illustrated in Fig. \ref{fig:conic_constraint_depiction}.

A conic constraint $\mathcal{U}_i \subset \mathbb{S}^n$ plays a role analogous to that of a spherical obstacle in classical navigation problems posed in Euclidean space. 
The conic constraint $\mathcal{U}_i$ represents a region on $\mathbb{S}^n$ that excludes all unit vectors whose angular separation from $\mathbf{g}_i$ is not more than $\xi_i$ in much the same way that a spherical obstacle excludes a neighborhood of a point in $\mathbb{R}^n$.

Therefore, inspired by the navigation function formulation for the sphere world in \cite{koditschek1990robot}, we construct a scalar function $W(\mathbf{x})$ such that steering $\mathbf{x}$ on $\mathbb{S}^n$ along the orthogonal projection of its negative gradient ensures safety and almost global asymptotic stability of the point $\mathbf{x}_d \in \mathbb{S}^n$ over the free space $\mathcal{M}_0$.
The scalar function $W(\mathbf{x})$ is given by
\begin{equation}\label{navigation_function_on_sphere}
    W(\mathbf{x}) = \frac{k_1d_s(\mathbf{x}, \mathbf{x}_d)}{d_s(\mathbf{x}, \mathbf{x}_d) + \beta(\mathbf{x})},
\end{equation}
where $k_1 > 0$, and $d_s(\mathbf{x},\mathbf{x}_d)$ is defined as in \eqref{definition:spherical_distance_function}.
The scalar function $\beta(\mathbf{x})$ is defined as
\begin{equation}\label{obstacle-function-defnition}
    \beta(\mathbf{x}) = \begin{cases}
        h_i(d_s(\mathbf{x}, \mathbf{g}_i)), & \mathbf{x}\in\mathcal{N}_{\epsilon}(\mathcal{U}_i),\\
        1, & \mathbf{x}\notin\mathcal{N}_{\epsilon}(\mathcal{U}),
    \end{cases}
\end{equation}
where the unit vector $\mathbf{g}_i$ is defined in \eqref{definition:conic_constraints}.

For each $i\in\mathbb{I}$ with $0\leq \underline{\Delta}_i < \overline{\Delta}_i$, the scalar mapping $h_i:\left[\underline{\Delta}_i, \overline{\Delta}_i\right]\to[0, 1]$ in \eqref{obstacle-function-defnition} is strictly increasing and real analytic
over $\left[\underline{\Delta}_i, \overline{\Delta}_i\right]$, and satisfies the following properties: $h_i(\underline{\Delta}_i) = 0$, $h_i(\overline{\Delta}_i) = 1$, $h_i'(\overline{\Delta}_i) = 0$, $h_i''(\overline{\Delta}_i) = 0$, and $h_i'(p) > 0, \forall p\in[\underline{\Delta}_i, \overline{\Delta}_i)$.\footnote{An example of such a function is $h_i(p) = 1 + \frac{\left(p - \overline{\Delta}_i\right)^3}{\left(\overline{\Delta}_i - \underline{\Delta}_i\right)^3}$. 
Since its derivative, $h_i'(p) = \frac{3\left(p - \overline{\Delta}_i\right)^2}{\left(\overline{\Delta}_i - \underline{\Delta}_i\right)^3}$, is positive for all $p\in\big[\underline{\Delta}_i, \overline{\Delta}_i\big)$, $h_i(p)$ is strictly increasing over $\left[\underline{\Delta}_i, \overline{\Delta}_i\right]$.}
For each $i\in\mathbb{I}$, the scalar parameters $\underline{\Delta}_i$ and $\overline{\Delta}_i$ are the values of $d_s(\mathbf{x}, \mathbf{g}_i)$ on the constraint boundary $\partial\mathcal{U}_i$ and on the outer boundary of the $\epsilon$-neighborhood $\partial\mathcal{N}_{\epsilon}(\mathcal{U}_i)\setminus\partial\mathcal{U}_i$, respectively.\footnote{For each $i\in\mathbb{I}$, $d_s(\mathbf{x}, \mathbf{g}_i)$ is constant on the constraint boundary $\partial\mathcal{U}_i$ and the outer boundary of the $\epsilon$-neighborhood $\partial\mathcal{N}_{\epsilon}(\mathcal{U}_i)\setminus\partial\mathcal{U}_i$, taking the distinct values $\underline{\Delta}_i = 1 - \cos(\xi_i)$ and $\overline{\Delta}_i = 1 - \cos\left(\xi_i + \arccos(1 - \epsilon)\right)$, respectively, where for any $p\in[-1, 1]$, the inverse cosine function $\arccos(p)$ is restricted to $[0, \pi]$.}

The parameter $\epsilon$ is chosen as $\epsilon\in\left(0, \min\{\bar{\epsilon}, 1-\cos(\delta)\}\right)$, where $\delta$ is defined in Assumption \ref{assumption:non-overlapping-constraints}, and $\bar{\epsilon}$ is a strictly positive scalar such that $\mathbf{x}_d\notin\mathcal{N}_{\bar{\epsilon}}(\mathcal{U})$.
Since $\mathbf{x}_d\in\mathcal{M}_0^{\circ}$, one has $\mathbf{x}_d\notin\mathcal{U}$, and the existence of $\bar{\epsilon} > 0$ such that $\mathbf{x}_d\notin\mathcal{N}_{\bar{\epsilon}}(\mathcal{U})$ is straightforward to establish.
The index $i$ in \eqref{obstacle-function-defnition} refers to the closest unsafe region $\mathcal{U}_i$ such that $\mathbf{x}\in\mathcal{N}_{\epsilon}(\mathcal{U}_i)$, where setting $\epsilon < 1 - \cos(\delta)$ ensures that the regions $\mathcal{D}_{\epsilon}(\mathcal{U}_i)$ and $\mathcal{D}_{\epsilon}(\mathcal{U}_j)$ are disjoint for every $i, j\in\mathbb{I}$ with $i\ne j$, \textit{i.e.}, $\mathcal{D}_{\epsilon}(\mathcal{U}_i)\cap\mathcal{D}_{\epsilon}(\mathcal{U}_j) = \emptyset$, as stated in the next lemma.
\begin{lemma}\label{lemma:constraint_separation}
    Let Assumption \ref{assumption:non-overlapping-constraints} hold for the free space $\mathcal{M}_0$, where $\mathcal{M}_0$ is obtained by replacing $p$ with $0$ in \eqref{eroded_free_space_definition}. If $\epsilon \in (0,  1 - \cos(\delta))$, then $\mathcal{D}_{\epsilon}(\mathcal{U}_i)\cap\mathcal{D}_{\epsilon}(\mathcal{U}_j) = \emptyset$ for every $i, j\in\mathbb{I}$ with $i\ne j$.
\end{lemma}
\proof{See Appendix \ref{proof:lemma:constraint_separation}}.

In the next lemma, we state some of the key properties of the scalar function $W(\mathbf{x})$ defined in \eqref{navigation_function_on_sphere}.
\begin{lemma}\label{lemma:properties_scalar_function}
    The scalar mapping $W:\mathcal{M}_0\to[0, k_1]$ defined in \eqref{navigation_function_on_sphere} satisfies the following properties:
    \begin{enumerate}
        \item \label{claim1:properties}It is positive definite with respect to $\mathbf{x}_d$ on $\mathcal{M}_0$.
        \item \label{claim2:properties}It attains its maximum value of $k_1$ if and only if $\mathbf{x}\in\partial\mathcal{M}_0$.
        \item \label{claim3:properties}It is twice continuously differentiable on $\mathcal{M}_0$.
    \end{enumerate}
\end{lemma}
\proof{See Appendix \ref{proof:lemma:properties_scalar_function}.}

The proposed feedback control law is the negative gradient of $W(\mathbf{x})$ with respect to $\mathbf{x}$ and is given as
\begin{equation}\label{negative_gradient_control_law}
    \mathbf{u}(\mathbf{x}) = -\nabla_{\mathbf{x}}W(\mathbf{x}).
\end{equation}

In the next theorem, we show that for the closed-loop system \eqref{system_dynamics}-\eqref{negative_gradient_control_law}, the set $\mathcal{M}_0$ is forward invariant and the desired point $\mathbf{x}_d$ is almost globally asymptotically stable over $\mathcal{M}_0$.
\begin{theorem}\label{maintheorem:conic_constraints}
    For the closed-loop system \eqref{system_dynamics}-\eqref{negative_gradient_control_law} under Assumption \ref{assumption:non-overlapping-constraints}, the following statements are valid:
    \begin{enumerate}
        \item \label{claim:conic1}The set $\mathcal{M}_0$ is forward invariant, where $\mathcal{M}_0$ is obtained by replacing $p$ with $0$ in \eqref{eroded_free_space_definition}. In other words, if $\mathbf{x}(0)\in\mathcal{M}_0$, then $\mathbf{x}(t)\in\mathcal{M}_0$ for all $t\geq 0$.
        \item \label{claim:conic2}The target point $\mathbf{x}_d$ is almost globally asymptotically stable over $\mathcal{M}_0$.
    \end{enumerate}
\end{theorem}
\proof{See Appendix \ref{proof:maintheorem:conic_constraints}.}

Since $\epsilon < 1 - \cos(\delta)$,  it follows from Lemma \ref{lemma:constraint_separation} that the $\epsilon$-dilated constraints $\mathcal{D}_{\epsilon}(\mathcal{U}_i), i\in\mathbb{I},$ are pairwise disjoint \textit{i.e.}, $\mathcal{D}_{\epsilon}(\mathcal{U}_i)\cap\mathcal{D}_{\epsilon}(\mathcal{U}_j) = \emptyset$ for all  $i, j\in\mathbb{I}, i\ne j$.
Consequently, for any $\mathbf{x}\in\mathcal{N}_{\epsilon}(\mathcal{U})$, the control input \eqref{negative_gradient_control_law} enforces safety with respect to a unique closest unsafe region and admits the following representation:
\begin{equation}\label{conic_control_expression}
    \mathbf{u}(\mathbf{x}) = \begin{cases}
    \frac{k_1}{(h_i(d_s(\mathbf{x}, \mathbf{g}_i))+d_s(\mathbf{x}, \mathbf{x}_d))^2}\mathbf{u}_i^c(\mathbf{x}), &\mathbf{x}\in\mathcal{N}_{\epsilon}(\mathcal{U}_i),\\
    \frac{k_1}{(1 + d_s(\mathbf{x}, \mathbf{x}_d))^2}\mathbf{x}_d, & \mathbf{x}\notin\mathcal{N}_{\epsilon}(\mathcal{U}),
    \end{cases}
\end{equation}
where for each $i\in\mathbb{I}$, $\mathbf{u}_i^c(\mathbf{x})$ is given by
\begin{equation}\label{conic_individual_control_expression}
    \mathbf{u}_i^c(\mathbf{x}) = h_i(d_s(\mathbf{x}, \mathbf{g}_i))\mathbf{x}_d -h_i'(d_s(\mathbf{x}, \mathbf{g}_i))d_s(\mathbf{x}, \mathbf{x}_d)\mathbf{g}_i.
\end{equation}

It follows from \eqref{conic_control_expression} and \eqref{conic_individual_control_expression} that
for any $\mathbf{x}\in\mathcal{M}_0$, the control input \eqref{negative_gradient_control_law} can be expressed as a linear combination of at most two unit vectors, namely $\mathbf{x}_d$ and $-\mathbf{g}_i$ for some $i\in\mathbb{I}$, with non-negative coefficients.
As $\mathbf{x}$ approaches the boundary $\partial\mathcal{U}_i$ for some $i\in\mathbb{I}$, the control input \eqref{negative_gradient_control_law} tends to $-\frac{k_1h_i'(d_s(\mathbf{x}, \mathbf{g}_i))}{d_s(\mathbf{x}, \mathbf{x}_d)}\mathbf{g}_i$,
and steers $\mathbf{x}$ in a direction aligned with the geodesic $\mathcal{G}(\mathbf{x}, -\mathbf{g}_i)$ toward $-\mathbf{g}_i$.
Consequently, for each $i\in\mathbb{I}$, since $\mathcal{G}(\mathbf{x}, -\mathbf{g}_i)\cap\mathcal{U}_i^{\circ} = \emptyset$ for all $\mathbf{x}\in\mathcal{N}_{\epsilon}(\mathcal{U}_i)$, the control law \eqref{negative_gradient_control_law} ensures the forward invariance of $\mathcal{M}_0$ for the closed-loop system \eqref{system_dynamics}-\eqref{negative_gradient_control_law}.


\begin{figure}[ht]
\centering\includegraphics[width =0.47\linewidth]{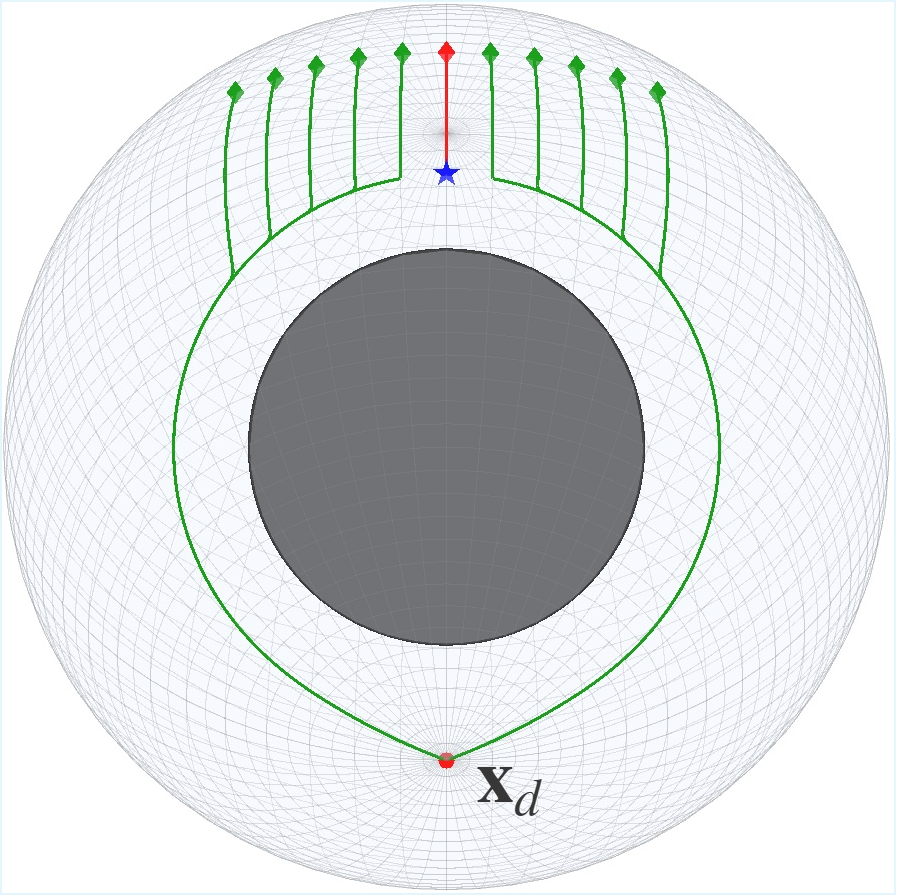}
\caption{$\mathbf{x}$-trajectories under the gradient-based control law~\eqref{negative_gradient_control_law} in the presence of a conic constraint.}\label{diagram_motivation_for_star_shaped}
\end{figure}

Due to the constraint $\mathcal{U}$, in general no continuous time-invariant feedback law, including the negative gradient control law \eqref{negative_gradient_control_law}, can render the target point $\mathbf{x}_d$ globally asymptotically stable over $\mathcal{M}_0$.
Under the assumption that $\mathcal{U}$ is the union of pairwise non-overlapping conic constraints $\mathcal{U}_i$, $i\in\mathbb{I}$, as stated in Assumption~\ref{assumption:non-overlapping-constraints}, the closed-loop system \eqref{system_dynamics}–\eqref{negative_gradient_control_law} admits $\mathbf{x}_d$ as an almost globally asymptotically stable equilibrium on $\mathcal{M}_0$, as stated in Theorem \ref{maintheorem:conic_constraints}. 
In particular, the set of initial conditions in $\mathcal{M}_0$ from which convergence to $\mathbf{x}_d$ fails has zero Lebesgue measure on $\mathbb{S}^n$.
This behaviour is illustrated in Fig.~\ref{diagram_motivation_for_star_shaped}, where an undesired unstable equilibrium is indicated by a star symbol. 
Except for the trajectory initialized at the red diamond, all other trajectories safely converge to $\mathbf{x}_d$. 

The gradient-based controller developed above does not directly extend to general star-shaped constraint sets on $\mathbb{S}^n$.
In particular, unlike the case of conic constraints, for a general star-shaped constraint set $\mathcal{U}_i$ on $\mathbb{S}^n$, there need not exist $\mathbf{g}_i\in\mathcal{U}_i^{\circ}$ such that the angle between $\mathbf{g}_i$ and every $\mathbf{x}\in\partial\mathcal{U}_i$ is constant.
Hence, $\partial\mathcal{U}_i$ cannot generally be represented by a single level set of $d_s(\mathbf{x}, \mathbf{g}_i)$, as required for $h_i(d_s(\mathbf{x}, \mathbf{g}_i))$ in \eqref{obstacle-function-defnition} to vanish on the entire constraint boundary.
This observation motivates the need for the design of a feedback control law that guarantees safe and almost global asymptotic stabilization of the target point $\mathbf{x}_d$ on $\mathcal{M}_0$ in the presence of general star-shaped constraint sets on $\mathbb{S}^n$.

\section{Constrained stabilization under star-shaped constraints}\label{section:feedback_control_design}

Let $\mathcal{U}_i$ denote the star-shaped set on $\mathbb{S}^n$ for each $i\in\mathbb{I}$, where a star-shaped set on $\mathbb{S}^n$ is defined in Section \ref{notations}.
According to Lemma \ref{lemma:reverse_geodesic_always_stay_outside}, for each $i\in\mathbb{I}$, if one selects $\mathbf{g}_i \in \sigma(\mathcal{U}_i)\cap\mathcal{U}_i^{\circ}$, then for every $\mathbf{x}\in\partial\mathcal{U}_i$, the geodesic $\mathcal{G}(\mathbf{x},-\mathbf{g}_i)$ connecting $\mathbf{x}$ to $-\mathbf{g}_i$ does not intersect the interior of $\mathcal{U}_i$, where $\sigma(\mathcal{U}_i)$ is defined as in \eqref{sigma_set}.  
Consequently, we use a strategy analogous to the one used for conic constraint avoidance in Section \ref{section:conic_constraint_avoidance}, namely, steering $\mathbf{x}$ along the geodesic $\mathcal{G}(\mathbf{x},-\mathbf{g}_i)$ toward $-\mathbf{g}_i$ whenever $\mathbf{x}$ approaches the boundary $\partial\mathcal{U}_i$ for some 
$i\in\mathbb{I}$.
This ensures that $\mathbf{x}$ remains outside $\mathcal{U}_i^{\circ}$ and guarantees the forward 
invariance of $\mathcal{M}_0$.
The proposed feedback control $\mathbf{u}(\mathbf{x})$ is given by
\begin{equation}\label{n-sphere-control-law}
    \mathbf{u}(\mathbf{x}) = \begin{cases}
    k_1\mathbf{u}_i(\mathbf{x}), &\mathbf{x}\in\mathcal{N}_{\epsilon}(\mathcal{U}_i),\\
    k_1\mathbf{x}_d, & \mathbf{x}\notin\mathcal{N}_{\epsilon}(\mathcal{U}),
    \end{cases}
\end{equation}
where $k_1 > 0$.
Similar to Section \ref{section:conic_constraint_avoidance}, the parameter $\epsilon$ is chosen such that $\epsilon\in(0, \min\{\bar{\epsilon}, 1- \cos(\delta)\})$, where $\bar{\epsilon}$ is a strictly positive scalar such that $\mathbf{x}_d\notin\mathcal{N}_{\bar{\epsilon}}(\mathcal{U})$.
Since $\mathbf{x}_d\in\mathcal{M}_0^{\circ}$, one has $\mathbf{x}_d\notin\mathcal{U}$, and the existence of $\bar{\epsilon} > 0$ such that $\mathbf{x}_d\notin\mathcal{N}_{\bar{\epsilon}}(\mathcal{U})$ is straightforward to establish.
Selecting $\epsilon < 1-\cos(\delta)$ ensures that the sets $\mathcal{D}_{\epsilon}(\mathcal{U}_i)$, $i\in\mathbb{I}$, are pairwise disjoint, as established earlier in Lemma \ref{lemma:constraint_separation}.
The vector-valued function $\mathbf{u}_{i}(\mathbf{x})$ is given by
\begin{equation}\label{individual_control_input_vector_design}
    \mathbf{u}_i(\mathbf{x}) = 
    \frac{d_s(\mathbf{x}, \mathcal{U}_i)}{\epsilon}\mathbf{x}_d - \frac{1}{\kappa}\left(1 - \frac{d_s(\mathbf{x}, \mathcal{U}_i)}{\epsilon}\right)\mathbf{g}_i,
\end{equation}
where $\kappa > 0$.
For each $i\in\mathbb{I}$, the constant unit vector $\mathbf{g}_i$ is chosen such that $\mathbf{g}_i\in\sigma(\mathcal{U}_i)\cap\mathcal{U}_i^{\circ}$ and $\mathbf{g}_i\ne-\mathbf{x}_d$, where the set $\sigma(\mathcal{U}_i)$ is defined according to \eqref{sigma_set} and $\mathcal{U}_i^{\circ}$ denotes the interior of $\mathcal{U}_i$ on $\mathbb{S}^n$.\footnote{Selecting $\mathbf{g}_i\in\sigma(\mathcal{U}_i)\cap\mathcal{U}_i^{\circ}$ allows us to leverage Lemma \ref{lemma:reverse_geodesic_always_stay_outside} to establish the forward invariance of $\mathcal{M}_0$ for the closed-loop system \eqref{system_dynamics}-\eqref{n-sphere-control-law}, as discussed later in Lemma \ref{lemma:n-sphere-forward-invariance}.
Furthermore, ensuring $\mathbf{g}_i\ne-\mathbf{x}_d$ for every $i\in\mathbb{I}$ guarantees that the geodesics $\mathcal{G}(\mathbf{x}_d, \mathbf{g}_i)$ and $\mathcal{G}(-\mathbf{x}_d, -\mathbf{g}_i)$, which are used later in Section \ref{section:stability_analysis}, are well-defined.}

\begin{remark}
The parameter $\kappa > 0$ scales the magnitude of the repulsive vector field in the control law \eqref{n-sphere-control-law}.
Unlike the case of conic constraints, for a general star-shaped constraint set $\mathcal{U}_i$ on $\mathbb{S}^n$, there may exists $\mathbf{g}_i\in\sigma(\mathcal{U}_i)\cap\mathcal{U}_i^{\circ}$ such that $-\mathbf{g}_i\in\mathcal{N}_{\epsilon}(\mathcal{U}_i)\setminus\partial\mathcal{U}_i$.
In this situation, the set $\mathcal{G}(-\mathbf{g}_i, \mathbf{x}_d)\cap\mathcal{N}_{\epsilon}(\mathcal{U}_i)$ is non-empty, and the vector field $\mathbf{P}(\mathbf{x})\mathbf{u}(\mathbf{x})$ may vanish at certain points in this set, potentially introducing undesired equilibria.
To avoid this behaviour, later in Theorem \ref{main_theorem}, we show the existence of $\bar{\kappa} > 0$ such that, for all $\kappa > \bar{\kappa}$, the contribution of the repulsive vector field in \eqref{individual_control_input_vector_design} is sufficiently small to ensure $\dot{d}_s(\mathbf{x}, \mathbf{x}_d) < 0$ for every $\mathbf{x}\in\mathcal{G}(-\mathbf{g}_i, \mathbf{x}_d)\cap\mathcal{N}_{\epsilon}(\mathcal{U}_i)$, and hence $\mathbf{P}(\mathbf{x})\mathbf{u}(\mathbf{x})\ne \mathbf{0}$.
\end{remark}

\begin{remark}[Continuous control input]
Since $\epsilon < 1 - \cos(\delta)$, it follows from Lemma \ref{lemma:constraint_separation} that $\mathcal{N}_{\epsilon}(\mathcal{U}_i)\cap\mathcal{N}_{\epsilon}(\mathcal{U}_j) = \emptyset$ for all $i, j\in\mathbb{I}$ with $i\ne j$.
Consequently, if $\mathbf{x}\in\mathcal{N}_{\epsilon}(\mathcal{U})$, then there exists a unique index $i\in\mathbb{I}$ such that $\mathbf{x}\in\mathcal{N}_{\epsilon}(\mathcal{U}_i)$ and $\mathbf{u}(\mathbf{x}) = k_1\mathbf{u}_i(\mathbf{x})$.
Furthermore, $\mathbf{u}_i(\mathbf{x})$ is continuous for each $i\in\mathbb{I}$ and for all $\mathbf{x}\in\mathcal{N}_{\epsilon}(\mathcal{U}_i)$.
Moreover, for each $i\in\mathbb{I}$ and for every $\mathbf{x}\in\partial\mathcal{N}_{\epsilon}(\mathcal{U}_i)\cap\mathcal{M}_{\epsilon}$, $\mathbf{u}_i(\mathbf{x})$ simplifies to $\mathbf{u}_i(\mathbf{x}) = \mathbf{x}_d$, where the set $\mathcal{M}_{\epsilon}$ is obtained by replacing $p$ with $\epsilon$ in \eqref{eroded_free_space_definition}.
As a result, the proposed feedback control input $\mathbf{u}(\mathbf{x})$, defined in \eqref{n-sphere-control-law}, is continuous for all $\mathbf{x}\in\mathcal{M}_0$.
\end{remark}

Similar to the negative gradient-based control law \eqref{negative_gradient_control_law}, the feedback control law \eqref{n-sphere-control-law} steers $\mathbf{x}$ along the $\mathcal{G}(\mathbf{x}, \mathbf{x}_d)$ toward $\mathbf{x}_d$ whenever $d_s(\mathbf{x}, \mathcal{U}) \geq \epsilon$.
Within the $\epsilon$-neighborhood $\mathcal{N}_{\epsilon}(\mathcal{U}_i)$ for some $i\in\mathbb{I}$, the control input is a linear combination of the unit vectors $\mathbf{x}_d$ and $-\mathbf{g}_i$ with non-negative coefficients.
Moreover, as $\mathbf{x}$ approaches the boundary $\partial\mathcal{U}_i$, the control input \eqref{n-sphere-control-law} tends to $-\frac{k_1}{\kappa}\mathbf{g}_i$ and, similar to the negative gradient-based control law \eqref{negative_gradient_control_law}, steers $\mathbf{x}$ along the geodesic $\mathcal{G}(\mathbf{x},-\mathbf{g}_i)$ toward $-\mathbf{g}_i$.

\subsection{Safety and stability analysis}\label{section:stability_analysis}

\begin{figure}[ht]
    \centering
    \includegraphics[width=0.6\linewidth]{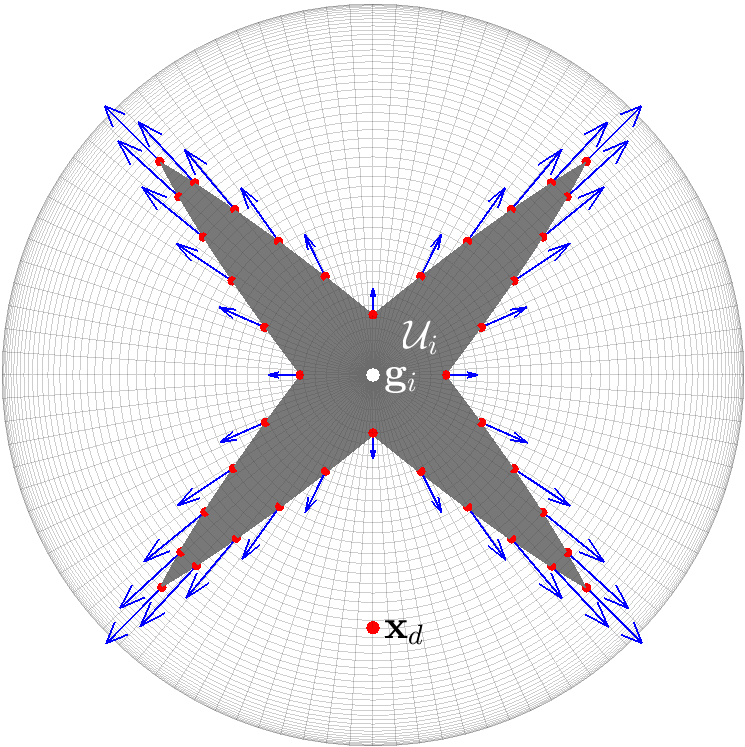}
    \caption{Representation of $-\mathbf{P}(\mathbf{x})\mathbf{g}_i$ for $\mathbf{x}\in\partial\mathcal{U}_i$.}
    \label{fig:forward_invaraince}
\end{figure}

First, we analyze the forward invariance of the safe region $\mathcal{M}_0$ for the closed-loop system \eqref{system_dynamics}-\eqref{n-sphere-control-law}.
According to Assumption \ref{assumption:non-overlapping-constraints}, if $\mathbf{x}\in\partial\mathcal{M}_0$, then $\mathbf{x}\in\partial\mathcal{U}_i$ for some $i\in\mathbb{I}$ and $\mathbf{x}\notin\partial\mathcal{U}_j$ for all $j\in\mathbb{I}$ with $j\ne i$.
According to \eqref{individual_control_input_vector_design}, if $\mathbf{x}\in\partial\mathcal{U}_i$ for some $i\in\mathbb{I}$, then the control input vector \eqref{n-sphere-control-law} simplifies to 
\begin{equation}\label{control_at_boundary}\mathbf{u}(\mathbf{x}) = \frac{-k_1}{\kappa}\mathbf{g}_i,\end{equation}
and steers $\mathbf{x}$ along the geodesic $\mathcal{G}(\mathbf{x}, -\mathbf{g}_i)$ toward $-\mathbf{g}_i$. 
Additionally, since $\mathcal{U}_i$ is a star-shaped constraint on the $n$-sphere and $\mathbf{g}_i \in \sigma(\mathcal{U}_i)\cap\mathcal{U}_i^{\circ}$, Lemma \ref{lemma:reverse_geodesic_always_stay_outside} implies that $\mathcal{G}(\mathbf{x}, -\mathbf{g}_i)\cap\mathcal{U}_i^{\circ} = \emptyset$. Consequently, when $\mathbf{x}\in\partial\mathcal{U}_i$, the vector $\mathbf{P}(\mathbf{x})\mathbf{u}(\mathbf{x})$, with $\mathbf{u}(\mathbf{x})$ given in \eqref{control_at_boundary} does not point to the interior of the unsafe region $\mathcal{U}_i$, as illustrated in Fig. \ref{fig:forward_invaraince}.
This behaviour allows us to establish the forward invariance of the set $\mathcal{M}_0$ for the closed-loop system \eqref{system_dynamics}-\eqref{n-sphere-control-law}, as stated in the next lemma.
\begin{lemma}\label{lemma:n-sphere-forward-invariance}
    For the closed-loop system \eqref{system_dynamics}-\eqref{n-sphere-control-law} under Assumption \ref{assumption:non-overlapping-constraints}, the set $\mathcal{M}_0$, defined according to \eqref{eroded_free_space_definition}, is forward invariant.
    In other words, if $\mathbf{x}(0)\in\mathcal{M}_0$, then $\mathbf{x}(t)\in\mathcal{M}_0$ for all $t\geq 0$.
\end{lemma}
\proof{See Appendix \ref{proof:lemma:n-sphere-forward-invariance}.}

Next, we show that $\mathbf{P}(\mathbf{x})\mathbf{u}(\mathbf{x})$ is locally Lipschitz over $\mathcal{M}_0$.
Combined with Lemma \ref{lemma:n-sphere-forward-invariance}, this will ensure that the solution to the closed-loop system \eqref{system_dynamics}-\eqref{n-sphere-control-law} is uniquely defined for each initial condition $\mathbf{x}(0)\in\mathcal{M}_0$ and exists for all $t \geq 0$.
\begin{lemma}\label{lemma:lipschitz_continuity}
    The continuous vector-valued function $\mathbf{P}(\mathbf{x})\mathbf{u}(\mathbf{x})$ is locally Lipschitz over $\mathcal{M}_0$.
\end{lemma}
\proof{See Appendix \ref{proof:lemma:lipschitz_continuity}.}

Next, we analyze the convergence properties of the proposed closed-loop system \eqref{system_dynamics}-\eqref{n-sphere-control-law}.
When $\mathbf{x}\in\mathcal{N}_{\epsilon}(\mathcal{U}_i)\setminus\{-\mathbf{x}_d\}$ for some $i\in\mathbb{I}$, the repulsive component $-\frac{k_1}{\kappa}\left(1 - \frac{d_s(\mathbf{x}, \mathcal{U}_i)}{\epsilon}\right)\mathbf{g}_i$ of the control input steers $\mathbf{x}$ along the geodesic $\mathcal{G}(\mathbf{x}, -\mathbf{g}_i)$ toward $-\mathbf{g}_i$. 
Meanwhile, the attractive component $k_1\frac{d_s(\mathbf{x}, \mathcal{U}_i)}{\epsilon}\mathbf{x}_d$ steers $\mathbf{x}$ along the geodesic $\mathcal{G}(\mathbf{x}, \mathbf{x}_d)$ toward $\mathbf{x}_d$.
This interaction leads to an increase in the cosine of the angle between the vectors $\mathbf{P}(\mathbf{g}_i)(\mathbf{x}-\mathbf{g}_i)$ and $\mathbf{P}(\mathbf{g}_i)(\mathbf{x}_d - \mathbf{g}_i)$ as long as $\mathbf{x}\in\mathcal{N}_{\epsilon}(\mathcal{U}_i)\setminus\left(\partial\mathcal{U}_i\cup\mathcal{Z}_i\cup\mathcal{V}_i\right)$, as established in the next lemma, where 
for each $i\in\mathbb{I}$, the set $\mathcal{Z}_i$ and $\mathcal{V}_i$ are defined as
\begin{equation}\label{set_definition_Vi_Zi}
\begin{aligned}
    \mathcal{Z}_i &= \mathcal{G}(\mathbf{g}_i, -\mathbf{x}_d)\cup\mathcal{G}(-\mathbf{g}_i, -\mathbf{x}_d),\\
    \mathcal{V}_i &= \mathcal{G}(\mathbf{g}_i, \mathbf{x}_d)\cup\mathcal{G}(-\mathbf{g}_i, \mathbf{x}_d).
    \end{aligned}
\end{equation}

\begin{lemma}\label{lemma:angle_keeps_changing}
Consider the closed-loop system \eqref{system_dynamics}–\eqref{n-sphere-control-law} under Assumption \ref{assumption:non-overlapping-constraints}. For each $i \in \mathbb{I}$, define the scalar function
\begin{equation}\label{the_only_positive_semidefinite_function}
    V_i(\mathbf{x}) = \left(\frac{\mathbf{P}(\mathbf{g}_i)(\mathbf{x}_d- \mathbf{g}_i)}{\|\mathbf{P}(\mathbf{g}_i)(\mathbf{x}_d- \mathbf{g}_i)\|}\right)^\top\left(\frac{\mathbf{P}(\mathbf{g}_i)(\mathbf{x}- \mathbf{g}_i)}{\|\mathbf{P}(\mathbf{g}_i)(\mathbf{x}- \mathbf{g}_i)\|}\right),
\end{equation}
over $\mathcal{F}_i$, where $\mathcal{F}_i = \left(\mathcal{N}_{\epsilon}(\mathcal{U}_i) \cup \mathcal{M}_{\epsilon}\right)\setminus\{-\mathbf{g}_i\}$. 
Then:
\begin{enumerate}
    \item\label{claim1:LemmaVi} $V_i(\mathbf{x})$ is well-defined for all $\mathbf{x}\in\mathcal{F}_i$,
    \item\label{claim2:LemmaVi} \label{scalar:claim3}$\dot{V}_i(\mathbf{x}) > 0$ for all $\mathbf{x} \in \mathcal{F}_i\setminus\left(\partial\mathcal{U}_i \cup \mathcal{Z}_i \cup \mathcal{V}_i\right)$,
    \item\label{claim3:lemmaVi}$\dot{V}_i(\mathbf{x}) = 0 $ for all $\mathbf{x}\in\mathcal{F}_i\cap\left(\partial\mathcal{U}_i \cup \mathcal{Z}_i \cup \mathcal{V}_i\right)$,
\end{enumerate}
where the sets $\mathcal{Z}_i$ and $\mathcal{V}_i$ are defined in \eqref{set_definition_Vi_Zi}.
\end{lemma}
\proof{See Appendix \ref{proof:lemma:angle_keeps_changing}.}

\begin{remark}\label{remark:possibilities}
For $i \in \mathbb{I}$ and $\mathbf{x} \in \mathcal{F}_i$, the function $V_i(\mathbf{x})$, defined in \eqref{the_only_positive_semidefinite_function}, represents the cosine of the angle between the projected vectors $\mathbf{P}(\mathbf{g}_i)(\mathbf{x} - \mathbf{g}_i)$ and $\mathbf{P}(\mathbf{g}_i)(\mathbf{x}_d - \mathbf{g}_i)$.
In $\mathcal{F}_i$, it attains its minimum value of $-1$ if and only if $\mathbf{x}\in\mathcal{Z}_i\cap\mathcal{F}_i$, and its maximum value of $1$ if and only if $\mathbf{x}\in\mathcal{V}_i\cap\mathcal{F}_i$.
Moreover, if there exists $t_1\geq 0$ such that $\mathbf{x}(t_1)\in\partial\mathcal{U}_i$, then the control input \eqref{n-sphere-control-law} becomes $\mathbf{u}(\mathbf{x}(t_1)) = -\frac{k_1}{\kappa}\mathbf{g}_i$, and it steers $\mathbf{x}$ in the direction aligned with the geodesic $\mathcal{G}(\mathbf{x}(t_1), -\mathbf{g}_i)$ towards $-\mathbf{g}_i$.
Furthermore, since $\mathcal{U}_i$ is a star-shaped set on $\mathbb{S}^n$ and $\mathbf{g}_i\in\sigma(\mathcal{U}_i)\cap\mathcal{U}_i^{\circ}$, one has $-\mathbf{g}_i\notin\mathcal{U}_i$, and by virtue of Lemma \ref{lemma:reverse_geodesic_always_stay_outside},
the geodesic $\mathcal{G}(\mathbf{x}, -\mathbf{g}_i)$
does not intersect with $\mathcal{U}_i^{\circ}$ for all $\mathbf{x}\in\partial\mathcal{U}_i$.
Consequently, since $\mathbf{u}(\mathbf{x}(t)) = -\frac{k_1}{\kappa}\mathbf{g}_i$ for all $t\geq t_1$ as long as $\mathbf{x}(t)\in\partial\mathcal{U}_i$, $\mathbf{x}(t)$ cannot remain on $\partial\mathcal{U}_i$ indefinitely and enters the interior of the free space $\mathcal{M}_0^{\circ}$ in finite time, \textit{i.e.}, there exists $t_2 > t_1$ such that $\mathbf{x}(t_2)\in\mathcal{M}_0^{\circ}$ and $\mathbf{x}(t)\in\partial\mathcal{U}_i$ for all $t\in[t_1, t_2)$.
By Claim \ref{claim3:lemmaVi} of Lemma \ref{lemma:angle_keeps_changing}, the value of $V_i(\mathbf{x}(t))$ does not change for all $t\in[t_1, t_2)$.
Consequently, it follows from Claim \ref{scalar:claim3} of Lemma \ref{lemma:angle_keeps_changing} that if there exists $t_1\geq 0$ such that $\mathbf{x}(t_1)\in\mathcal{F}_i\setminus\mathcal{Z}_i$, then one of the following statements hold:
\begin{enumerate}
    \item There exists $s_1 > t_1$ such that $\mathbf{x}(s_1)\in\mathcal{M}_0\setminus\mathcal{F}_i$.
    \item  $\Lim_{t\to\infty}d_s(\mathbf{x}(t), \mathcal{V}_i) = 0$ and $\mathbf{x}(t)\in\mathcal{F}_i\setminus(\mathcal{V}_i\cup\mathcal{Z}_i)$ for all $t\geq t_1$.
    \item $\mathbf{x}(t)\in\mathcal{F}_i\cap\mathcal{V}_i$ for all $t\geq t_1$.
\end{enumerate}
This behaviour of a solution $\mathbf{x}(t)$ helps us in establishing the almost global asymptotic stability of $\mathbf{x}_d$ for the closed-loop system \eqref{system_dynamics}-\eqref{n-sphere-control-law} over $\mathcal{M}_0$, as stated later in Theorem \ref{main_theorem}.
\end{remark}

\begin{remark}
    The scalar function $V_i(\mathbf{x})$, defined in \eqref{the_only_positive_semidefinite_function}, is undefined at $\mathbf{x} = \mathbf{g}_i$ and $\mathbf{x} = -\mathbf{g}_i$ since $\mathbf{P}(\mathbf{g}_i)\mathbf{g}_i = \mathbf{0}$.
    However, $\mathbf{g}_i\in\mathcal{U}_i^{\circ}$ for each $i\in\mathbb{I}$, and therefore $\mathbf{g}_i\notin\bar{\mathcal{F}}_i$ for any $i\in\mathbb{I}$.
    On the other hand, since $\mathcal{U}_i$ is a star-shaped set on $\mathbb{S}^n$ and $\mathbf{g}_i\in\sigma(\mathcal{U}_i)$, by \eqref{sigma_set}, one has $-\mathbf{g}_i\notin\mathcal{U}_i$ for each $i\in\mathbb{I}$.
    Therefore, it is possible that $-\mathbf{g}_i\in\bar{\mathcal{F}}_i$ for some $i\in\mathbb{I}$.
    In this case, by \eqref{system_dynamics} and \eqref{n-sphere-control-law}, the closed-loop vector field at $\mathbf{x} = -\mathbf{g}_i$ satisfies 
    \[\mathbf{P}(-\mathbf{g}_i)\mathbf{u}(-\mathbf{g}_i) = k_1b_i(-\mathbf{g}_i)\mathbf{P}(-\mathbf{g}_i)\mathbf{x}_d,\] where $b_i(-\mathbf{g}_i) = \min\left\{\frac{d_s(-\mathbf{g}_i, \mathcal{U}_i)}{\epsilon}, 1\right\}$.
    Since $\mathbf{g}_i\ne\mathbf{x}_d$ and $\mathbf{g}_i\ne-\mathbf{x}_d$, the control input $\mathbf{P}(-\mathbf{g}_i)\mathbf{u}(-\mathbf{g}_i)$ is non-zero.
    
    Consequently, when $-\mathbf{g}_i\in\bar{\mathcal{F}}_i$ for some $i\in\mathbb{I}$, and a solution $\mathbf{x}(t)$ satisfies $\mathbf{x}(t_1) = -\mathbf{g}_i$ for some $t_1\geq 0$, then immediately after $t_1$, the control input steers $\mathbf{x}$ in the direction aligned with the geodesic $\mathcal{G}(\mathbf{x}, \mathbf{x}_d)$ towards $\mathbf{x}_d$, and $\mathbf{x}(t)$ enters the set $\mathcal{F}_i$ \textit{i.e.} there exists $\gamma > 0$ such that $\mathbf{x}(t)\in\mathcal{F}_i$ for all $t\in(t_1, t_1 + \gamma]$.
    Loosely speaking, even though $V_i$ is undefined at some isolated point $-\mathbf{g}_i$ in the set $\bar{\mathcal{F}}_i$, any solution $\mathbf{x}(t)$ that reaches $-\mathbf{g}_i$ immediately enters $\mathcal{F}_i$ and is steered towards the desired point $\mathbf{x}_d$. 
\end{remark}

According to Lemma \ref{lemma:angle_keeps_changing} and Remark \ref{remark:possibilities}, if there exists $t_1\geq 0$ such that $\mathbf{x}(t_1)\in\mathcal{F}_i\setminus\left(\mathcal{Z}_i\cup\mathcal{V}_i\right)$, then the control input vector \eqref{n-sphere-control-law} 
drives $\mathbf{x}$ away from $\mathcal{Z}_i\cap\mathcal{F}_i$ and toward $\mathcal{V}_i\cap\mathcal{F}_i$ for all $t\geq t_1$ as long as $\mathbf{x}(t)\in\mathcal{F}_i\setminus\left(\mathcal{Z}_i\cup\mathcal{V}_i\right)$. 
It is possible that $\mathbf{x}(t)$ exits $\mathcal{N}_{\epsilon}(\mathcal{U}_i)$ and enters the set $\mathcal{M}_\epsilon^{\circ}$.
In such a case, the trajectory $\mathbf{x}(t)$ may continue toward another neighborhood $\mathcal{N}_{\epsilon}(\mathcal{U}_j)$, with $j \in \mathbb{I} \setminus \{i\}$, where the new entry point $\mathbf{h}_j$ to $\mathcal{N}_{\epsilon}(\mathcal{U}_j)$ is farther from $\mathbf{x}_d$ than the previous entry point $\mathbf{h}_i$ to  $\mathcal{N}_{\epsilon}(\mathcal{U}_i)$ \textit{i.e.}, $d_s(\mathbf{x}_d, \mathbf{h}_j) > d_s(\mathbf{x}_d, \mathbf{h}_i)$.
This behaviour introduces the possibility of closed trajectories, which prevents us from establishing almost global asymptotic convergence to the desired point $\mathbf{x}_d$ for the closed-loop system \eqref{system_dynamics}-\eqref{n-sphere-control-law}. 
To avoid such cases, we require that the unsafe regions $\mathcal{U}_i$, where $i \in \mathbb{I}$, be sufficiently separated, as described next.

Let $\mathbb{I}_a$ be a subset of $\mathbb{I}$ such that for every $i\in\mathbb{I}_a$, $-\mathbf{x}_d\notin\mathcal{D}_{\epsilon}(\mathcal{U}_i)$, as defined below
\begin{equation}\label{definition:Ia}
    \mathbb{I}_a = \{i\in\mathbb{I}\mid-\mathbf{x}_d\notin\mathcal{D}_{\epsilon}(\mathcal{U}_i)\}.
\end{equation}
The set $\mathbb{I}\setminus\mathbb{I}_a$ is either a singleton set or an empty set.
For each $i\in\mathbb{I}_a$, the set $\mathcal{S}_i(\mathbf{x}_d)$ is the union of all geodesics $\mathcal{G}(\mathbf{x}, \mathbf{x}_d)$ with $\mathbf{x}\in\mathcal{D}_{\epsilon}(\mathcal{U}_i)$, defined as follows:
\begin{equation}\label{definition:Si}
    \mathcal{S}_i(\mathbf{x}_d) = \bigcup_{\mathbf{x}\in\mathcal{D}_{\epsilon}(\mathcal{U}_i)}\mathcal{G}(\mathbf{x}, \mathbf{x}_d).
\end{equation}

\begin{figure}
    \centering
    \includegraphics[width=0.7\linewidth]{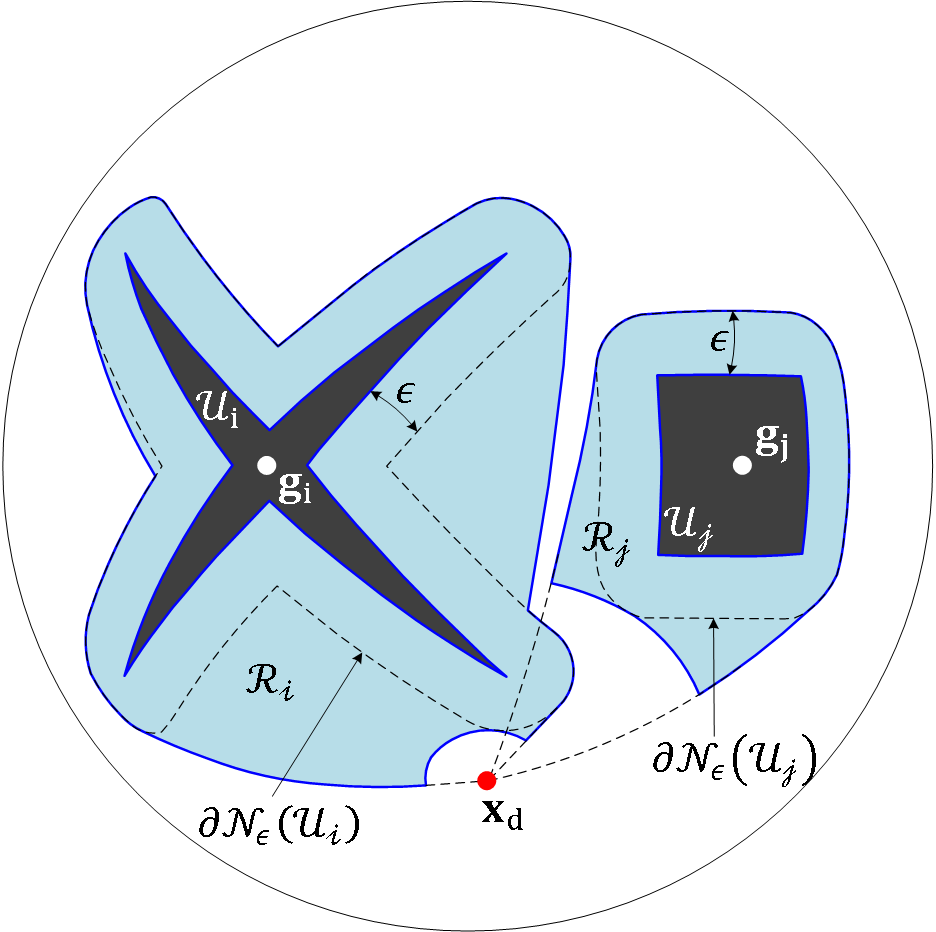}
    \caption{Illustration of mutually exclusive sets $\mathcal{R}_i$, where $i\in\mathbb{I}_a$.}
    \label{fig:mutually_exclusive_Ri_sets}
\end{figure}

\begin{figure}
    \centering
    \includegraphics[width=0.7\linewidth]{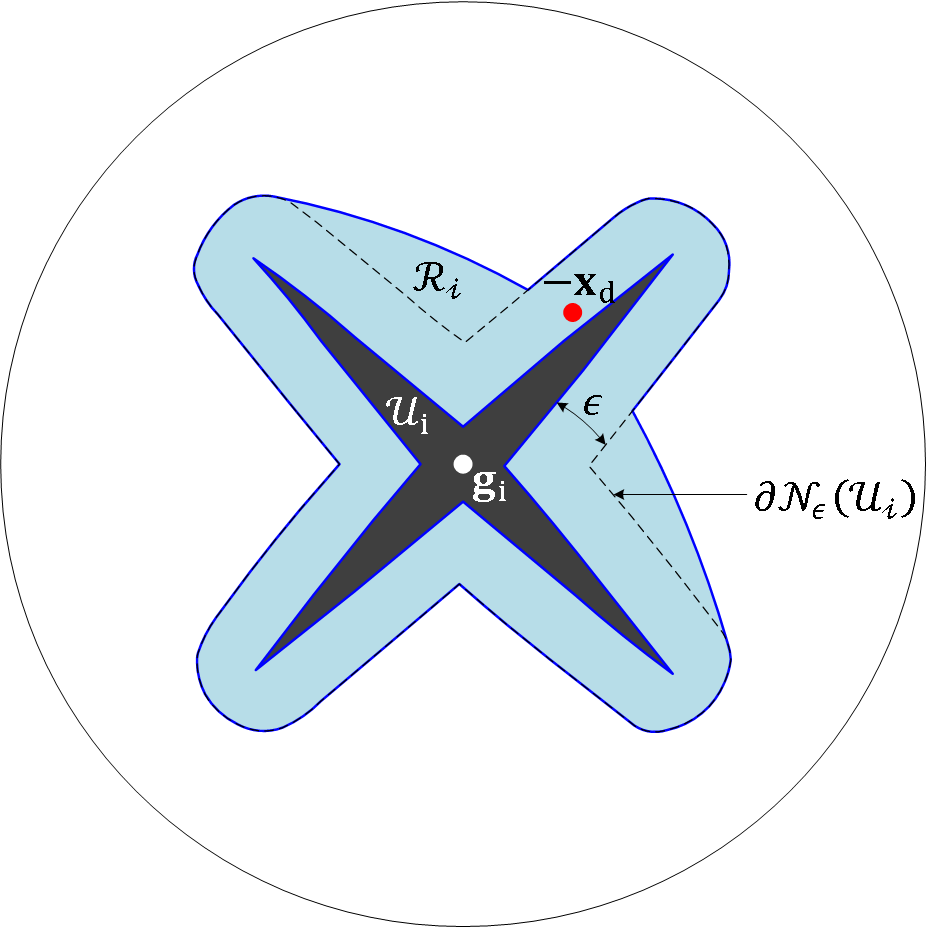}
    \caption{Illustration of the set $\mathcal{R}_i$, where $i\in\mathbb{I}\setminus\mathbb{I}_a$.}
    \label{fig:the_separate_Ri}
\end{figure}

For each $i\in\mathbb{I}_a$, the region $\mathcal{R}_i$ is defined as
\begin{equation}\label{definition:Ri}
    \mathcal{R}_i = \{\mathbf{x}\in\mathcal{S}_i(\mathbf{x}_d)\setminus\mathcal{U}_i^{\circ}\mid d_s(\mathbf{x}, \mathbf{x}_d)\geq d_s(\mathbf{x}_d, \mathcal{D}_{\epsilon}(\mathcal{U}_i))\},
\end{equation}
as illustrated in Fig. \ref{fig:mutually_exclusive_Ri_sets}.
Moreover, if $i\in\mathbb{I}\setminus\mathbb{I}_a$, then define $\mathcal{R}_i$ as follows:
\begin{equation}\label{definition:Ri_forspecialindex}\mathcal{R}_i = \mathcal{S}_i(-\mathbf{x}_d)\setminus\mathcal{U}_i^{\circ},\end{equation} as depicted in Fig. \ref{fig:the_separate_Ri}, where the set $\mathcal{S}_i(-\mathbf{x}_d)$ is obtained using \eqref{definition:Si} by replacing $\mathbf{x}_d$ with $-\mathbf{x}_d$.

For the constraint sets $\mathcal{U}_i$, $i\in\mathbb{I}$, satisfying $-\mathbf{x}_d\notin\mathcal{D}_{\epsilon}(\mathcal{U}_i)$, the set $\mathcal{R}_i$ can be interpreted geometrically as the subset of the free space $\mathcal{M}_0$ consisting of points $\mathbf{x}$ that are lying on geodesics connecting the target point $\mathbf{x}_d$ to points in $\mathcal{D}_{\epsilon}(\mathcal{U}_i)$ and satisfy $d_s(\mathbf{x}_d, \mathbf{x})\geq d_s(\mathbf{x}_d, \mathcal{D}_{\epsilon}(\mathcal{U}_i))$.
If $-\mathbf{x}_d\in\mathcal{D}_{\epsilon}(\mathcal{U}_i)$ for some $i\in\mathbb{I}$, then the corresponding set $\mathcal{R}_i$ is the subset of the free space $\mathcal{M}_0$ consisting of points lying on the geodesics connecting $-\mathbf{x}_d$ to points in $\mathcal{D}_{\epsilon}(\mathcal{U}_i)$.

We require that for each $i, j\in\mathbb{I}$ with $i\ne j$, the sets $\mathcal{R}_i$ and $\mathcal{R}_j$ have no common element, as mentioned in the next assumption.
\begin{assumption}\label{assumption:sufficient_separation}
    The sets $\mathcal{R}_i$ and $\mathcal{R}_j$ are mutually exclusive for all $i, j\in\mathbb{I}$ with $i\ne j$. 
    In other words, for all $i, j\in\mathbb{I}$ with $i\ne j$, $\mathcal{R}_i\cap\mathcal{R}_j = \emptyset$.
\end{assumption}

Assumption \ref{assumption:sufficient_separation} allows us to ensure that if any solution $\mathbf{x}(t)$ to the closed-loop system \eqref{system_dynamics}-\eqref{n-sphere-control-law} is first driven into $\mathcal{R}_i$ at some $\mathbf{h}_i\in\mathcal{R}_i$, where $i\in\mathbb{I}$, and subsequently to $\mathcal{R}_j$ at some $\mathbf{h}_j\in\mathcal{R}_j$, where $j\in\mathbb{I}\setminus\{i\}$, then $d_s(\mathbf{h}_j, \mathbf{x}_d) < d_s(\mathbf{h}_i, \mathbf{x}_d)$.
This behaviour supports the guarantee of almost global asymptotic stability of $\mathbf{x}_d$ for the closed-loop system \eqref{system_dynamics}-\eqref{n-sphere-control-law} over $\mathcal{M}_0$, as stated in the next theorem. 
\begin{theorem}\label{main_theorem}
    For the closed-loop system \eqref{system_dynamics}-\eqref{n-sphere-control-law} under Assumptions \ref{assumption:non-overlapping-constraints} and \ref{assumption:sufficient_separation}, the following statements hold:
    \begin{enumerate}
        \item \label{theorem:claim1}The set $\mathcal{M}_0$ is forward invariant.
        \item \label{theorem:claim2}There exists $\bar{\kappa} > 0$ such that if $\kappa > \bar{\kappa}$, then the desired equilibrium point $\mathbf{x}_d$ is almost globally asymptotically stable over $\mathcal{M}_0$.
    \end{enumerate}
\end{theorem}
\proof{See Appendix \ref{proof:main_theorem}.}

\begin{remark}\label{remark:which_one_is_better_and_why}
    Since a conic constraint $\mathcal{U}_i$, as defined in \eqref{definition:conic_constraints}, is a star-shaped set on $\mathbb{S}^n$ with $\mathbf{g}_i\in\sigma(\mathcal{U}_i)\cap\mathcal{U}_i^{\circ}$, the proposed feedback control law \eqref{n-sphere-control-law} is applicable for the stabilization of $\mathbf{x}_d$ on $\mathbb{S}^n$ in the presence of conic constraints.
    The control law \eqref{n-sphere-control-law} is designed to achieve safe stabilization to $\mathbf{x}_d$ over $\mathcal{M}_0$ in the presence of general star-shaped constraint sets on $\mathbb{S}^n$.
    In this setting, Assumption \ref{assumption:sufficient_separation} provides a sufficient condition under which $\mathbf{x}_d$ is rendered almost globally asymptotically stable over $\mathcal{M}_0$.
    In contrast, when the constraint sets are restricted to conic constraints, the negative gradient-based control law \eqref{negative_gradient_control_law} guarantees safe almost global asymptotic stability of $\mathbf{x}_d$ over $\mathcal{M}_0$ under Assumption \ref{assumption:non-overlapping-constraints}, without requiring the stronger separation condition imposed in Assumption \ref{assumption:sufficient_separation}.
\end{remark}

\section{Application to Constrained Attitude Stabilization}\label{section:application}
We present two illustrative applications of the proposed constrained stabilization framework on $\mathbb{S}^n$.
The first one consists of the constrained stabilization on $\mathbb{S}^2$ (also known as the reduced attitude stabilization) at the kinematic level, where the angular velocity of a rigid body is designed to stabilize the direction of an axis on $\mathbb{S}^2$ to a desired direction while avoiding a predefined set of forbidden directions. 
This problem is related, for instance, to the stabilization of a pointing direction of a satellite \cite{bullo1995control} or the stabilization of a spherical pendulum \cite{berkane2021constrained}. 
The second example consists of the constrained attitude stabilization at the kinematic level, where the attitude is described by the unit-quaternion representation on $\mathbb{S}^3$ and the angular velocity is designed to steer the attitude of the rigid body to a desired attitude while avoiding a set of forbidden attitudes. Although these examples are dealt with at the kinematic level, they could be used to design the torque input at the dynamic level in a cascaded scheme where the designed angular velocity is used as a virtual control to be tracked by the actual control torque. 

\subsection{Reduced attitude control}
Let $\mathbf{x} = \mathbf{R}^\top\mathbf{e}_3\in\mathbb{S}^2$ be the pointing direction of a rigid body, corresponding to the inertial direction $\mathbf{e}_3 = [0, 0, 1]^\top$ in the body frame.
The matrix $\mathbf{R}\in\mathrm{SO}(3)$ represents the orientation of the rigid body attached frame with respect to the inertial frame.
The control objective is to align the pointing direction $\mathbf{x}$ with the desired direction $\mathbf{x}_d\in\mathbb{S}^2$ while avoiding unsafe directions represented as star-shaped constraint sets on $\mathbb{S}^2$.
The rigid body kinematics is given by
\begin{equation}\label{rotation_matrix_dynamics}\dot{\mathbf{R}} = \mathbf{R}\boldsymbol{\omega}^{\times},\end{equation}
where $\boldsymbol{\omega}\in\mathbb{R}^3$ is the angular velocity of the rigid body in the body frame and $\boldsymbol{\omega}^{\times}\in\mathbb{R}^{3\times 3}$ is a skew symmetric matrix associated with $\boldsymbol{\omega}$ such that $\boldsymbol{\omega}^{\times}\mathbf{y} = \boldsymbol{\omega}\times\mathbf{y}$ for any $\boldsymbol{\omega}, \mathbf{y}\in\mathbb{R}^3$ with $\times$ being the vector cross product.

Taking the time derivative of $\mathbf{x} = \mathbf{R}^{\top}\mathbf{e}_3$ and using \eqref{rotation_matrix_dynamics}, one obtains
\begin{equation}\label{xdot_expression_rotation_matrix}
    \dot{\mathbf{x}} = -\boldsymbol{\omega}^{\times}\mathbf{R}^{\top}\mathbf{e}_3 = \mathbf{x}^{\times}\boldsymbol{\omega}.
\end{equation}
Choosing $\boldsymbol{\omega}$ as follows:
\begin{equation}\label{angular_velocity_reduced_attitude_control}
    \boldsymbol{\omega} = -\mathbf{x}^{\times}\mathbf{u},
\end{equation}
and using the fact that $-\left(\mathbf{x}^{\times}\right)^2 = \mathbf{I} - \mathbf{x}\mathbf{x}^\top = \mathbf{P}(\mathbf{x})$, leads to $\dot{\mathbf{x}} = \mathbf{P}(\mathbf{x})\mathbf{u}$ which corresponds to model \eqref{system_dynamics}.
Therefore, the control laws developed in the present paper (\eqref{negative_gradient_control_law} or \eqref{n-sphere-control-law}) can be used to obtain $\mathbf{u}$, and then the angular velocity can be obtained using \eqref{angular_velocity_reduced_attitude_control}.

\subsection{Full attitude control using the unit-quaternion}
The attitude of a rigid body with respect to the inertial frame can be described by a four-parameter representation, namely a unit-quaternion.
To denote the unit-quaternion, we use $\mathbf{x} = \begin{bmatrix}\eta,\mathbf{q}^\top\end{bmatrix}\in\mathbb{S}^3$, where $\eta\in\mathbb{R}$ and $\mathbf{q}\in\mathbb{R}^3$.
The quaternion kinematics is given by
\begin{equation}\label{quaternion_kinematics}
    \dot{\mathbf{x}} = \frac{1}{2}\mathbf{A}(\mathbf{x})\boldsymbol{\omega} = \frac{1}{2}\begin{bmatrix}-\mathbf{q}^\top\\\eta\mathbf{I}_3 + \mathbf{q}^{\times}\end{bmatrix}\boldsymbol{\omega},
\end{equation}
where $\boldsymbol{\omega}\in\mathbb{R}^3$ is the angular velocity of the rigid body in the body frame.
Picking 
\begin{equation}\label{omega_representation_using_control_u}
    \boldsymbol{\omega}=2\mathbf{A}(\mathbf{x})^\top \mathbf{u},
\end{equation}
and using the fact that $\mathbf{A}(\mathbf{x})^\top \mathbf{A}(\mathbf{x})=\mathbf{I}_3$, and $\mathbf{A}(\mathbf{x})\mathbf{A}(\mathbf{x})^\top = \mathbf{P}(\mathbf{x})$ for all $\mathbf{x}\in\mathbb{S}^3$, one can show that $\dot{\mathbf{x}}=\mathbf{P}(\mathbf{x})\mathbf{u}$ which corresponds to model \eqref{system_dynamics}. Therefore, the control laws developed in the present paper (\eqref{negative_gradient_control_law} or \eqref{n-sphere-control-law}) can be used to obtain $\mathbf{u}$, and then the angular velocity can be obtained using \eqref{omega_representation_using_control_u}.

\section{Simulation results}
\label{section:simulation}
First, we provide a geometric procedure for the construction of a star-shaped set on $\mathbb{S}^n$ by projecting an $n$-dimensional star-shaped set embedded in $n+1$-dimensional Euclidean space onto $\mathbb{S}^n$.

\subsection{Geometric construction of a star-shaped set \texorpdfstring{$\mathcal{U}_i$}{} on the \texorpdfstring{$n$}{}-sphere}
Consider a line segment $\mathcal{L}_s(\mathbf{a}, \mathbf{b})$, defined as in Section \ref{notations}, connecting any two points $\mathbf{a}, \mathbf{b}\in\mathbb{R}^{n+1}$ such that $\mathbf{0}\notin\mathcal{L}_s(\mathbf{a}, \mathbf{b})$.
Define a set $\mathcal{Q}(\mathbf{a}, \mathbf{b})$ as follows:
\begin{equation}\label{Q_set_definition}
    \mathcal{Q}(\mathbf{a}, \mathbf{b}) = \{\mathbf{x}\in\mathbb{S}^n\mid\mathbf{x} = \psi(\mathbf{p}), \mathbf{p}\in\mathcal{L}_s(\mathbf{a}, \mathbf{b})\},
\end{equation}
where the mapping $\psi:\mathbb{R}^{n+1}\setminus\{\mathbf{0}\}\to\mathbb{S}^n$ is given by
\begin{equation}\label{psi_function_definition}
    \psi(\mathbf{p}) = \frac{\mathbf{p}}{\|\mathbf{p}\|}.
\end{equation}
Since $\mathbf{0}\notin\mathcal{L}_s(\mathbf{a}, \mathbf{b})$, the set $\mathcal{Q}(\mathbf{a}, \mathbf{b})$ is well-defined.
In the next lemma, we show that for any $\mathbf{a}, \mathbf{b}\in\mathbb{R}^{n+1}$ with $\mathbf{0}\notin\mathcal{L}_s(\mathbf{a}, \mathbf{b})$, the set $\mathcal{Q}(\mathbf{a}, \mathbf{b})$ coincides with the geodesic $\mathcal{G}(\psi(\mathbf{a}), \psi(\mathbf{b}))$.

\begin{lemma}\label{lemma:Q_set_is_geodesic}
    Let $\mathbf{a}, \mathbf{b}\in\mathbb{R}^{n+1}$ and $\mathbf{0}\notin\mathcal{L}_s(\mathbf{a}, \mathbf{b})$. Then, $\mathcal{G}(\psi(\mathbf{a}), \psi(\mathbf{b})) = \mathcal{Q}(\mathbf{a}, \mathbf{b})$, where the geodesic $\mathcal{G}(\psi(\mathbf{a}), \psi(\mathbf{b}))$ and the set $\mathcal{Q}(\mathbf{a}, \mathbf{b})$ are defined in \eqref{geodesic_expression} and \eqref{Q_set_definition}, respectively.
\end{lemma}
\proof{See Appendix \ref{proof:lemma:Q_set_is_geodesic}.}

\begin{figure}
    \centering
    \includegraphics[width=0.8\linewidth]{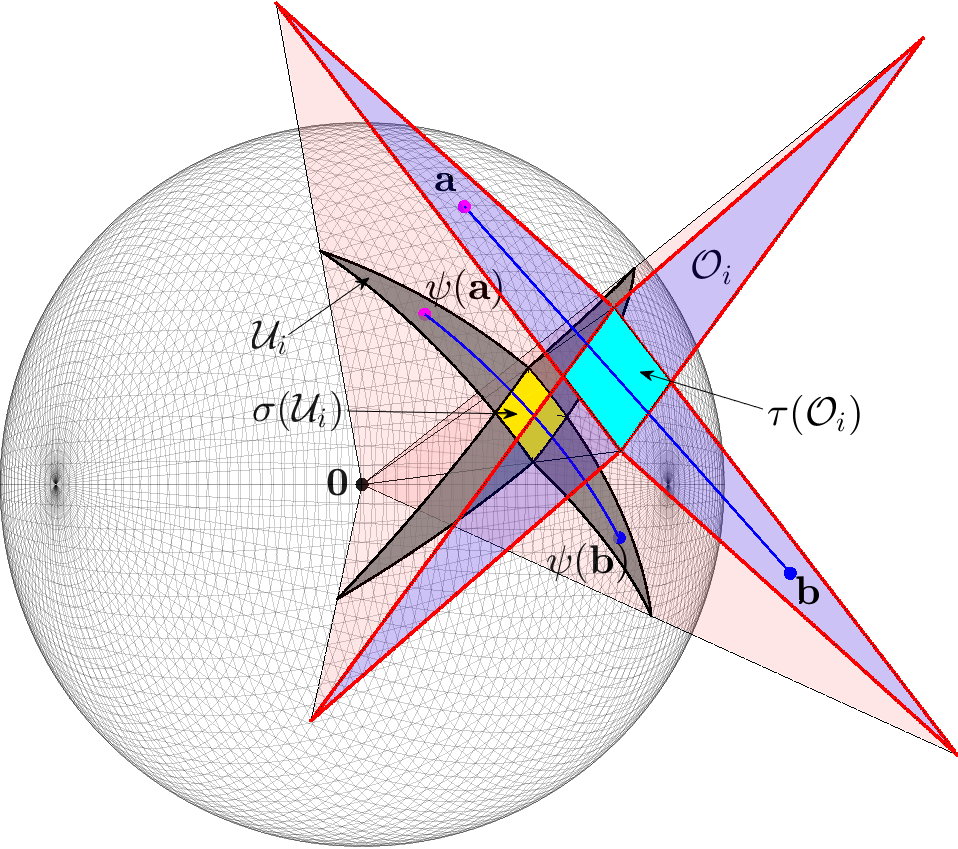}
    \caption{Geometric construction of a star-shaped set $\mathcal{U}_i$ on $\mathbb{S}^2$.}
    \label{fig:g-star_set-construction}
\end{figure}

Lemma \ref{lemma:Q_set_is_geodesic} states that if a line segment $\mathcal{L}_s(\mathbf{a}, \mathbf{b})$ does not pass through $\mathbf{0}$ for some $\mathbf{a}, \mathbf{b}\in\mathbb{R}^{n+1}$, then the curve $\mathcal{Q}(\mathbf{a}, \mathbf{b})$, obtained by projecting $\mathcal{L}_s(\mathbf{a}, \mathbf{b})$ onto the $n$-sphere, coincides with the unique geodesic connecting $\psi(\mathbf{a})$ and $\psi(\mathbf{b})$.
Consequently, if two line segments $\mathcal{L}_s(\mathbf{a}_1, \mathbf{b}_1)$ and $\mathcal{L}_s(\mathbf{a}_2, \mathbf{b}_2)$ satisfy $\mathbf{0}\notin\mathcal{L}_s(\mathbf{a}_1, \mathbf{b}_1)$ and $\mathbf{0}\notin\mathcal{L}_s(\mathbf{a}_2, \mathbf{b}_2)$, and intersect each other in $\mathbb{R}^{n+1}$, then the corresponding geodesics $\mathcal{G}(\psi(\mathbf{a}_1), \psi(\mathbf{b}_1))$ and $\mathcal{G}(\psi(\mathbf{a}_2), \psi(\mathbf{b}_2))$, which coincide with $\mathcal{Q}(\mathbf{a}_1, \mathbf{b}_1)$ and $\mathcal{Q}(\mathbf{a}_2, \mathbf{b}_2)$, respectively, also intersect each other on the $n$-sphere.
This property allows us to construct a star-shaped set $\mathcal{U}_i$ on the $n$-sphere by projecting every point of a given $n$-dimensional star-shaped set $\mathcal{O}_i$ embedded in $\mathbb{R}^{n+1}$ on the $n$-sphere, provided that $\mathbf{0}\notin\mathcal{O}_i$, as discussed next.

\begin{figure}[ht]
    \centering
    \includegraphics[width=0.7\linewidth]{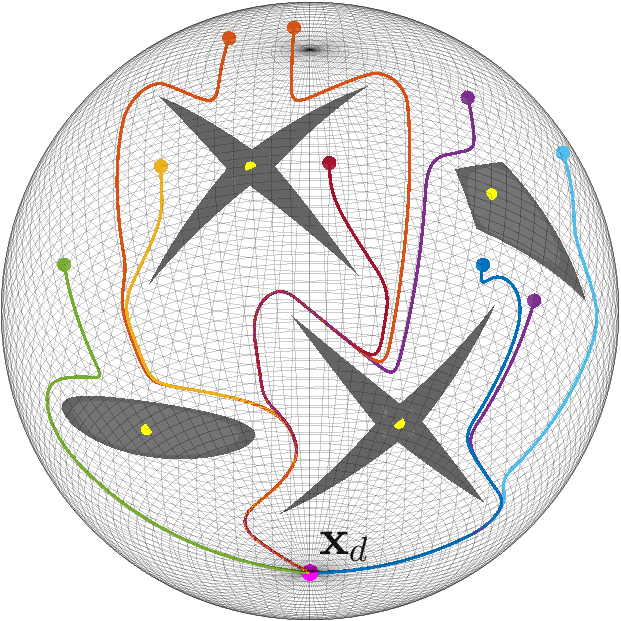}
    \caption{$\mathbf{x}$-trajectories safely converging to $\mathbf{x}_d$.}
    \label{fig:2-sphere-star-shaped-CS}
\end{figure}

\begin{figure}[ht]
    \centering
    \includegraphics[width=1\linewidth]{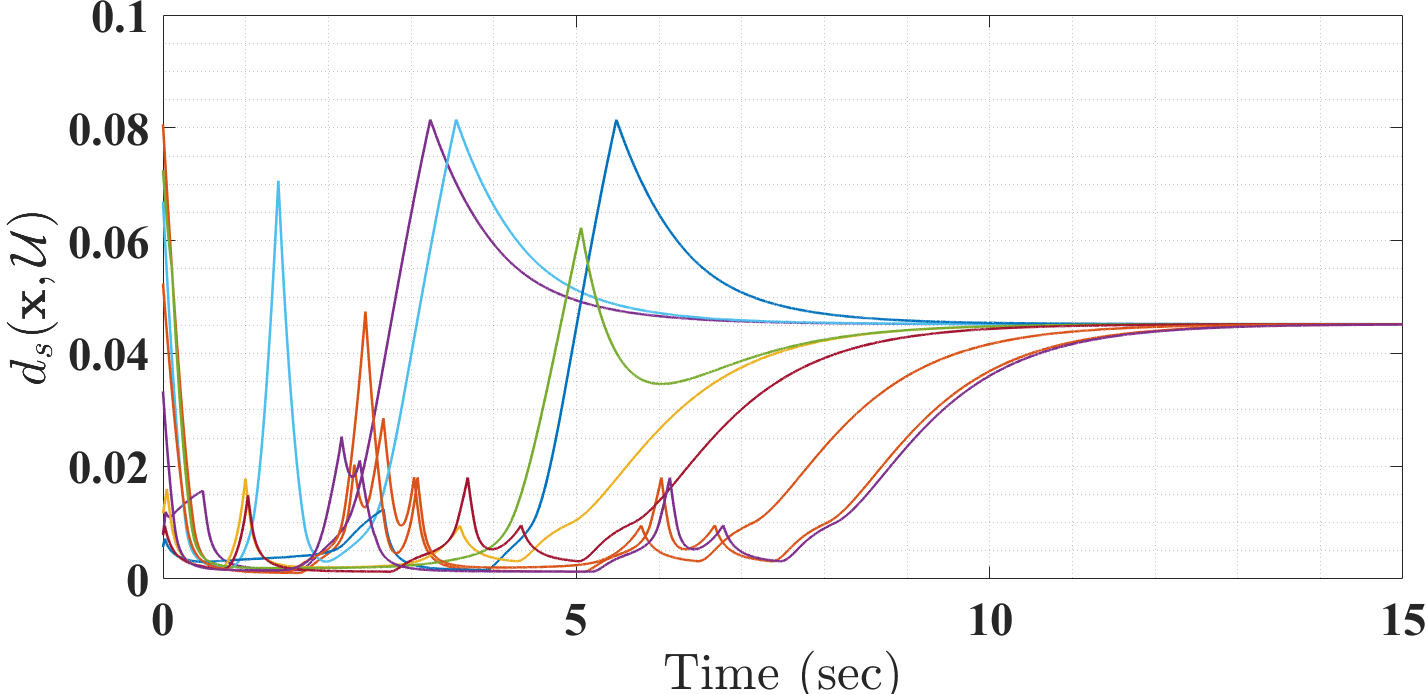}
    \caption{ $d_s(\mathbf{x}, \mathcal{U})$ versus time.}
    \label{fig:2-sphere-star-shaped-safety}
\end{figure}


\begin{figure*}[ht]
\centering
\subfloat[][]{\includegraphics[width =0.32\linewidth]{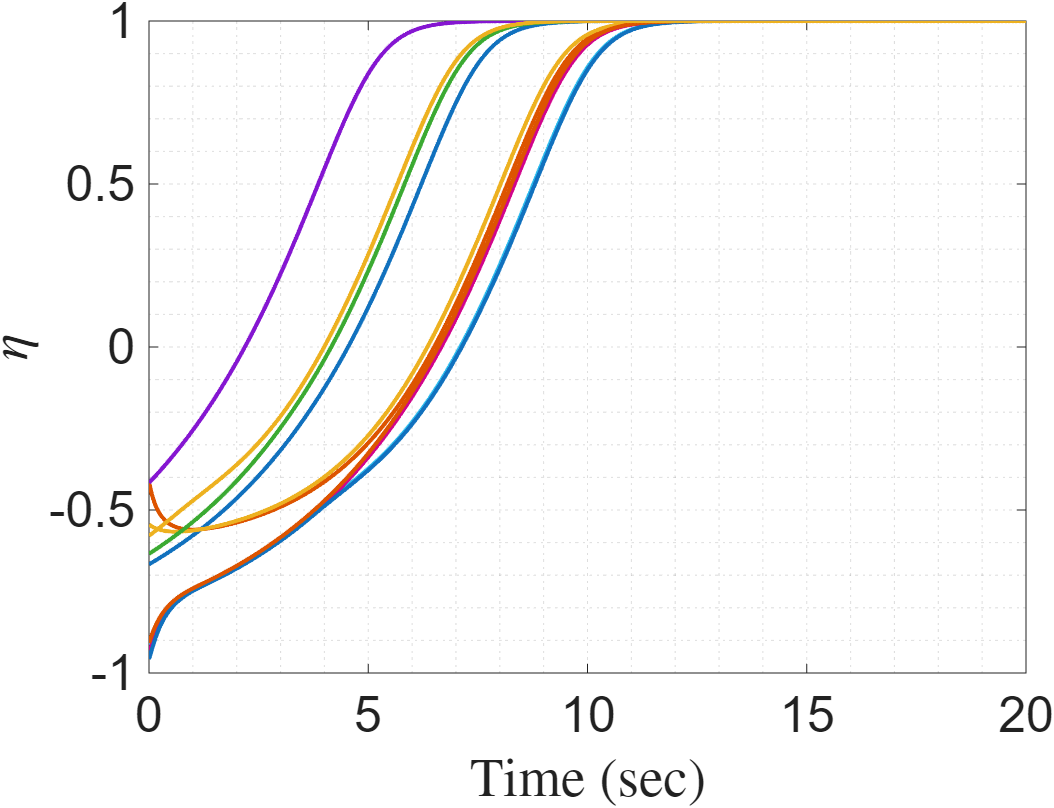}\label{fig:eta_plot}}\hspace{0.2cm}
\subfloat[][]{\includegraphics[width =0.32\linewidth]{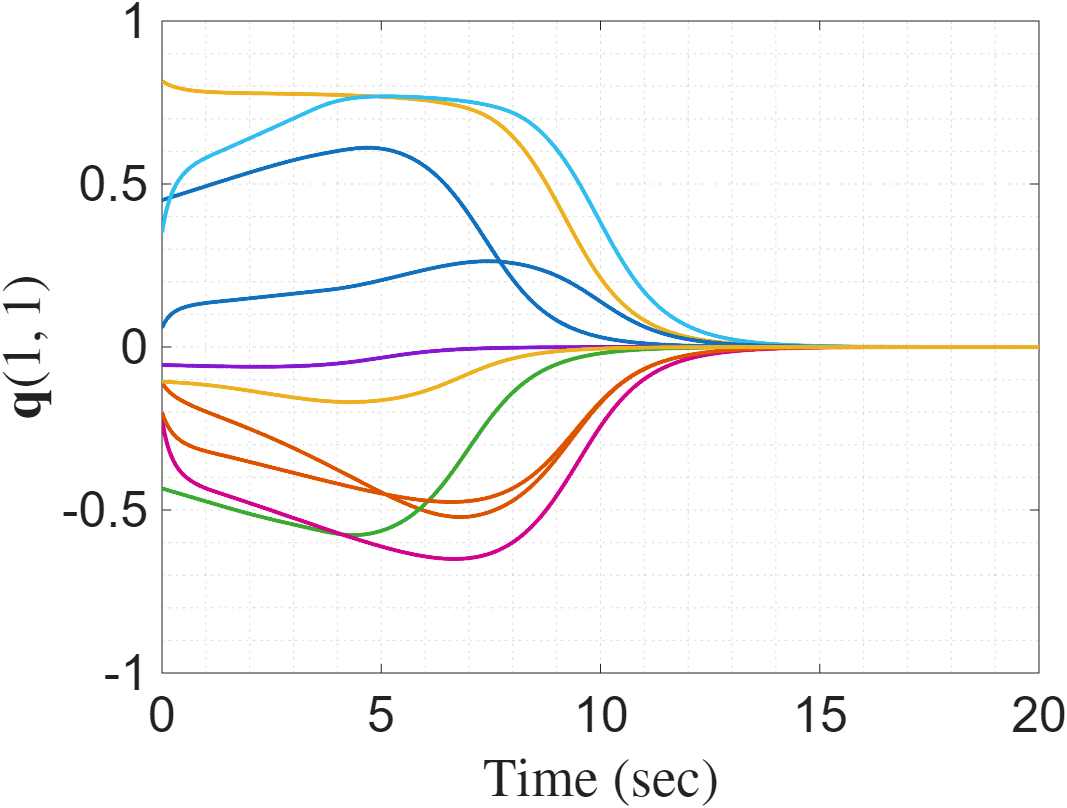}\label{fig:q1_plot}}\hspace{0.2cm}
\subfloat[][]{\includegraphics[width =0.32\linewidth]{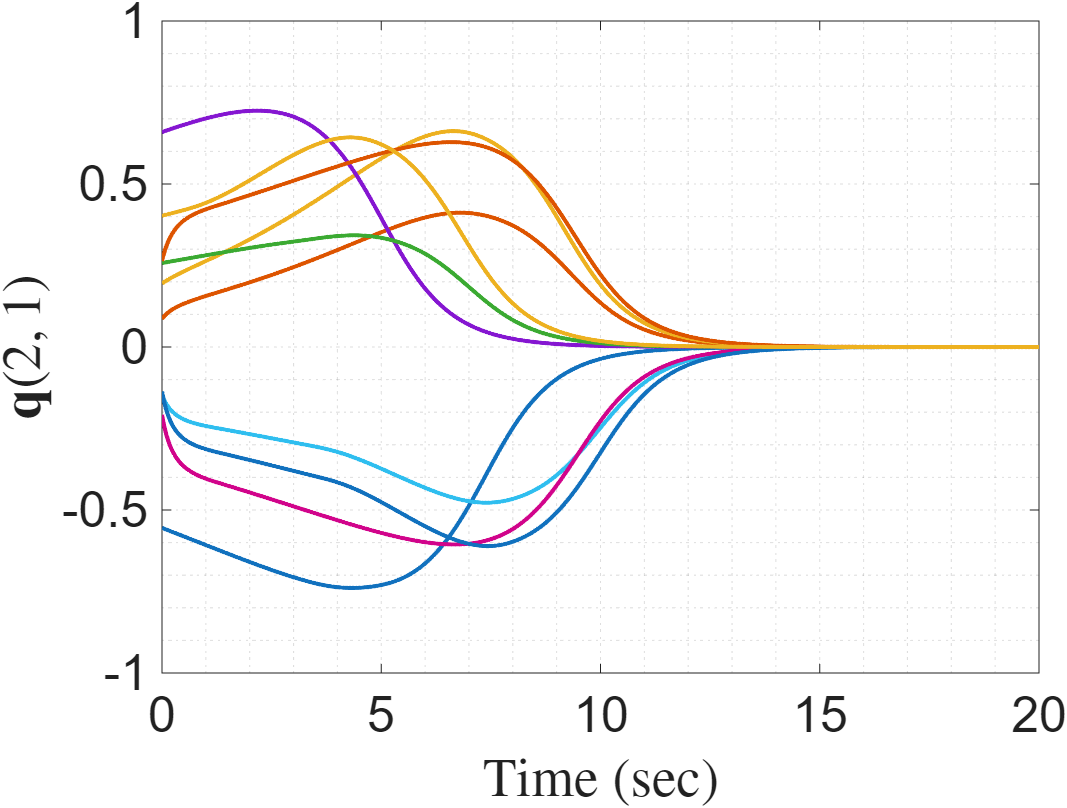}\label{fig:q2_plot}}\\
\subfloat[][]{\includegraphics[width =0.32\linewidth]{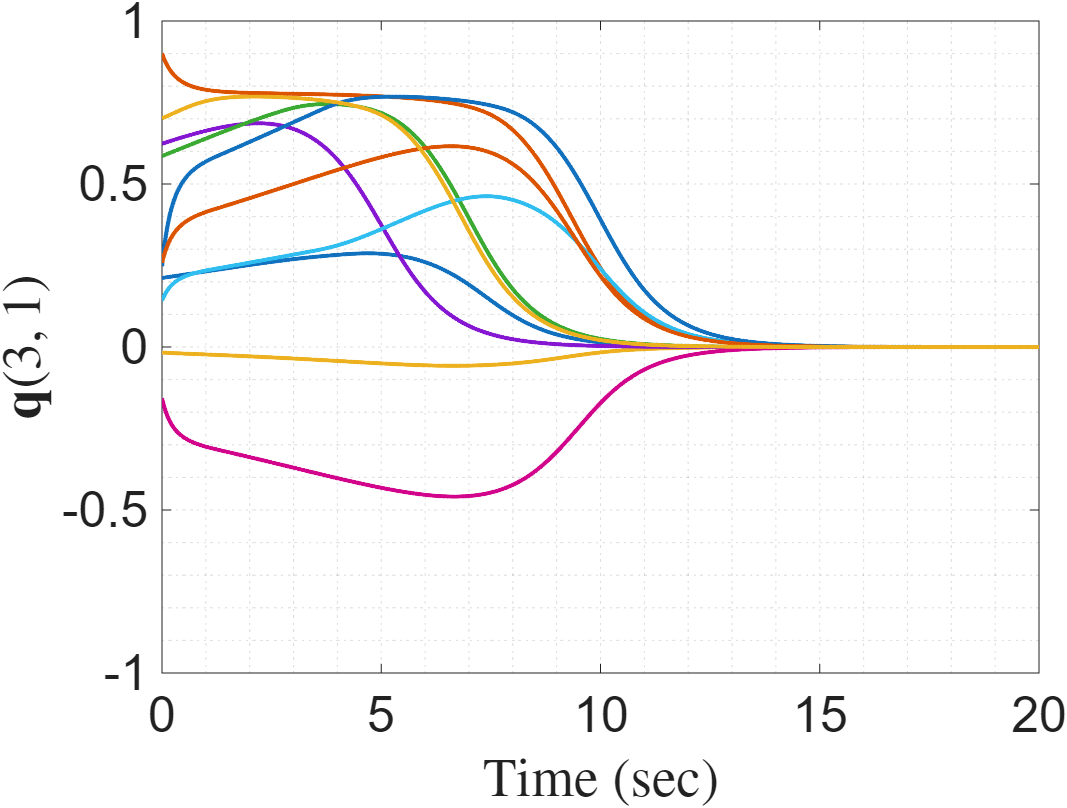}\label{fig:q3_plot}}\hspace{0.2cm}
\subfloat[][]{\includegraphics[width =0.32\linewidth]{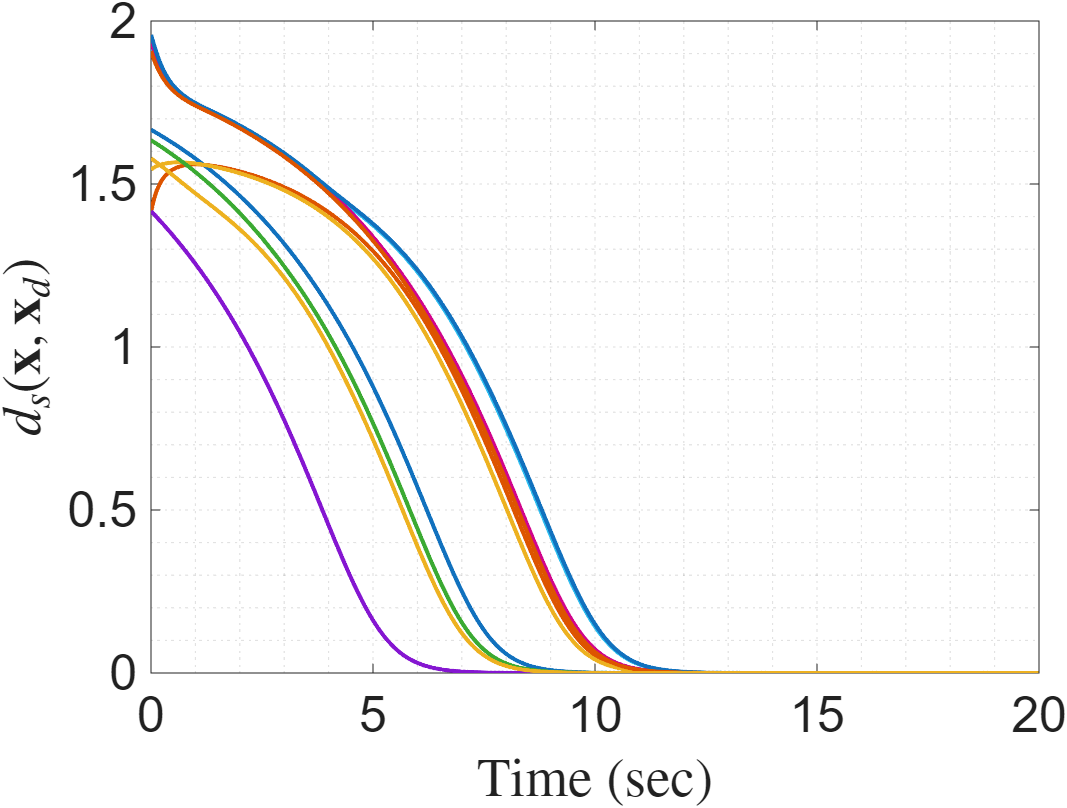}\label{fig:dxd_plot}}\hspace{0.2cm}
\subfloat[][]{\includegraphics[width =0.32\linewidth]{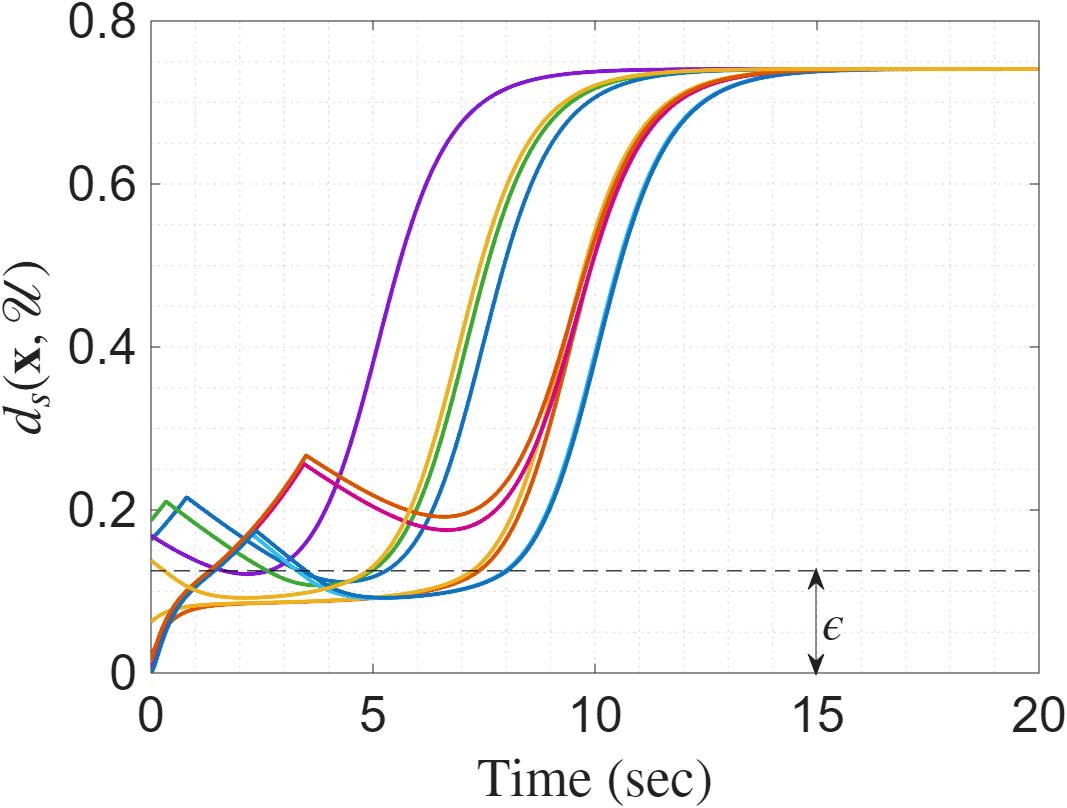}\label{fig:dunsafe_plot}}
\caption{Implementation of the closed-loop system \eqref{quaternion_kinematics}-\eqref{omega_representation_using_control_u} with $\mathbf{u}$ defined in \eqref{negative_gradient_control_law}. (a)-(d)  $\mathbf{x}$-trajectories converging to $\mathbf{x}_d =[1, 0, 0, 0]^\top$, (e) $d_s(\mathbf{x}, \mathbf{x}_d)$ versus time, (f) $d_s(\mathbf{x}, \mathcal{U})$ versus time.}\label{result:quaternion_kinematic_implementation}
\end{figure*}

In view of Lemma \ref{lemma:Q_set_is_geodesic}, we construct a star-shaped set $\mathcal{U}_i$ on $\mathbb{S}^n$ as follows:
\begin{equation}\label{construction_of_star_shaped_set}
    \mathcal{U}_i = \{\mathbf{x}\in\mathbb{S}^n\mid\mathbf{x} = \psi(\mathbf{p}), \mathbf{p}\in\mathcal{O}_i\},
\end{equation}
where $\mathcal{O}_i$ is an $n$-dimensional star-shaped set\footnote{A set $\mathcal{A}\subset\mathbb{R}^n$ is a star-shaped set, if there exists $\mathbf{a}\in\mathcal{A}$ such that $\mathcal{L}_s(\mathbf{x}, \mathbf{a})\subset\mathcal{A}$ for all $\mathbf{x}\in\mathcal{A}$.} embedded in the Euclidean space $\mathbb{R}^{n+1}$ such that $\mathbf{0}\notin\mathcal{O}_i$.
Since $\mathcal{O}_i$ is a star-shaped set, analogous to $\sigma(\mathcal{U}_i)$ \eqref{sigma_set}, one can define $\tau(\mathcal{O}_i)$ as follows:
\begin{equation*}
    \tau(\mathcal{O}_i) := \{\mathbf{a}\in\mathcal{O}_i\mid\forall \mathbf{x}\in\mathcal{O}_i, \mathcal{L}_s(\mathbf{a}, \mathbf{x})\subset\mathcal{O}_i\},
\end{equation*}
which is a subset of $\mathcal{O}_i$ such that for every $\mathbf{a}\in\tau(\mathcal{O}_i)$ the line segments $\mathcal{L}_s(\mathbf{a}, \mathbf{x})$ connecting $\mathbf{a}$ to any other point $\mathbf{x}$ in $\mathcal{O}_i$ always belong to $\mathcal{O}_i$.
Since $\psi(\cdot)$ maps every point in $\mathcal{O}_i$ to $\mathbb{S}^n$ while preserving direction, it follows that if $\mathcal{L}_s(\mathbf{a}, \mathbf{b})\subset\mathcal{O}_i$ for any pair $\mathbf{a}, \mathbf{b}\in\mathcal{O}_i$, then $\mathcal{G}(\psi(\mathbf{a}), \psi(\mathbf{b}))\subset\mathcal{U}_i$, as illustrated in Fig. \ref{fig:g-star_set-construction}.
Therefore, using $\tau(\mathcal{O}_i)$, the set $\sigma(\mathcal{U}_i)$ can be identified as
\begin{equation*}
    \sigma(\mathcal{U}_i) = \{\psi(\mathbf{a})\in\mathbb{S}^n\mid\mathbf{a}\in\tau(\mathcal{O}_i)\}.
\end{equation*}

\subsection{Constrained stabilization on \texorpdfstring{$2$-sphere}{}}
We consider $\mathbb{S}^2$ with $4$ star-shaped constraints, as shown in Fig. \ref{fig:2-sphere-star-shaped-CS}. 
The location of constant unit vectors $\mathbf{g}_i$ is denoted using yellow dots.
The scalar parameters $k_1, \kappa$ and $\epsilon$ are set to $1, 1$ and $0.01$, respectively.
The $\mathbf{x}$-trajectories are initialized at $9$ different initial locations and asymptotically converge to the target point at $\mathbf{x}_d$, as depicted in Fig. \ref{fig:2-sphere-star-shaped-CS}.
The proposed feedback controller \eqref{n-sphere-control-law} ensures safety
for all time $t\geq 0$, as illustrated in Fig. \ref{fig:2-sphere-star-shaped-safety}.

\subsection{Constrained stabilization on \texorpdfstring{$3$}{}-sphere}
We consider $\mathbb{S}^3$ with $7$ conic constraints, defined as in \eqref{definition:conic_constraints}, where the constant unit vectors $\mathbf{g}_i$ are set to $[0, 1, 0, 0]^\top$, $[0, 0, 1, 0]^\top$, $[0, 0, 0, 1]^\top$, $[0, -1, 0, 0]^\top$, $[0, 0, -1, 0]^\top$, $[0, 0, 0, -1]^\top$ and $[-1, 0, 0, 0]^\top$. 
For each $i\in\mathbb{I}$, the parameters $\xi_i$ are set to $\frac{\pi}{12}$ rad.
Notice that the unsafe regions $\mathcal{U}_i$ satisfy Assumption \ref{assumption:non-overlapping-constraints} with $\delta = \frac{\pi}{6}$ rad.
The target location $\mathbf{x}_d$ is set to $[1, 0, 0, 0]^\top$.
The parameters $k_1$ and $\epsilon$, used in \eqref{negative_gradient_control_law}, are set to $1$ and $0.125$, respectively. 
The closed-loop system \eqref{quaternion_kinematics}-\eqref{omega_representation_using_control_u} is initialized at $10$ different initial conditions $\mathbf{x}(0)\in\mathcal{M}_0$.
The $\mathbf{x}$-trajectories asymptotically converge to $\mathbf{x}_d$, as illustrated in Fig. \ref{fig:eta_plot}-\ref{fig:q3_plot}.
The proposed feedback controller \eqref{negative_gradient_control_law}, used in \eqref{omega_representation_using_control_u}, ensures safety
for all time $t\geq 0$, as depicted in Fig. \ref{fig:dunsafe_plot}.

For the next simulation, we consider $\mathbb{S}^3$ with a star-shaped constraint $\mathcal{U}_1$, which is constructed from $\mathcal{O}_1$ using \eqref{construction_of_star_shaped_set}, where $\mathcal{O}_1$ is a three-dimensional set embedded in $\mathbb{R}^4$, and it is given by
\begin{equation*}
    \mathcal{O}_1 = \{\mathbf{y}\in\mathbb{R}^4\mid \mathbf{y} = \mathbf{g}_1 + \mathbf{p}, \mathbf{p}\in\mathcal{O}_0\},
\end{equation*}
with $\mathbf{g}_1\in\mathbb{S}^3$.
The $3$-dimensional star-shaped set $\mathcal{O}_0$ embedded in $\mathbb{R}^4$, as illustrated in Fig. \ref{fig:3S_obstacle_plot}, is defined as
\begin{equation}\label{set_O_0_definition}
    \mathcal{O}_0 = \left\{\mathbf{p}\in\mathbb{R}^4\mid |p_1|^{0.4} + |p_2|^{0.4} + |p_3|^{0.4} \leq 1.5, p_4 = 0\right\},  
\end{equation}
$\mathbf{p} = [p_1, p_2, p_3, p_4]^\top$.
The unit vector $\mathbf{g}_1$ and the target location $\mathbf{x}_d$ are set to $[-0.5, -0.5, -0.5, -0.5]^\top$ and $[1, 0, 0, 0]^\top$, respectively.
The parameters $k_1$, $\kappa$ and $\epsilon$, used in \eqref{n-sphere-control-law} and \eqref{individual_control_input_vector_design} are chosen as $1$, $1$ and $0.1$, respectively.
The $\mathbf{x}$-trajectories initialized at $10$ different initial conditions $\mathbf{x}(0)\in\mathcal{M}_0$ asymptotically converge to $\mathbf{x}_d$, as depicted in Fig. \ref{fig:3S_eta_plot}-\ref{fig:3S_q3_plot}.
The proposed feedback controller \eqref{n-sphere-control-law}, used in \eqref{omega_representation_using_control_u}, ensures safety
for all $t\geq 0$, as shown in Fig. \ref{fig:3S_dunsafe_plot}.

\section{Conclusion}
\label{section:conclusion}
In this work, we proposed a feedback control law for the constrained stabilization problem on the $n$-sphere. 
Unlike much of the existing literature \cite{berkane2021constrained}, in which the unsafe regions are typically modeled as conic constraints, we model the unsafe region as a union of star-shaped constraints on the $n$-sphere. 
This offers a more flexible characterization of the unsafe region, potentially enabling a larger safe region for stabilization purposes.
The proposed feedback control law combines an attractive vector field, which guides the system state $\mathbf{x}$ along the geodesic toward the target, with a repulsive vector field that steers $\mathbf{x}$ away from the unsafe region. 
Almost global asymptotic stability of the target location is rigorously proven for the closed-loop system \eqref{system_dynamics}-\eqref{n-sphere-control-law}.
An interesting future work would be the constrained dynamic stabilization on $\mathbb{S}^n$, where the input would be the acceleration rather than the velocity.
In addition, future work may incorporate hybrid control techniques \cite{goebel2012hybrid} to achieve global asymptotic stability of the desired target point.

\appendix
\subsection{Proof of Lemma \ref{lemma:reverse_geodesic_always_stay_outside}}\label{proof:lemma:reverse_geodesic_always_stay_outside}
Since $\mathbf{g}\in\mathcal{A}^{\circ}$ and $\mathbf{x}\in\partial\mathcal{A}$, one has $\mathbf{x}\ne\mathbf{g}$, and the geodesic $\mathcal{G}(\mathbf{x}, -\mathbf{g})$ exists and is unique. 
We proceed by contradiction.
Assume that there exists $\mathbf{p}\in\mathcal{G}(\mathbf{x}, -\mathbf{g})$ such that $\mathbf{p}\in\mathcal{A}^{\circ}$.
Since $\mathbf{p}\in\mathcal{A}^{\circ}$, there exists $\mu > 0$ such that $\mathcal{D}_{\mu}(\mathbf{p})\subset\mathcal{A}$.
Since $\mathcal{A}$ is a star-shaped set on $\mathbb{S}^n$ and $\mathbf{g}\in\sigma(\mathcal{A})$, $\mathcal{D}_{\mu}(\mathbf{p})\subset\mathcal{A}$ implies that $\mathcal{A}_{\mu}(\mathbf{p}, \mathbf{g})\subset\mathcal{A}$, where the set $\mathcal{A}_{\mu}(\mathbf{p}, \mathbf{g})$ is defined as
\[\mathcal{A}_{\mu}(\mathbf{p}, \mathbf{g}) = \{\mathbf{a}\in\mathbb{S}^n\mid\mathbf{a}\in\mathcal{G}(\mathbf{q}, \mathbf{g}), \mathbf{q}\in\mathcal{D}_{\mu}(\mathbf{p})\}.\]

Now, since $\mathbf{p}\in\mathcal{G}(\mathbf{x}, -\mathbf{g})\setminus\{\mathbf{x}, -\mathbf{g}\}$, one has $\mathbf{x}\in\mathcal{G}(\mathbf{p}, \mathbf{g})$ and it follows that $\mathbf{x}\in\mathcal{A}_{\mu}(\mathbf{p}, \mathbf{g})$.
Owing to the positive sectional curvature of $\mathbb{S}^n$ \cite[Ch. 6, Ex. 2.8]{do1992riemannian}, one can show that $\mathbf{x}\in\left(\mathcal{A}_{\mu}(\mathbf{p}, \mathbf{g})\right)^{\circ}$.
Consequently, since $\mathcal{A}_{\mu}(\mathbf{p}, \mathbf{g})\subset\mathcal{A}$, it follows that $\mathbf{x}\in\mathcal{A}^{\circ}$.
However, this contradicts the fact that $\mathbf{x}\in\partial\mathcal{A}$, and the proof is complete.

\subsection{Proof of Lemma \ref{lemma:constraint_separation}}\label{proof:lemma:constraint_separation}
For each $i\in\mathbb{I}$ and for any $\mathbf{x}\in\partial\mathcal{N}_{\epsilon}(\mathcal{U}_i)\setminus\partial\mathcal{U}_i$, one has $d_s(\mathbf{x}, \mathcal{U}_i) = \epsilon$.
Therefore, the minimum positive angle between any $\mathbf{x}\in\partial\mathcal{N}_{\epsilon}(\mathcal{U}_i)\setminus\partial\mathcal{U}_i$ and any unit vector in $\mathcal{U}_i$ is given by $\theta_{\epsilon} = \arccos(1 - \epsilon)$, where for any $p\in[-1, 1]$, the inverse cosine function $\arccos(p)$ is restricted to $[0, \pi]$. 
To guarantee $\mathcal{D}_{\epsilon}(\mathcal{U}_i)\cap\mathcal{D}_{\epsilon}(\mathcal{U}_j) = \emptyset$ for all $i, j\in\mathbb{I}$ with $i\ne j$, it suffices to choose $\epsilon > 0$ such that $\theta_{\epsilon} \in (0,  \delta)$.
Consequently, setting $\arccos(1-\epsilon) < \delta$, taking cosine on both sides and rearranging the terms, one obtains $\epsilon < 1 - \cos(\delta)$. This completes the proof.

\begin{figure*}[ht]
\centering
\subfloat[][]{\includegraphics[width =0.32\linewidth]{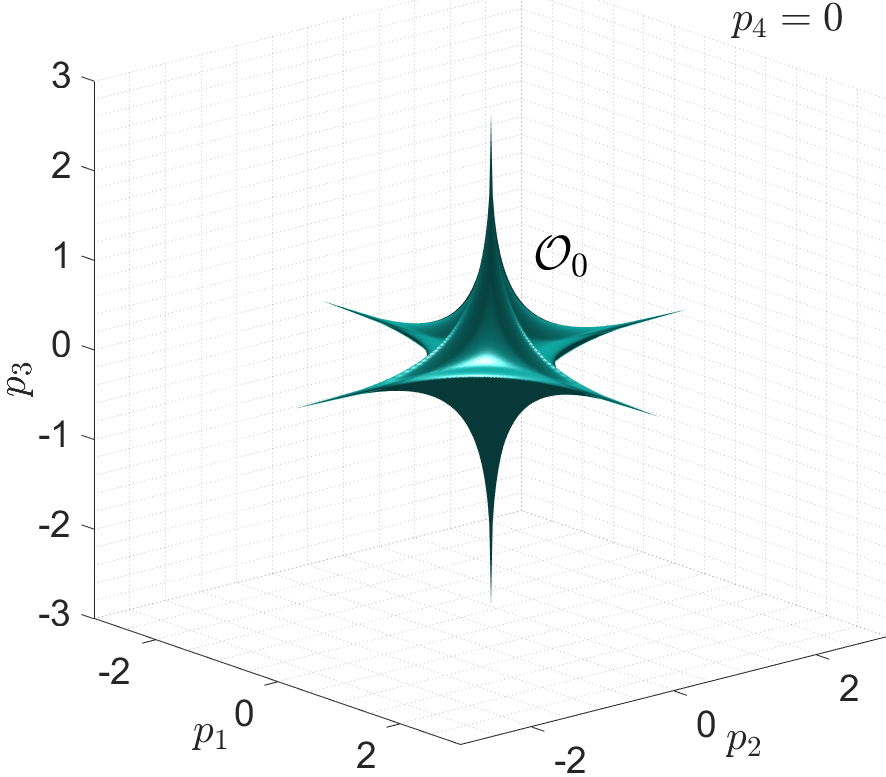}\label{fig:3S_obstacle_plot}}\hspace{0.2cm}
\subfloat[][]{\includegraphics[width =0.32\linewidth]{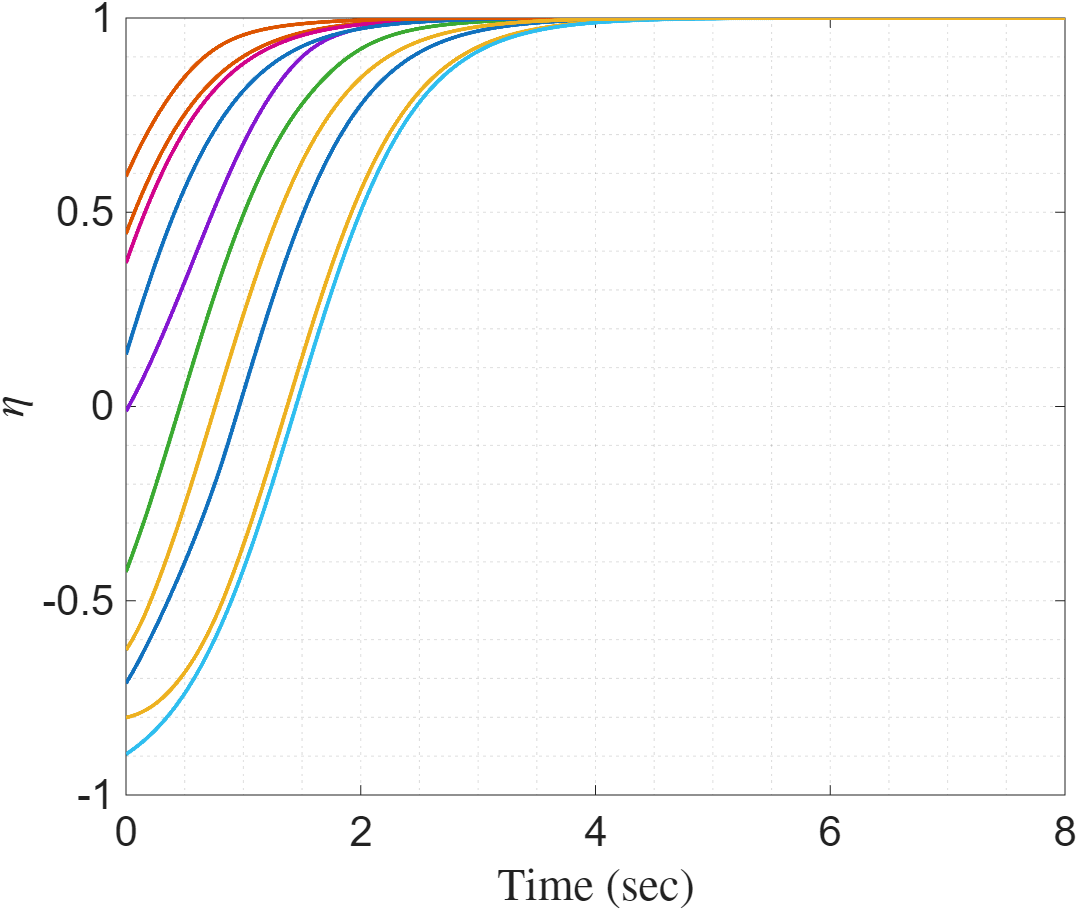}\label{fig:3S_eta_plot}}\hspace{0.2cm}
\subfloat[][]{\includegraphics[width =0.32\linewidth]{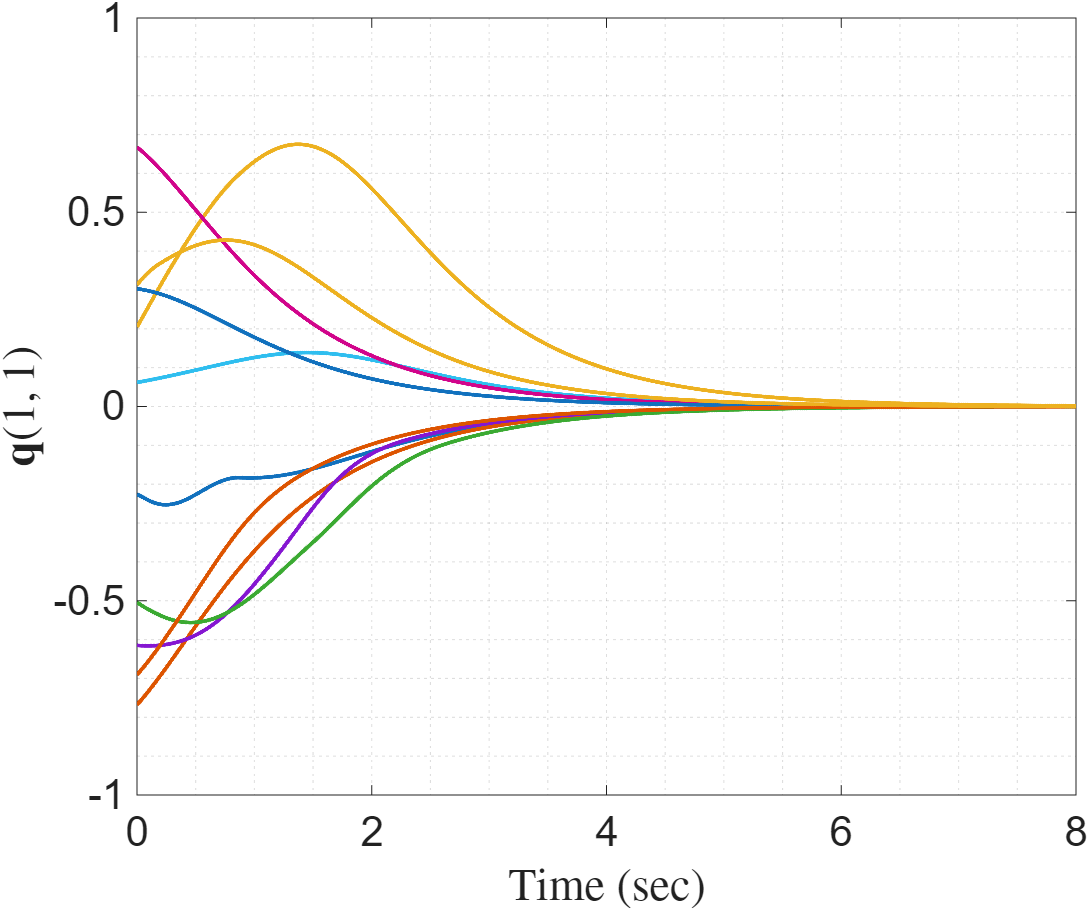}\label{fig:3s_q1_plot}}\\
\subfloat[][]{\includegraphics[width =0.32\linewidth]{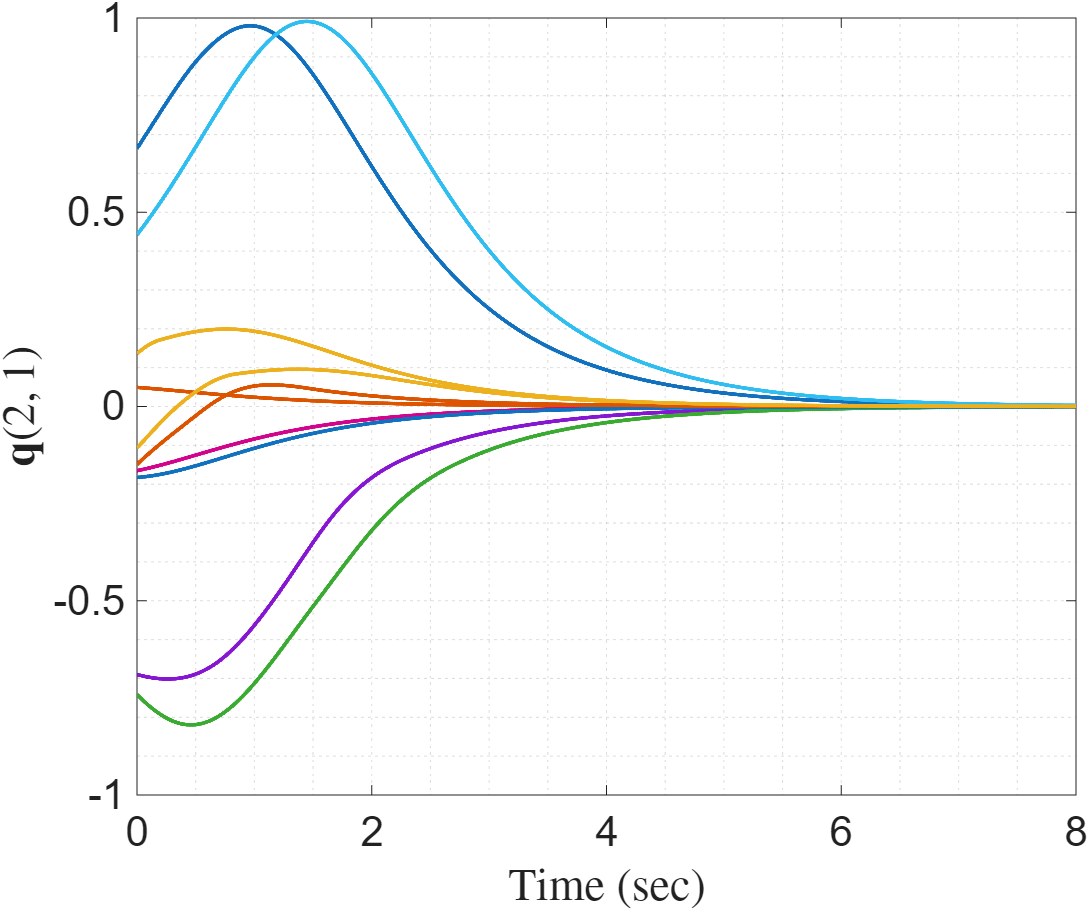}\label{fig:3S_q2_plot}}\hspace{0.2cm}
\subfloat[][]{\includegraphics[width =0.32\linewidth]{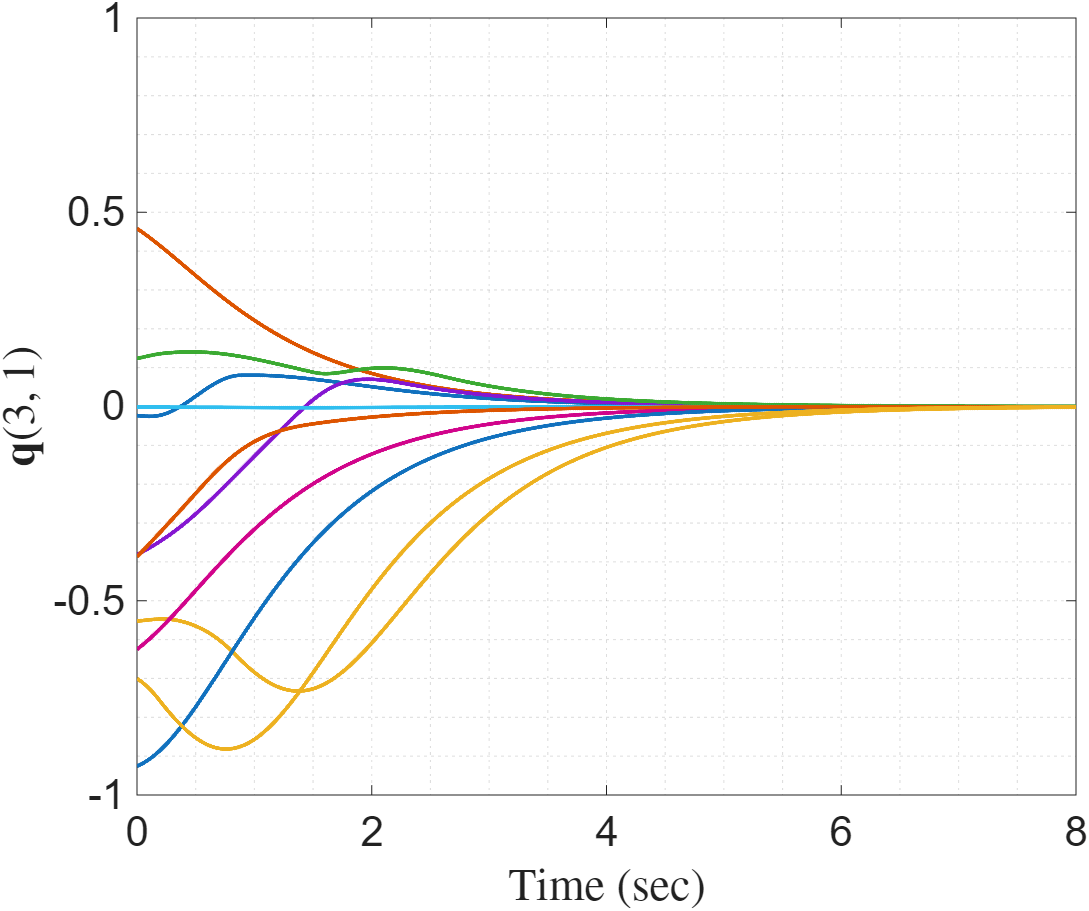}\label{fig:3S_q3_plot}}\hspace{0.2cm}
\subfloat[][]{\includegraphics[width =0.32\linewidth]{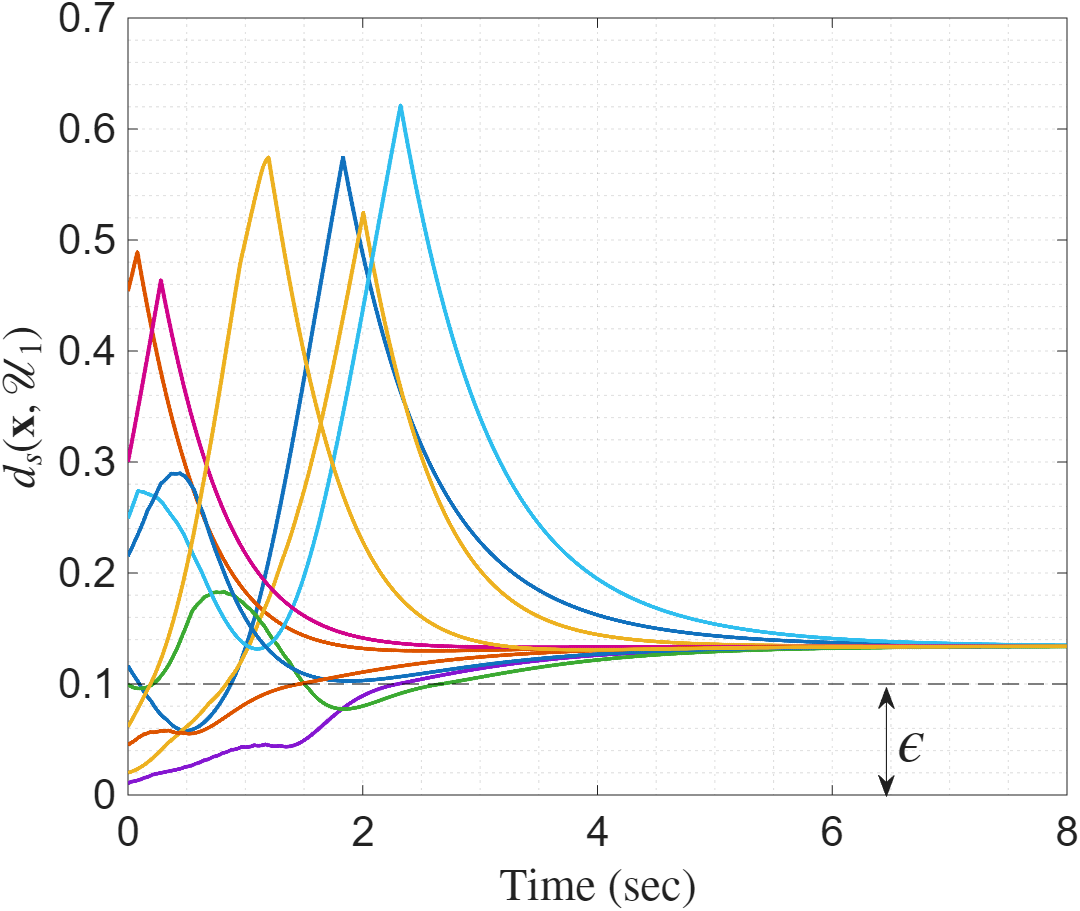}\label{fig:3S_dunsafe_plot}}
\caption{Implementation of the closed-loop system \eqref{quaternion_kinematics}-\eqref{omega_representation_using_control_u} with $\mathbf{u}$ defined in \eqref{n-sphere-control-law}. (a) Set $\mathcal{O}_0$ \eqref{set_O_0_definition}, (b)-(e)  $\mathbf{x}$-trajectories converging to $\mathbf{x}_d =[1, 0, 0, 0]^\top$, (f) $d_s(\mathbf{x}, \mathcal{U}_1)$ versus time.}\label{result:quaternion_kinematic_implementation_star_obstacle}
\end{figure*}

\subsection{Proof of Lemma \ref{lemma:properties_scalar_function}}\label{proof:lemma:properties_scalar_function}
We show that $W(\mathbf{x})$ is well-defined for all $\mathbf{x}\in\mathcal{M}_0$.
Since $\mathbf{x}\in\mathbb{S}^n$ and $\mathbf{x}_d\in\mathbb{S}^n$, it follows from \eqref{definition:spherical_distance_function} that $d_s(\mathbf{x}, \mathbf{x}_d)\in[0, 2]$. 
Moreover, since $\epsilon < 1-\cos(\delta)$, Lemma \ref{lemma:constraint_separation} implies that if $\mathbf{x}\in\mathcal{N}_{\epsilon}(\mathcal{U}_i)$ for some $i\in\mathbb{I}$, then $\mathbf{x}\notin\mathcal{N}_{\epsilon}(\mathcal{U}_j)$ for all $j\in\mathbb{I}\setminus\{i\}$.
Therefore, by definition of $h_i(\cdot)$ in \eqref{obstacle-function-defnition}, the scalar function $\beta(\mathbf{x})$ in \eqref{navigation_function_on_sphere} is single-valued for all $\mathbf{x}\in\mathcal{M}_0$, and satisfies $\beta(\mathbf{x})\in[0, 1]$.
As a result, $W(\mathbf{x})$ is well-defined over $\mathcal{M}_0$.
Consequently, by \eqref{navigation_function_on_sphere}, Claims \ref{claim1:properties} and \ref{claim2:properties} of Lemma \ref{lemma:properties_scalar_function} follow.

We proceed to prove Claim \ref{claim3:properties} of Lemma \ref{lemma:properties_scalar_function}.
The scalar function $d_s(\mathbf{x}, \mathbf{x}_d)$ is smooth with respect to $\mathbf{x}$. 
Moreover, since $h_i$ is real analytic, the function $h_i(d_s(\mathbf{x}, \mathbf{g}_i))$ is twice continuously differentiable on $\mathcal{N}_{\epsilon}(\mathcal{U}_i)$ for each $i\in\mathbb{I}$.
In addition, since $\epsilon < 1-\cos(\delta)$, according to Lemma \ref{lemma:constraint_separation}, $\mathcal{N}_{\epsilon}(\mathcal{U}_i)\cap\mathcal{N}_{\epsilon}(\mathcal{U}_j) = \emptyset$ for all $i, j\in\mathbb{I}, i\ne j$.
Consequently, $\beta(\mathbf{x})$, used in \eqref{navigation_function_on_sphere}, is single-valued for all $\mathbf{x}\in\mathcal{M}_0$.
Furthermore, for $\mathbf{x}\in\mathcal{N}_{\epsilon}(\mathcal{U}_i)\cap\mathcal{M}_{\epsilon}$, one has $d_s(\mathbf{x}, \mathbf{g}_i) = \overline{\Delta}_i$.
Therefore, conditions $h_i(\overline{\Delta}_i) = 1, h_i'(\overline{\Delta}_i)=0$, and $h_i''(\overline{\Delta}_i) = 0$, stated after \eqref{obstacle-function-defnition}, ensure that $\beta$ is twice continuously differentiable over $\mathcal{M}_0$.
Therefore, $W(\mathbf{x})$ is twice continuously differentiable on $\mathcal{M}_0$.

\subsection{Proof of Theorem \ref{maintheorem:conic_constraints}}\label{proof:maintheorem:conic_constraints}
\subsubsection{Proof of Claim \ref{claim:conic1}}
According to Lemma \ref{lemma:properties_scalar_function}, the scalar function $W(\mathbf{x})$ defined in \eqref{navigation_function_on_sphere} is twice continuously differentiable on $\mathcal{M}_0$.
Taking the time derivative of $W(\mathbf{x})$, one obtains \begin{equation}\label{time_derivative_of_W}\dot{W}(\mathbf{x}) = -\nabla_{\mathbf{x}}W(\mathbf{x})^\top\mathbf{P}(\mathbf{x})\nabla_{\mathbf{x}}W(\mathbf{x}).\end{equation} 
Since $\mathbf{P}(\mathbf{x})$ is a positive semidefinite matrix for all $\mathbf{x}\in\mathbb{S}^n$, one has $\dot{W}(\mathbf{x}) \leq 0$ over $\mathcal{M}_0$.
In other words, $W(\mathbf{x}(t)) \leq W(\mathbf{x}(0))$ for all $t\geq 0$.
Additionally, by virtue of Lemma \ref{lemma:properties_scalar_function}, $W(\mathbf{x})$ attains its maximum value $k_1$ if and only if $\mathbf{x}\in\partial\mathcal{M}_0$.
Moreover, if $\mathbf{x}\in\partial\mathcal{M}_0$, then by Assumption \ref{assumption:non-overlapping-constraints}, there exists a unique $i\in\mathbb{I}$ such that $\mathbf{x}\in\partial\mathcal{U}_i$, and the control input becomes $\mathbf{u}(\mathbf{x}) = -\frac{k_1h_i'(\underline{\Delta}_i)}{d_s(\mathbf{x}, \mathbf{x}_d)}\mathbf{g}_i$.
Since $h_i'(\underline{\Delta}_i)>0$, the resulting projected vector field steers $\mathbf{x}$ toward $-\mathbf{g}_i$ and into the interior of $\mathcal{M}_0$.
Therefore $\mathcal{M}_0$ is forward invariant for the closed-loop system \eqref{system_dynamics}-\eqref{negative_gradient_control_law}.

\subsubsection{Proof of Claim \ref{claim:conic2}}
According to Lemma \ref{lemma:properties_scalar_function}, $W(\mathbf{x})$ is positive definite with respect to $\mathbf{x}_d$ over $\mathcal{M}_0$, and $\dot{W}(\mathbf{x}) \le 0$ for all $\mathbf{x}\in\mathcal{M}_0$, as established earlier in the proof of Claim \ref{claim:conic1} of Theorem \ref{maintheorem:conic_constraints}.
Therefore, $\mathbf{x}_d$ is a stable equilibrium of the closed-loop system \eqref{system_dynamics}–\eqref{negative_gradient_control_law}.
To establish almost global asymptotic stability of $\mathbf{x}_d$ on $\mathcal{M}_0$, it remains to show that $\mathbf{x}_d$ is almost globally attractive on $\mathcal{M}_0$.

It follows from LaSalle's invariance principle \cite[Theorem 4.4]{khalil2002nonlinear} that $\mathbf{x}$ will converge to the largest invariant set characterized by $\dot{W}(\mathbf{x}) = 0$.
It follows from \eqref{negative_gradient_control_law} and \eqref{time_derivative_of_W} that $\dot{W}(\mathbf{x}) = 0$ if and only if $\mathbf{P}(\mathbf{x})\mathbf{u}(\mathbf{x}) = \mathbf{0}$, where $\mathbf{P}(\mathbf{x})$ is defined in \eqref{orthogonal_projection_operator_formula}. Therefore the set characterized by $\dot{W}(\mathbf{x}) = 0$ is given by \begin{equation}\label{largest_invariant_set_defnition}
    \mathcal{E}:= \left\{\mathbf{x}\in\mathcal{M}_0\mid  \mathbf{P}(\mathbf{x} )\mathbf{u}(\mathbf{x}) = \mathbf{0} \right\}.
\end{equation}

We proceed to identify the elements of $\mathcal{E}$. 
First, we show that $\mathbf{u}(\mathbf{x})\ne \mathbf{0}$ for all $\mathbf{x}\in\mathcal{M}_0$.
It follows from \eqref{conic_control_expression}, \eqref{conic_individual_control_expression}, and the definition of $\beta(\cdot)$ in \eqref{obstacle-function-defnition} that, for any $\mathbf{x}\in\mathcal{M}_0$, the control input $\mathbf{u}(\mathbf{x})$ can be expressed as a linear combination of at most two unit vectors, namely $\mathbf{x}_d$ and $-\mathbf{g}_i$ for some $i\in\mathbb{I}$, with non-negative coefficients, at least one of which is strictly positive.
Moreover, since $\mathbf{x}_d\in\mathcal{M}_0^{\circ}$ and $\mathbf{g}_i\in\mathcal{U}_i^{\circ}$ for any $i\in\mathbb{I}$, it follows that $\mathbf{x}_d\ne c\mathbf{g}_i$ for any $c\geq 0$, any $\mathbf{x}\in\mathcal{M}_0$ and any $i\in\mathbb{I}$.
Therefore, $\mathbf{u}(\mathbf{x})\ne\mathbf{0}$ for all $\mathbf{x}\in\mathcal{M}_0$.
Consequently, the set $\mathcal{E}$ can be represented as follows:
\begin{equation}\label{alternate_expression_for_equilibrium_set}
    \mathcal{E} = \{\mathbf{x}\in\mathcal{M}_0\mid \mathbf{u}(\mathbf{x}) = c\mathbf{x} \text{ for some }c\in\mathbb{R}\setminus\{0\}\}.
\end{equation}

Since $\mathbf{u}(\mathbf{x}_d) = k_1\mathbf{x}_d$, it follows that $\mathbf{x}_d\in\mathcal{E}$.
Furthermore, $\mathbf{u}(-\mathbf{x}_d) = \frac{k_1}{9}\mathbf{x}_d$ if and only if $-\mathbf{x}_d\in\mathcal{M}_{\epsilon}$, where $\mathcal{M}_{\epsilon}$ is obtained by replacing $p$ with $\epsilon$ in \eqref{eroded_free_space_definition}.
Therefore, $\left(\{-\mathbf{x}_d\}\cap\mathcal{M}_{\epsilon}\right)\in\mathcal{E}$.
Since $\mathbf{u}(\mathbf{x}) = \frac{k_1}{(1 + d_s(\mathbf{x}, \mathbf{x}_d))^2}\mathbf{x}_d$ for all $\mathbf{x}\in\mathcal{M}_{\epsilon}$, there are no equilibrium points of the closed-loop system \eqref{system_dynamics}-\eqref{negative_gradient_control_law} in $\mathcal{M}_{\epsilon}\setminus\{\mathbf{x}_d, -\mathbf{x}_d\}$. 
Furthermore, since $h_i(d_s(\mathbf{x}, \mathbf{g}_i))=0$, $h_i'(d_s(\mathbf{x}, \mathbf{g}_i))>0$, and $\mathbf{P}(\mathbf{x})\mathbf{g}_i\ne \mathbf{0}$ for every $i\in\mathbb{I}$ and $\mathbf{x}\in\partial\mathcal{U}_i$, \eqref{conic_control_expression} implies that the closed-loop system \eqref{system_dynamics}-\eqref{negative_gradient_control_law} has no equilibrium point on $\partial\mathcal{U}$.
Consequently, if there exists $\mathbf{x}^*\in\mathcal{E}\setminus\{\mathbf{x}_d, -\mathbf{x}_d\}$, then, by Lemma \ref{lemma:constraint_separation}, there exists $i\in\mathbb{I}$ such that $\mathbf{x}^*\in\left(\mathcal{N}_{\epsilon}(\mathcal{U}_i)\right)^{\circ}$.
Suppose there exists an equilibrium point $\mathbf{x}^*$ of the closed-loop system \eqref{system_dynamics}-\eqref{negative_gradient_control_law} in $\left(\mathcal{N}_{\epsilon}(\mathcal{U}_i)\right)^{\circ}$ for some $i\in\mathbb{I}$.
We use the following fact to show that $\mathbf{x}^*$ 
belongs to the convex cone  $\mathcal{C}(-\mathbf{x}_d, \mathbf{g}_i)$ connecting $-\mathbf{x}_d$ and $\mathbf{g}_i$, where $\mathbf{g}_i$ is defined in \eqref{definition:conic_constraints} and the convex cone $\mathcal{C}(-\mathbf{x}_d, \mathbf{g}_i)$ is defined as in Section \ref{notations}.
\begin{fact}\label{fact:undesired}
    For the closed-loop system \eqref{system_dynamics}-\eqref{negative_gradient_control_law} under Assumption \ref{assumption:non-overlapping-constraints}, if 
    $\mathbf{x}^*\in\mathcal{E}\cap\left(\mathcal{N}_{\epsilon}(\mathcal{U}_i)\right)^{\circ}$ for some $i\in\mathbb{I}$, then
    $\mathbf{x}^*\in\mathcal{N}_{\mathcal{E}}^i\setminus\{-\mathbf{x}_d\}$ and $\mathbf{x}^*$ is an isolated equilibrium point, where
    \begin{equation}\label{expression_for_NE}
    \mathcal{N}_{\mathcal{E}}^i := \left(\mathcal{N}_{\epsilon}(\mathcal{U}_i)\right)^{\circ}\cap\mathcal{C}(-\mathbf{x}_d, \mathbf{g}_i).
    \end{equation}
\end{fact}
\proof{See Appendix \ref{proof:fact:undesired}.}

As a result, the equilibrium set $\mathcal{E}$ satisfies
\begin{equation}\label{equilibrium_point_set_inclusion}
    \mathcal{E}\subset \{\mathbf{x}_d\}\cup\left(\{-\mathbf{x}_d\}\cap\mathcal{M}_{\epsilon}\right)\cup\bigcup_{i\in\mathbb{I}}\mathcal{N}_{\mathcal{E}}^i.
\end{equation}
According to \eqref{largest_invariant_set_defnition} and \eqref{equilibrium_point_set_inclusion}, the desired point $\mathbf{x}_d$ is an equilibrium point of the closed-loop system \eqref{system_dynamics}-\eqref{negative_gradient_control_law}.
Additionally, if $-\mathbf{x}_d\in\mathcal{M}_{\epsilon}$, then $-\mathbf{x}_d$ is also an equilibrium point.
Furthermore, if there exists $\mathbf{x}^*\in\mathcal{E}\setminus\{\mathbf{x}_d, -\mathbf{x}_d\}$, then it belongs to a set $\mathcal{N}_{\mathcal{E}}^i$ for some $i\in\mathbb{I}$, and is an isolated equilibrium point.

To guarantee almost global attractivity of $\mathbf{x}_d$ for the closed-loop system \eqref{system_dynamics}-\eqref{negative_gradient_control_law} over $\mathcal{M}_0$, it is sufficient to show that the stable manifold associated with the undesired equilibrium set $\mathcal{E}\setminus\{\mathbf{x}_d\}$ has zero Lebesgue measure on $\mathbb{S}^n$.

The Jacobian of the closed-loop system \eqref{system_dynamics}-\eqref{negative_gradient_control_law} is given by
\begin{equation}\label{jacobian_expression_navigation}
    \mathbf{J}(\mathbf{x}) = \mathbf{P}(\mathbf{x})\frac{\partial\mathbf{u}(\mathbf{x})}{\partial\mathbf{x}} - \mathbf{x}\mathbf{u}(\mathbf{x})^\top - \mathbf{x}^\top\mathbf{u}(\mathbf{x})\mathbf{I}_{n+1},
\end{equation}
where using the control input representation in \eqref{conic_control_expression}, the matrix $\frac{\partial\mathbf{u}(\mathbf{x})}{\partial\mathbf{x}}$ is evaluated as
{\begin{equation*}
    \frac{\partial\mathbf{u}(\mathbf{x})}{\partial\mathbf{x}} = \begin{cases}
        \mathbf{u}_i^c(\mathbf{x})\nabla_{\mathbf{x}}\Psi_i(\mathbf{x})^\top + \Psi_i(\mathbf{x})\frac{\partial\mathbf{u}_i^c(\mathbf{x})}{\partial\mathbf{x}}, &\mathbf{x}\in\mathcal{N}_{\epsilon}(\mathcal{U}_i),\\
        \frac{2k_1}{(1 + d_s(\mathbf{x}, \mathbf{x}_d))^3}\mathbf{x}_d\mathbf{x}_d^\top, &\mathbf{x}\notin\mathcal{N}_{\epsilon}(\mathcal{U}),
    \end{cases}
\end{equation*}}
where $\Psi_i(\mathbf{x}) = \frac{k_1}{(d_s(\mathbf{x}, \mathbf{x}_d) + h_i(d_s(\mathbf{x}, \mathbf{g}_i)))^2}$, $\mathbf{u}_i^c(\mathbf{x})$ is given in \eqref{conic_individual_control_expression}, and 
\begin{equation*}
\begin{aligned}
    \frac{\partial\mathbf{u}_i^c(\mathbf{x})}{\partial\mathbf{x}} &=  h_i'(d_s(\mathbf{x}, \mathbf{g}_i))\left(\mathbf{g}_i\mathbf{x}_d^\top - \mathbf{x}_d\mathbf{g}_i^\top\right)\\
    &+ h_i''(d_s(\mathbf{x}, \mathbf{g}_i))d_s(\mathbf{x}, \mathbf{x}_d)\mathbf{g}_i\mathbf{g}_i^\top.
    \end{aligned}
\end{equation*}

Since $\epsilon < \bar{\epsilon}$, as stated in Section \ref{section:conic_constraint_avoidance}, one has $\mathbf{x}_d\notin\mathcal{N}_{\epsilon}(\mathcal{U})$.
Therefore, by \eqref{jacobian_expression_navigation}, the Jacobian matrix $\mathbf{J}(\mathbf{x}_d)$ for the closed-loop system \eqref{system_dynamics}-\eqref{negative_gradient_control_law} evaluated at $\mathbf{x}_d$ is given by
\begin{equation*}
    \mathbf{J}(\mathbf{x}_d) = -k_1\left(\mathbf{I}_{n+1} + \mathbf{x}_d\mathbf{x}_d^\top \right).
\end{equation*}
The matrix $\mathbf{J}(\mathbf{x}_d)$ has one eigenvalue equal to $-2k_1$ and an eigenvalue $-k_1$ with algebraic multiplicity $n$. Since all eigenvalues of $\mathbf{J}(\mathbf{x}_d)$ are negative, it follows that $\mathbf{x}_d$ is asymptotically stable for the closed-loop system \eqref{system_dynamics}–\eqref{negative_gradient_control_law}.

Next, we show that if $-\mathbf{x}_d\in\mathcal{M}_{\epsilon}$, then $-\mathbf{x}_d$ is an unstable node for the closed-loop system \eqref{system_dynamics}-\eqref{negative_gradient_control_law}.
If $-\mathbf{x}_d\in\mathcal{M}_{\epsilon}^{\circ}$, then there exists $\varrho > 0$ such that $\mathcal{B}_g(-\mathbf{x}_d, \varrho)\subset\mathcal{M}_{\epsilon}$ and $\mathbf{u}(\mathbf{x}) = \frac{k_1}{(1 + d_s(\mathbf{x}, \mathbf{x}_d))^2}\mathbf{x}_d$ for all $\mathbf{x}\in\mathcal{B}_g(-\mathbf{x}_d, \varrho)$, where $\mathcal{B}_g(-\mathbf{x}_d, \varrho)$ is defined as
\begin{equation*}\label{defiinition:geodesic_ball}\mathcal{B}_g(-\mathbf{x}_d, \varrho) = \{\mathbf{x}\in\mathbb{S}^n\mid d_s(\mathbf{x}, -\mathbf{x}_d)\leq \varrho\}.
\end{equation*}
Therefore, using \eqref{jacobian_expression_navigation}, the Jacobian matrix $\mathbf{J}(-\mathbf{x}_d)$ for the closed-loop system \eqref{system_dynamics}-\eqref{negative_gradient_control_law} evaluated at $-\mathbf{x}_d\in\mathcal{M}_{\epsilon}^{\circ}$ is given by
\begin{equation}\label{jacobian_at_minus_xd}
    \mathbf{J}(-\mathbf{x}_d) = \frac{k_1}{9}\left(\mathbf{I}_{n+1} + \mathbf{x}_d\mathbf{x}_d^\top\right).
\end{equation}
The matrix $\mathbf{J}(-\mathbf{x}_d)$ has one eigenvalue equal to $\frac{2k_1}{9}$ and an eigenvalue $\frac{k_1}{9}$ with algebraic multiplicity $n$. Since all eigenvalues of $\mathbf{J}(-\mathbf{x}_d)$ are positive, it follows that if $-\mathbf{x}_d\in\mathcal{M}_{\epsilon}^{\circ}$, then $-\mathbf{x}_d$ is an unstable node for the closed-loop system \eqref{system_dynamics}–\eqref{negative_gradient_control_law}.

Now, we consider the case where $-\mathbf{x}_d\in\partial\mathcal{M}_{\epsilon}$.
Since $\epsilon < 1-\cos(\delta)$, it follows from Lemma \ref{lemma:constraint_separation} that there exists a unique $i\in\mathbb{I}$ such that $-\mathbf{x}_d\in\partial\mathcal{N}_{\epsilon}(\mathcal{U}_i)\cap\mathcal{M}_{\epsilon}$.
Since $d_s(-\mathbf{x}_d, \mathcal{U}_i) = \epsilon$, it follows from the definition of the scalar function $h_i(\cdot)$ in \eqref{obstacle-function-defnition} that $h_i(d_s(-\mathbf{x}_d, \mathbf{g}_i)) = 1$, $h_i'(d_s(-\mathbf{x}_d, \mathbf{g}_i)) = 0$ and $h_i''(d_s(-\mathbf{x}_d, \mathbf{g}_i)) = 0$.
Consequently, if $-\mathbf{x}_d\in\partial\mathcal{M}_{\epsilon}$, then $\mathbf{J}(-\mathbf{x}_d)$ is given by \eqref{jacobian_at_minus_xd}, thereby ensuring that if $-\mathbf{x}_d\in\partial\mathcal{M}_{\epsilon}$, then $-\mathbf{x}_d$ is an unstable node for the closed-loop system \eqref{system_dynamics}-\eqref{negative_gradient_control_law}.

Finally, we show that if there exists $\mathbf{x}^*\in\left(\mathcal{E}\cap\mathcal{N}_{\mathcal{E}}^i\right)\setminus\{\mathbf{x}_d, -\mathbf{x}_d\}$ for some $i\in\mathbb{I}$, then the Jacobian matrix $\mathbf{J}(\mathbf{x}^*)$ for the closed-loop system \eqref{system_dynamics}-\eqref{negative_gradient_control_law} has at least $n-1$ eigenvalues with positive real parts, with corresponding eigenvectors lying in the tangent space $\mathbf{T}_{\mathbf{x}^*}(\mathbb{S}^n)$, where $\mathcal{N}_{\mathcal{E}}^i$ is defined in \eqref{expression_for_NE}.
Consequently, by the center manifold theorem \cite[Section 2.7, Pg. 116]{perko2013differential}, this will imply that the stable manifold at $\mathbf{x}^*$ on $\mathbb{S}^n$ is at most one-dimensional and has zero Lebesgue measure in $\mathbb{S}^n$.

Let $\mathbf{x}^*\in\left(\mathcal{E}\cap\mathcal{N}_{\mathcal{E}}^i\right)\setminus\{\mathbf{x}_d, -\mathbf{x}_d\}$ for some $i\in\mathbb{I}$. 
Since $\mathbf{P}(\mathbf{x}^*)\mathbf{u}(\mathbf{x}^*) = \mathbf{0}$, it follows from \eqref{conic_control_expression} and \eqref{jacobian_expression_navigation} that the Jacobian matrix $\mathbf{J}(\mathbf{x}^*)$ is given by
\begin{equation*}
    \begin{aligned}
        \mathbf{J}(\mathbf{x}^*) &=  \Psi_i(\mathbf{x}^*)\Big(h_i'(d_s(\mathbf{x}^*, \mathbf{g}_i))\left(\mathbf{g}_i\mathbf{x}_d^\top - \mathbf{x}_d\mathbf{g}_i^\top\right)\\
        & - h_i(d_s(\mathbf{x}^*, \mathbf{g}_i))\mathbf{x}^*\mathbf{x}_d^\top + h_i'(d_s(\mathbf{x}^*, \mathbf{g}_i))d_s(\mathbf{x}^*, \mathbf{x}_d)\mathbf{x}^*\mathbf{g}_i^\top\\
        & + h_i''(d_s(\mathbf{x}^*, \mathbf{g}_i))d_s(\mathbf{x}^*, \mathbf{x}_d)\mathbf{P}(\mathbf{x}^*)\mathbf{g}_i\mathbf{g}_i^\top\\
        &-\mathbf{x}^{*\top}\mathbf{u}_i^c(\mathbf{x}^*)\mathbf{I}_{n+1}\Big),
    \end{aligned}
\end{equation*}
where, for each $i\in\mathbb{I}$,
$\Psi_i(\mathbf{x}^*) = \frac{k_1}{(d_s(\mathbf{x}^*, \mathbf{x}_d) + h_i(d_s(\mathbf{x}^*, \mathbf{g}_i)))^2}$.

Define a $(n-1)$-dimensional subspace $\mathbb{U}_i(\mathbf{x}^*)$ of the $n$-dimensional tangent space $\mathbf{T}_{\mathbf{x}^*}(\mathbb{S}^n)$ as follows:
\[\mathbb{U}_i(\mathbf{x}^*) = \big\{\mathbf{p}\in\mathbf{T}_{\mathbf{x}^*}(\mathbb{S}^n)\mid\mathbf{p}^\top\mathbf{g}_i = 0\big\}.\]
Let $\boldsymbol{\eta}\in\mathbb{U}_i(\mathbf{x}^*)\setminus\{\mathbf{0}\}$ and evaluate $\mathbf{J}(\mathbf{x}^*)\boldsymbol{\eta}$ as follows:
\begin{equation}\label{intermediate_expression_for_eigenvalue}
    \begin{aligned}
        \mathbf{J}(\mathbf{x}^*)\boldsymbol{\eta} &=  \Psi_i(\mathbf{x}^*)\Big(h_i'(d_s(\mathbf{x}^*, \mathbf{g}_i))\mathbf{P}(\mathbf{x}^*)\mathbf{g}_i\mathbf{x}_d^\top\boldsymbol{\eta}\\
    & - h_i(d_s(\mathbf{x}^*, \mathbf{g}_i))\mathbf{x}^*\mathbf{x}_d^\top\boldsymbol{\eta} -\mathbf{x}^{*\top}\mathbf{u}_i^c(\mathbf{x}^*)\boldsymbol{\eta}\Big),
    \end{aligned}
\end{equation}
where we used the fact that $\boldsymbol{\eta}^\top\mathbf{g}_i = 0$.

We show that for any $\boldsymbol{\eta}\in\mathbb{U}_i(\mathbf{x}^*)\setminus\{\mathbf{0}\}$, $\mathbf{x}_d^\top\boldsymbol{\eta} = 0$.
Since $\mathbf{x}^*\in(\mathcal{E}\cap\mathcal{N}_{\mathcal{E}}^i)\setminus\{\mathbf{x}_d, -\mathbf{x}_d\}$, one has $\mathbf{x}^*\ne\mathbf{g}_i$ and by \eqref{expression_for_NE}, $\mathbf{x}_d$ lies in the span of $\mathbf{x}^*$ and $\mathbf{g}_i$.
Consequently, for any $\boldsymbol{\eta}\in\mathbb{U}_i(\mathbf{x}^*)\setminus\{\mathbf{0}\}$ the conditions $\boldsymbol{\eta}^\top\mathbf{x}^* = 0$ and $\boldsymbol{\eta}^\top\mathbf{g}_i = 0$ imply that $\boldsymbol{\eta}^\top\mathbf{x}_d = 0$, and
\eqref{intermediate_expression_for_eigenvalue} reduces to
\begin{equation}\label{final_expression_of_eigenvalue_analysis}
    \begin{aligned}
        \mathbf{J}(\mathbf{x}^*)\boldsymbol{\eta} &=  -\Psi_i(\mathbf{x}^*)\mathbf{x}^{*\top}\mathbf{u}_i^c(\mathbf{x}^*)\boldsymbol{\eta}.
    \end{aligned}
\end{equation}
Since \eqref{final_expression_of_eigenvalue_analysis} holds for all $\boldsymbol{\eta}\in\mathbb{U}_i(\mathbf{x}^*)\setminus\{\mathbf{0}\}$ and $\mathbb{U}_i(\mathbf{x}^*)$ is an $(n-1)$-dimensional subspace of the tangent space $\mathbf{T}_{\mathbf{x}^*}(\mathbb{S}^n)$, it follows that $-\Psi_i(\mathbf{x}^*)\mathbf{x}^{*\top}\mathbf{u}_i^c(\mathbf{x}^*)$ is an eigenvalue of $\mathbf{J}(\mathbf{x}^*)$ with geometric multiplicity $n-1$ and corresponding eigenspace $\mathbb{U}_i(\mathbf{x}^*)$.

We proceed to show that for any $\mathbf{x}^*\in(\mathcal{E}\cap\mathcal{N}_{\mathcal{E}}^i)\setminus\{\mathbf{x}_d, -\mathbf{x}_d\}$, the eigenvalue $-\Psi_i(\mathbf{x}^*)\mathbf{x}^{*\top}\mathbf{u}_i^c(\mathbf{x}^*)$ is positive, where $i\in\mathbb{I}$ and the set $\mathcal{N}_{\mathcal{E}}^i$ is defined in \eqref{expression_for_NE}. 
Since $\mathbf{x}^*\ne \mathbf{x}_d$, one has $d_s(\mathbf{x}^*, \mathbf{x}_d) > 0$.
Furthermore, $h_i(d_s(\mathbf{x}^*, \mathbf{g}_i))\geq 0$ since $\mathbf{x}^*\in\mathcal{N}_{\mathcal{E}}^i$. 
Therefore, since $k_1 > 0$, it follows that $\Psi_i(\mathbf{x}^*) > 0$, where for each $i\in\mathbb{I}$,
$\Psi_i(\mathbf{x}^*) = \frac{k_1}{(d_s(\mathbf{x}^*, \mathbf{x}_d) + h_i(d_s(\mathbf{x}^*, \mathbf{g}_i)))^2}$.
Therefore, to show that $-\Psi_i(\mathbf{x}^*)\mathbf{x}^{*\top}\mathbf{u}_i^c(\mathbf{x}^*) > 0$, it is sufficient to show that $\mathbf{x}^{*\top}\mathbf{u}_i^c(\mathbf{x}^*) < 0$.

First, we show that for any $\mathbf{x}^*\in\mathcal{N}_{\mathcal{E}}^i$, one has $\mathbf{x}^{*\top}\mathbf{u}_i^c(\mathbf{x}^*) < 0$.
Since $\mathbf{x}^*\in\mathcal{N}_{\mathcal{E}}^i$, it follows from \eqref{expression_for_NE} that $\mathbf{x}^*\in\mathcal{C}(-\mathbf{x}_d, \mathbf{g}_i)$, where the convex cone $\mathcal{C}(-\mathbf{x}_d, \mathbf{g}_i)$ is defined as in Section \ref{notations}.
Moreover, since $\mathbf{x}^*\in\mathcal{E}\cap\mathcal{N}_{\mathcal{E}}^i$, it follows from \eqref{conic_control_expression}, \eqref{conic_individual_control_expression} and \eqref{alternate_expression_for_equilibrium_set} that $\mathbf{u}_i^c(\mathbf{x}^*)\in\mathcal{C}(\mathbf{x}_d, -\mathbf{g}_i)\setminus\{\mathbf{0}\}$.
Furthermore, since $\mathbf{x}^*\in\mathcal{E}$, by \eqref{alternate_expression_for_equilibrium_set}, one has $\mathbf{u}_i^c(\mathbf{x}^*) = c\mathbf{x}^*$ for some $c\in\mathbb{R}\setminus\{0\}$.
Consequently, since the vectors $\mathbf{x}^*$ and $\mathbf{u}_i^c(\mathbf{x}^*)$ are linearly dependent and point in the opposite direction with respect to each other, it follows that $\mathbf{x}^{*\top}\mathbf{u}_i^c(\mathbf{x}^*) < 0$.
This completes the proof of Claim \ref{claim:conic2} of Theorem \ref{maintheorem:conic_constraints}.

\subsection{Proof of Fact \ref{fact:undesired}}\label{proof:fact:undesired}
Since $\mathbf{x}^*\in\mathcal{E}$, it follows from \eqref{alternate_expression_for_equilibrium_set} that $\mathbf{u}(\mathbf{x}^*) = c\mathbf{x}^*$ for some $c \in\mathbb{R}\setminus\{0\}$.
Additionally, according to \eqref{conic_control_expression} and \eqref{conic_individual_control_expression}, for any $\mathbf{x}\in\left(\mathcal{N}_{\epsilon}(\mathcal{U}_i)\right)^{\circ}$, the control input satisfies $\mathbf{u}(\mathbf{x})\in\mathcal{C}(\mathbf{x}_d, -\mathbf{g}_i)$, where the convex cone $\mathcal{C}(\mathbf{x}_d, -\mathbf{g}_i)$ is defined as in Section \ref{notations}.
Therefore, $\mathbf{x}^*$ satisfies
\[\mathbf{x}^*\in\left(\mathcal{N}_{\epsilon}(\mathcal{U}_i)\right)^{\circ}\cap\left(\mathcal{C}(-\mathbf{x}_d, \mathbf{g}_i)\cup\mathcal{C}(\mathbf{x}_d,-\mathbf{g}_i)\right).\]

We show that for any $i\in\mathbb{I}$, $\left(\mathcal{N}_{\epsilon}(\mathcal{U}_i)\right)^{\circ}\cap\mathcal{C}(\mathbf{x}_d, -\mathbf{g}_i) = \emptyset$.
Since $\epsilon < \bar{\epsilon}$, as defined in Section \ref{section:conic_constraint_avoidance}, one has $\mathbf{x}_d\notin\mathcal{D}_{\epsilon}(\mathcal{U})$.
Additionally, since $\mathcal{D}_{\epsilon}(\mathcal{U}_i) \ne \mathbb{S}^n$ for any $i\in\mathbb{I}$, one has $-\mathbf{g}_i\notin\mathcal{D}_{\epsilon}(\mathcal{U}_i)$ for each $i\in\mathbb{I}$.
Consequently, since $\mathcal{U}_i$ is a conic constraint \eqref{definition:conic_constraints}, one can verify that $\left(\mathcal{N}_{\epsilon}(\mathcal{U}_i)\right)^{\circ}\cap\mathcal{C}(\mathbf{x}_d, -\mathbf{g}_i) = \emptyset$ for any $i\in\mathbb{I}$.
Therefore, $\mathbf{x}^*\in\mathcal{N}_{\mathcal{E}}^i$.

Furthermore, suppose that $-\mathbf{x}_d\in\left(\mathcal{N}_{\epsilon}(\mathcal{U}_i)\right)^{\circ}$ for some $i\in\mathbb{I}$.
It follows from \eqref{conic_control_expression} and \eqref{conic_individual_control_expression} that 
\begin{equation}\mathbf{P}(-\mathbf{x}_d)\mathbf{u}(-\mathbf{x}_d) = \frac{-2k_1h_i'(d_s(-\mathbf{x}_d, \mathbf{g}_i))}{(2+h_i(d_s(-\mathbf{x}_d, \mathbf{g}_i)))^2}\mathbf{P}(-\mathbf{x}_d)\mathbf{g}_i.\end{equation} 
Since $d_s(-\mathbf{x}_d, \mathbf{g}_i)\in(\underline{\Delta}_i, \overline{\Delta}_i)$, $h_i'(d_s(-\mathbf{x}_d, \mathbf{g}_i)) > 0$ and $\mathbf{P}(-\mathbf{x}_d)\mathbf{g}_i\ne\mathbf{0}$, it follows that $-\mathbf{x}_d\in\left(\mathcal{N}_{\epsilon}(\mathcal{U}_i)\right)^{\circ}$ is not an equilibrium point.
Consequently, if $\mathbf{x}^*\in\mathcal{E}\cap\left(\mathcal{N}_{\epsilon}(\mathcal{U}_i)\right)^{\circ}$ for some $i\in\mathbb{I}$, then $\mathbf{x}^*\in\mathcal{N}_{\mathcal{E}}^i\setminus\{-\mathbf{x}_d\}$.

We proceed to show that $\mathbf{x}^*\in\mathcal{N}_{\mathcal{E}}^i\setminus\{-\mathbf{x}_d\}$ is an isolated equilibrium point of the closed-loop system \eqref{system_dynamics}-\eqref{negative_gradient_control_law}.
Let $\widehat{\mathcal{N}}_{\mathcal{E}}^i = \mathcal{N}_{\mathcal{E}}^i\setminus\{-\mathbf{x}_d\}$.
For $\mathbf{x}\in\widehat{\mathcal{N}}_{\mathcal{E}}^i$, define the scalar function $f_i(\mathbf{x}) = \mathbf{g}_i^\top\mathbf{P}(\mathbf{x})\mathbf{u}(\mathbf{x})$.
Since $\mathbf{P}(\mathbf{x})\mathbf{u}(\mathbf{x}) = \mathbf{0}$ implies that $f_i(\mathbf{x}) = 0$, it is sufficient to show that $\mathbf{x}^*$ is an isolated zero of $f_i$ on $\widehat{\mathcal{N}}_{\mathcal{E}}^i$.

Using \eqref{conic_control_expression}, for every $\mathbf{x}\in\widehat{\mathcal{N}}_{\mathcal{E}}^i$, one gets
\begin{equation}
    f_i(\mathbf{x}) = \Psi_i(\mathbf{x})f_i^c(\mathbf{x}), \qquad f_i^c(\mathbf{x}) :=\mathbf{g}_i^\top\mathbf{P}(\mathbf{x})\mathbf{u}_i^c(\mathbf{x}),
\end{equation}
where $\Psi_i(\mathbf{x}) = \frac{k_1}{(d_s(\mathbf{x}, \mathbf{x}_d) + h_i(d_s(\mathbf{x}, \mathbf{g}_i)))^2}$.
Since $\Psi_i(\mathbf{x}) > 0$ for all $\mathbf{x}\in\widehat{\mathcal{N}}_{\mathcal{E}}^i$, $f_i$ and $f_i^c$ have the same zeros on $\widehat{\mathcal{N}}_{\mathcal{E}}^i$.
Therefore, it is sufficient to show that $\mathbf{x}^*$ is an isolated zero of $f_i^c$.

Substituting \eqref{conic_individual_control_expression} into $f_i^c$, one obtains
\begin{equation}\label{function_f_i_c}
\begin{aligned}
    f_i^c(\mathbf{x}) &= h_i(d_s(\mathbf{x}, \mathbf{g}_i))\mathbf{g}_i^\top\mathbf{P}(\mathbf{x})\mathbf{x}_d\\ 
    &- h_i'(d_s(\mathbf{x}, \mathbf{g}_i))d_s(\mathbf{x}, \mathbf{x}_d)\|\mathbf{P}(\mathbf{x})\mathbf{g}_i\|^2.
    \end{aligned}
\end{equation}
Since $\widehat{\mathcal{N}}_{\mathcal{E}}^i\subset\mathcal{G}(\mathbf{g}_i, -\mathbf{x}_d)$ is a continuous, one-dimensional curve, there exist $0\leq a_i < b_i \leq 1$ such that
\begin{equation}
    \widehat{\mathcal{N}}_{\mathcal{E}}^i=\{\vartheta_i(\lambda)\mid\lambda\in I_i\}, \qquad \vartheta_i(\lambda) := g(\lambda;\mathbf{g}_i, -\mathbf{x}_d),
\end{equation}
where $g$ is defined as in \eqref{geodesic_parameterization}, and $I_i = (a_i, b_i)$ denotes the open interval parameterizing $\widehat{\mathcal{N}}_{\mathcal{E}}^i$.

Define $\widehat{f}_i^c:I_i\to\mathbb{R}$ by $\widehat{f}_i^c(\lambda) = f_i^c(\vartheta_i(\lambda))$.
Since $h_i$, defined in \eqref{obstacle-function-defnition}, is real analytic on $(\underline{\Delta}_i, \overline{\Delta}_i)$, and $\vartheta_i, d_s$, and $\mathbf{P}$ are analytic, it follows from \eqref{function_f_i_c} that $\widehat{f}_i^c$ is real analytic on $I_i$.
By \cite[Corollary 1.2.7]{krantz2002primer}, the zeros of $\widehat{f}_i^c$ are isolated provided that there exists $\lambda_0\in I_i$ such that $\widehat{f}_i^c(\lambda_0)\ne 0$.\footnote{In the application of \cite[Corollary 1.2.7]{krantz2002primer}, the two real analytic functions are $\widehat{f}_i^c$ and the zero function defined over the open interval $I_i$.}

Let $\underline{\mathbf{x}}_i := \vartheta_i(a_i)\in\partial\mathcal{U}_i$.
Since $d_s(\underline{\mathbf{x}}_i, \mathbf{g}_i) = \underline{\Delta}_i$, by definition of $h_i$ in Section \ref{section:conic_constraint_avoidance}, one has $h_i(\underline{\Delta}_i) = 0$ and $h_i'(\underline{\Delta}_i)> 0$.
Furthermore, since $\mathbf{x}_d\notin\partial\mathcal{U}$, one as $d_s(\underline{\mathbf{x}}_i, \mathbf{x}_d) > 0$.
Therefore, it follows from \eqref{function_f_i_c} that
\begin{equation}
    f_i^c(\underline{\mathbf{x}}_i) = -h_i'(\underline{\Delta}_i)d_s(\underline{\mathbf{x}}_i, \mathbf{x}_d)\|\mathbf{P}(\underline{\mathbf{x}}_i)\mathbf{g}_i\|^2 < 0.
\end{equation}
Since $\widehat f_i^c=f_i^c\circ\vartheta_i$ extends continuously to $\lambda = a_i$,
there exists $\lambda_0\in I_i$, close to $a_i$, such that $\widehat{f}_i^c(\lambda_0) < 0$.
Consequently, the zeros of $\widehat{f}_i^c$ are isolated.
Since $\mathbf{x}^*$ lies in $\widehat{\mathcal{N}}_{\mathcal{E}}^i$, there exists $\lambda^*\in I_i$ such that $\mathbf{x}^* = \vartheta_i(\lambda^*)$.
Therefore, $\mathbf{x}^*$ is an isolated zero of $f_i^c$, and hence of $f_i$.
It follows that $\mathbf{x}^*$ is an isolated equilibrium point of the closed-loop system \eqref{system_dynamics}-\eqref{negative_gradient_control_law}.
This completes the proof.

\subsection{Proof of Lemma \ref{lemma:n-sphere-forward-invariance}}\label{proof:lemma:n-sphere-forward-invariance}
By Assumption \ref{assumption:non-overlapping-constraints}, if $\mathbf{x}\in\partial\mathcal{M}_0$, then $\mathbf{x}\in\partial\mathcal{U}_i$ for some $i\in\mathbb{I}$ and $\mathbf{x}\notin\partial\mathcal{U}_j$ for all $j\in\mathbb{I}$ with $j\ne i$.
According to \eqref{individual_control_input_vector_design}, 
if $\mathbf{x}\in\partial\mathcal{U}_i$ for some $i\in\mathbb{I}$, then the control input vector \eqref{n-sphere-control-law} simplifies to $\mathbf{u}(\mathbf{x}) = \frac{-k_1}{\kappa}\mathbf{g}_i$, and the resulting vector field $\mathbf{P}(\mathbf{x})\mathbf{u}(\mathbf{x})$ steers $\mathbf{x}$ in the direction aligned with the geodesic $\mathcal{G}(\mathbf{x}, -\mathbf{g}_i)$ towards $-\mathbf{g}_i$.
Moreover, since for each $i\in\mathbb{I}$, $\mathcal{U}_i$ is a star-shaped set on $\mathbb{S}^n$ and  $\mathbf{g}_i\in\sigma(\mathcal{U}_i)\cap\mathcal{U}_i^{\circ}$, it follows from Lemma \ref{lemma:reverse_geodesic_always_stay_outside} that $\mathcal{G}(\mathbf{x}, -\mathbf{g}_i)\cap\mathcal{U}_i^{\circ} = \emptyset$ for all $\mathbf{x}\in\partial\mathcal{U}_i$.
Consequently, for each $i\in\mathbb{I}$ and for every $\mathbf{x}\in\partial\mathcal{U}_i$, the vector field $\mathbf{P}(\mathbf{x})\mathbf{u}(\mathbf{x})$ belongs to the tangent cone\footnote{The tangent cone to the set $\mathcal{M}_0\subset\mathbb{S}^n$ at $\mathbf{x}\in\mathbb{S}^n$, denoted by $\mathcal{T}_{\mathbf{x}}(\mathcal{M}_0)$, is defined as in \cite[Def. 4.1.1]{aubin2009setvalued}.} $\mathcal{T}_{\mathbf{x}}(\mathcal{M}_0)$ to $\mathcal{M}_0$ on $\mathbb{S}^n$ at $\mathbf{x}$, \textit{i.e.,} $\mathbf{P}(\mathbf{x})\mathbf{u}(\mathbf{x})\in\mathcal{T}_{\mathbf{x}}(\mathcal{M}_0)$ for all $\mathbf{x}\in\partial\mathcal{M}_0$.
Furthermore, since $\mathcal{M}_0\subset\mathbb{S}^n$, for all $\mathbf{x}\in\mathcal{M}_0^{\circ}$, one has $\mathcal{T}_{\mathbf{x}}(\mathcal{M}_0) = \mathbf{T}_{\mathbf{x}}(\mathbb{S}^n)$.
Therefore, for all $\mathbf{x}\in\mathcal{M}_0$, $\mathbf{P}(\mathbf{x})\mathbf{u}(\mathbf{x})\in\mathcal{T}_{\mathbf{x}}(\mathcal{M}_0)$.
As a result, by Nagumo's theorem \cite[Theorem 11.2.3]{aubin2011viability}, the free space $\mathcal{M}_0$ is forward invariant for the closed-loop system \eqref{system_dynamics}-\eqref{n-sphere-control-law}. This completes the proof.

\subsection{Proof of Lemma \ref{lemma:lipschitz_continuity}}\label{proof:lemma:lipschitz_continuity}
For any given $\mathbf{a}\in\mathbb{S}^n$, the scalar function $d_s(\mathbf{x}, \mathbf{a})$, defined as in \eqref{definition:spherical_distance_function}, is globally Lipschitz in $\mathbf{x}$ over $\mathbb{S}^n$ with Lipschitz constant $1$.
Since, for any $\mathbf{x}\in\mathbb{S}^n$ and any closed set $\mathcal{U}_i\subset\mathbb{S}^n$, where $i\in\mathbb{I}$, the scalar function $d_s(\mathbf{x}, \mathcal{U}_i)$
is the pointwise minimum of Lipschitz functions $d_s(\mathbf{x}, \mathbf{a})$ with $\mathbf{a}\in\mathcal{U}_i$, it follows that $d_s(\mathbf{x}, \mathcal{U}_i)$ is locally Lipschitz in $\mathbf{x}$ over $\mathbb{S}^n$ for every $i\in\mathbb{I}$. 
Since the control input vector $\mathbf{u}(\mathbf{x})$ \eqref{n-sphere-control-law} is obtained through addition and scalar multiplication of locally Lipschitz functions, it is locally Lipschitz in $\mathbf{x}$ over $\mathcal{M}_0$. 
Moreover, since $\mathbf{P}(\mathbf{x})$ is continuously differentiable for all $\mathbf{x}\in\mathbb{S}^n$, it follows that $\mathbf{P}(\mathbf{x})\mathbf{u}(\mathbf{x})$ is locally Lipschitz in $\mathbf{x}$ over $\mathcal{M}_0$, and the proof is complete.

\subsection{Proof of Lemma \ref{lemma:angle_keeps_changing}}\label{proof:lemma:angle_keeps_changing}
\subsubsection{Proof of Claim \ref{claim1:LemmaVi}}
Using the fact $\mathbf{P}(\mathbf{x})^2=\mathbf{P}(\mathbf{x})$ and $\mathbf{P}(\mathbf{x})\mathbf{x} = \mathbf{0}$ for all $\mathbf{x}\in\mathbb{S}^n$, the scalar function $V_i(\mathbf{x})$, defined in \eqref{the_only_positive_semidefinite_function}, can be re-written as:
\begin{equation}\label{scalar_function_Vi}
    V_i(\mathbf{x}) = \frac{\mathbf{x}_d^\top\mathbf{P}(\mathbf{g}_i)\mathbf{x}}{\|\mathbf{P}(\mathbf{g}_i)\mathbf{x}_d\|\|\mathbf{P}(\mathbf{g}_i)\mathbf{x}\|}.
\end{equation}
We know that $\mathbf{g}_i \ne -\mathbf{x}_d$ for every $i \in \mathbb{I}$, as stated in Section~\ref{section:feedback_control_design}. 
Moreover, since $\mathbf{x}_d \in \mathcal{M}_{0}^{\circ}$, it follows that $\mathbf{g}_i \ne \mathbf{x}_d$ for each $i \in \mathbb{I}$. 
Therefore, $\mathbf{P}(\mathbf{g}_i)\mathbf{x}_d \ne \mathbf{0}$. 
Furthermore, since $-\mathbf{g}_i \notin \mathcal{F}_i$, it follows that $\mathbf{P}(\mathbf{g}_i)\mathbf{x} \ne \mathbf{0}$ for all $\mathbf{x} \in \mathcal{F}_i$. 
Consequently, $V_i(\mathbf{x})$ is well-defined for all $\mathbf{x} \in \mathcal{F}_i$, where $\mathcal{F}_i$ is defined in Lemma \ref{lemma:angle_keeps_changing}.

\subsubsection{Proof of Claim \ref{claim2:LemmaVi}}
Taking the time derivative of $V_i(\mathbf{x})$ at $\mathbf{x}\in\mathcal{F}_i\setminus\left(\partial\mathcal{U}_i\cup\mathcal{Z}_i\cup\mathcal{V}_i\right)$, one obtains 
\begin{equation*}
    \dot{V}_i(\mathbf{x}) = \frac{\mathbf{x}_d^\top\mathbf{P}(\mathbf{g}_i)\dot{\mathbf{x}}}{\|\mathbf{P}(\mathbf{g}_i)\mathbf{x}_d\|\|\mathbf{P}(\mathbf{g}_i)\mathbf{x}\|} - \frac{\mathbf{x}_d^\top\mathbf{P}(\mathbf{g}_i)\mathbf{x}\mathbf{x}^\top\mathbf{P}(\mathbf{g}_i)\dot{\mathbf{x}}}{\|\mathbf{P}(\mathbf{g}_i)\mathbf{x}_d\|\|\mathbf{P}(\mathbf{g}_i)\mathbf{x}\|^3}.
\end{equation*}
Since $\mathbf{P}(\mathbf{g}_i)\mathbf{x}_d \ne \mathbf{0}$ and $\mathbf{P}(\mathbf{g}_i)\mathbf{x} \ne \mathbf{0}$ for all $\mathbf{x} \in \mathcal{F}_i$, as noted earlier, and since $\dot{\mathbf{x}}$ is well-defined on $\mathcal{M}_0$, it follows that $\dot{V}(\mathbf{x})$ is well-defined for all $\mathbf{x} \in \mathcal{F}_i$.

To show that $\dot{V}_i(\mathbf{x}) > 0$ for all $\mathbf{x}\in\mathcal{F}_i\setminus\left(\partial\mathcal{U}_i\cup\mathcal{Z}_i\cup\mathcal{V}_i\right)$, it is sufficient to show that
\begin{equation}\label{sufficient_condition_to_show}
    \mathbf{w}_i(\mathbf{x})^\top\mathbf{P}(\mathbf{x})\mathbf{u}(\mathbf{x}) > 0, \text{ for all }\mathbf{x}\in\mathcal{F}_i\setminus\left(\partial\mathcal{U}_i\cup\mathcal{Z}_i\cup\mathcal{V}_i\right),
\end{equation}
where $\mathbf{w}_i(\mathbf{x})$ is given by
\begin{equation}\label{vector_w_definition}
    \mathbf{w}_i(\mathbf{x}) = \|\mathbf{P}(\mathbf{g}_i)\mathbf{x}\|^2\mathbf{P}(\mathbf{g}_i)\mathbf{x}_d - \mathbf{x}_d^\top\mathbf{P}(\mathbf{g}_i)\mathbf{x}\mathbf{P}(\mathbf{g}_i)\mathbf{x}.
\end{equation}
To proceed with the proof, we require the following fact:
\begin{fact}\label{fact:angle_condition_help}
Let $\mathbf{w}_i(\mathbf{x})$ be defined as in \eqref{vector_w_definition} for $\mathbf{x}\in\mathcal{F}_i$. Then, the following hold:
\begin{enumerate}
\item\label{claim1:factAngle} $\mathbf{w}_i(\mathbf{x})^\top \mathbf{P}(\mathbf{x}) \mathbf{x}_d > 0$ for all $\mathbf{x} \in \mathcal{F}_i \setminus (\mathcal{Z}_i \cup \mathcal{V}_i)$,
\item\label{claim2:factAngle} $\mathbf{w}_i(\mathbf{x})^\top \mathbf{P}(\mathbf{x}) \mathbf{g}_i = 0$ for all $\mathbf{x} \in \mathcal{F}_i$.
\end{enumerate}
\end{fact}
\proof{Using \eqref{vector_w_definition}, one obtains
\begin{equation}\label{W3_term_different_definition}
    \mathbf{w}_i(\mathbf{x})^\top \mathbf{P}(\mathbf{x})\mathbf{x}_d = \mathbf{x}^\top\mathbf{P}(\mathbf{g}_i)(\mathbf{x}\mathbf{x}_d^\top-\mathbf{x}_d\mathbf{x}^\top)\mathbf{P}(\mathbf{g}_i)\mathbf{P}(\mathbf{x})\mathbf{x}_d.
\end{equation}
Using \eqref{orthogonal_projection_operator_formula}, one gets $\mathbf{P}(\mathbf{g}_i)\mathbf{P}(\mathbf{x})\mathbf{x}_d = \mathbf{P}(\mathbf{g}_i)\mathbf{x}_d - \mathbf{x}^\top\mathbf{x}_d\mathbf{P}(\mathbf{g}_i)\mathbf{x}$, and \eqref{W3_term_different_definition} becomes
\begin{equation}\label{expression_11}
\begin{aligned}
    \mathbf{w}_i(\mathbf{x})^\top \mathbf{P}(\mathbf{x})\mathbf{x}_d &= \mathbf{x}^\top\mathbf{P}(\mathbf{g}_i)(\mathbf{x}\mathbf{x}_d^\top-\mathbf{x}_d\mathbf{x}^\top)\mathbf{P}(\mathbf{g}_i)\mathbf{x}_d\\
    &-\mathbf{x}^\top\mathbf{x}_d\mathbf{x}^\top\mathbf{P}(\mathbf{g}_i)(\mathbf{x}\mathbf{x}_d^\top-\mathbf{x}_d\mathbf{x}^\top)\mathbf{P}(\mathbf{g}_i)\mathbf{x}.
    \end{aligned}
\end{equation}
Since the matrix $\mathbf{x}\mathbf{x}_d^\top-\mathbf{x}_d\mathbf{x}^\top$ is skew symmetric, the second term in \eqref{expression_11} vanishes, and one obtains
\begin{equation}
\begin{aligned}
    \mathbf{w}_i(\mathbf{x})^\top \mathbf{P}(\mathbf{x})\mathbf{x}_d &= \|\mathbf{P}(\mathbf{g}_i)\mathbf{x}\|^2\|\mathbf{P}(\mathbf{g}_i)\mathbf{x}_d\|^2\\&-\left(\left(\mathbf{P}(\mathbf{g}_i)\mathbf{x}\right)^\top\mathbf{P}(\mathbf{g}_i)\mathbf{x}_d\right)^2.
    \end{aligned}
\end{equation}
It follows from Cauchy-Schwarz inequality that $\mathbf{w}_i(\mathbf{x})^\top\mathbf{P}(\mathbf{x})\mathbf{x}_d \geq 0$ for all $\mathbf{x}\in\mathcal{F}_i\setminus\left(\mathcal{Z}_i\cup\mathcal{V}_i\right)$.
In fact, $\mathbf{w}_i(\mathbf{x})^\top\mathbf{P}(\mathbf{x})\mathbf{x}_d = 0$ if and only if $\mathbf{P}(\mathbf{g}_i)\mathbf{x} = q\mathbf{P}(\mathbf{g}_i)\mathbf{x}_d$ for some $q\in\mathbb{R}$.
It can be shown that for any $\mathbf{x}\in\mathcal{F}_i\setminus\left(\mathcal{Z}_i\cup\mathcal{V}_i\right)$, there does not exist $q\in\mathbb{R}$ such that $\mathbf{P}(\mathbf{g}_i)\mathbf{x} = q\mathbf{P}(\mathbf{g}_i)\mathbf{x}_d$.
Consequently, $\mathbf{w}_i(\mathbf{x})^\top\mathbf{P}(\mathbf{x})\mathbf{x}_d > 0$ for all $\mathbf{x}\in\mathcal{F}_i\setminus\left(\mathcal{Z}_i\cup\mathcal{V}_i\right)$.

Now, we show that $\mathbf{w}_i(\mathbf{x})^\top\mathbf{P}(\mathbf{x})\mathbf{g}_i = 0$ for all $\mathbf{x}\in\mathcal{F}_i$.
Using \eqref{orthogonal_projection_operator_formula} and \eqref{vector_w_definition}, one obtains
\begin{equation*}
    \mathbf{w}_i(\mathbf{x})^\top\mathbf{P}(\mathbf{x})\mathbf{g}_i =-\mathbf{x}^\top\mathbf{g}_i \mathbf{x}^\top\mathbf{P}(\mathbf{g}_i)(\mathbf{x}\mathbf{x}_d^\top-\mathbf{x}_d\mathbf{x}^\top)\mathbf{P}(\mathbf{g}_i)\mathbf{x}.
\end{equation*}
Since the matrix $\mathbf{x}\mathbf{x}_d^\top-\mathbf{x}_d\mathbf{x}^\top$ is skew symmetric, it follows that $\mathbf{w}_i(\mathbf{x})^\top\mathbf{P}(\mathbf{x})\mathbf{g}_i = 0$ for all $\mathbf{x}\in\mathcal{F}_i$.
This completes the proof of Fact \ref{fact:angle_condition_help}.}

For any $\mathbf{x}\in\mathcal{F}_i\setminus\left(\partial\mathcal{U}_i\cup\mathcal{Z}_i\cup\mathcal{V}_i\right)$, using \eqref{n-sphere-control-law} and \eqref{individual_control_input_vector_design}, one has
\begin{equation*}
    \mathbf{P}(\mathbf{x})\mathbf{u}(\mathbf{x}) = a_i(\mathbf{x})\mathbf{P}(\mathbf{x})\mathbf{x}_d - \frac{1}{\kappa}\left(1-a_i(\mathbf{x})\right)\mathbf{P}(\mathbf{x})\mathbf{g}_i,
\end{equation*}
for some $a_i(\mathbf{x})\in(0, 1]$.
Consequently, it follows from Fact \ref{fact:angle_condition_help} that $\mathbf{w}_i(\mathbf{x})^\top\mathbf{P}(\mathbf{x})\mathbf{u}(\mathbf{x}) > 0$ for all $\mathbf{x}\in\mathcal{F}_i\setminus\left(\partial\mathcal{U}_i\cup\mathcal{Z}_i\cup\mathcal{V}_i\right)$.
As a result, $\dot{V}_i(\mathbf{x}) > 0$ for all $\mathbf{x}\in\mathcal{F}_i\setminus\left(\partial\mathcal{U}_i\cup\mathcal{Z}_i\cup\mathcal{V}_i\right)$.

\subsubsection{Proof of Claim \ref{claim3:lemmaVi}}
To show $\dot{V}_i(\mathbf{x}) = 0$ for any $\mathbf{x}\in\mathcal{F}_i\cap\left(\partial\mathcal{U}_i\cup\mathcal{Z}_i\cup\mathcal{V}_i\right)$, it is sufficient to show that \[\mathbf{w}_i(\mathbf{x})^\top\mathbf{P}(\mathbf{x})\mathbf{u}(\mathbf{x}) = 0\; \text{for all}\; \mathbf{x}\in\mathcal{F}_i\cap\left(\partial\mathcal{U}_i\cup\mathcal{Z}_i\cup\mathcal{V}_i\right),\] where $\mathbf{w}_i(\mathbf{x})$ is defined in \eqref{vector_w_definition}.
According to \eqref{set_definition_Vi_Zi}, for every $\mathbf{x}\in\mathcal{Z}_i\cup\mathcal{V}_i$, the vectors $\mathbf{P}(\mathbf{x})\mathbf{x}_d$ and $\mathbf{P}(\mathbf{x})\mathbf{g}_i$ are linearly dependent.
Moreover, by \eqref{n-sphere-control-law} and \eqref{individual_control_input_vector_design}, for all $\mathbf{x}\in\partial\mathcal{U}_i$, the control input becomes $\mathbf{u}(\mathbf{x}) = \frac{-k_1}{\kappa}\mathbf{g}_i$.
Furthermore, by construction $\mathcal{F}_i\cap\{\mathbf{g}_i, -\mathbf{g}_i\} = \emptyset$.
Therefore, for any $\mathbf{x}\in\mathcal{F}_i\cap\left(\partial\mathcal{U}_i\cup\mathcal{Z}_i\cup\mathcal{V}_i\right)$, the projected control input satisfies $\mathbf{P}(\mathbf{x})\mathbf{u}(\mathbf{x}) = a_i(\mathbf{x})\mathbf{P}(\mathbf{x})\mathbf{g}_i$ for some $a_i(\mathbf{x})\in\mathbb{R}$.
Consequently, by Claim \ref{claim2:factAngle} of Fact \ref{fact:angle_condition_help}, $\mathbf{w}_i(\mathbf{x})^\top\mathbf{P}(\mathbf{x})\mathbf{u}(\mathbf{x}) = 0$ for all $\mathbf{x}\in\mathcal{F}_i\cap\left(\partial\mathcal{U}_i\cup\mathcal{Z}_i\cup\mathcal{V}_i\right)$.
As a result, $\dot{V}_i(\mathbf{x}) = 0$ for all $\mathbf{x}\in\mathcal{F}_i\cap\left(\partial\mathcal{U}_i\cup\mathcal{Z}_i\cup\mathcal{V}_i\right)$.

\subsection{Proof of Theorem \ref{main_theorem}}\label{proof:main_theorem}

The forward invariance of $\mathcal{M}_0$ for the closed-loop system \eqref{system_dynamics}-\eqref{n-sphere-control-law} is established in Lemma \ref{lemma:n-sphere-forward-invariance}.
To show that the desired equilibrium point $\mathbf{x}_d$ is almost globally asymptotically stable for the closed-loop system \eqref{system_dynamics}-\eqref{n-sphere-control-law} over $\mathcal{M}_0$, it suffices to show that $\mathbf{x}_d$ is asymptotically stable and almost globally attractive in $\mathcal{M}_0$.

We first show that $\mathbf{x}_d$ is asymptotically stable for the closed-loop system \eqref{system_dynamics}-\eqref{n-sphere-control-law}.
Since $\epsilon < \bar{\epsilon}$, as stated in Section \ref{section:feedback_control_design}, one has $\mathbf{x}_d\notin\mathcal{D}_{\epsilon}(\mathcal{U})$.
Therefore, there exists $\varrho > 0$ such that $\mathcal{B}_g(\mathbf{x}_d, \varrho)\subset\mathcal{M}_{\epsilon}$ and $-\mathbf{x}_d\notin\mathcal{B}_g(\mathbf{x}_d, \varrho)$, where $\mathcal{M}_{\epsilon}$ is obtained by replacing $p$ with $\epsilon$ in \eqref{eroded_free_space_definition}, and the set $\mathcal{B}_g(\mathbf{x}_d, \varrho)$ is defined as 
\[\mathcal{B}_g(\mathbf{x}_d, \varrho) = \{\mathbf{x}\in\mathbb{S}^n\mid d_s(\mathbf{x}, \mathbf{x}_d)\leq \varrho\}.\]
The function $d_s(\mathbf{x}, \mathbf{x}_d)$ is positive definite with respect to $\mathbf{x}_d$ over $\mathcal{M}_0$.
It follows from \eqref{system_dynamics} and \eqref{n-sphere-control-law} that \[\dot{d}_s(\mathbf{x}, \mathbf{x}_d) = -k_1\mathbf{x}_d^\top\mathbf{P}(\mathbf{x})\mathbf{x}_d,\] which is negative definite with respect to $\mathbf{x}_d$ over $\mathcal{B}_g(\mathbf{x}_d, \varrho)$.
Consequently $\mathbf{x}_d$ is asymptotically stable for the closed-loop system \eqref{system_dynamics}-\eqref{n-sphere-control-law}.

We proceed to show that there exists $\bar{\kappa} > 0$ such that if $\kappa > \bar{\kappa}$, then $\mathbf{x}_d$ is almost globally attractive for the closed-loop system \eqref{system_dynamics}-\eqref{n-sphere-control-law} over $\mathcal{M}_0$.
In other words, we show that there exists $\bar{\kappa} > 0$ such that if $\kappa > \bar{\kappa}$, then the solution $\mathbf{x}(t)$ to the closed-loop system \eqref{system_dynamics}-\eqref{n-sphere-control-law}, initialized at any $\mathbf{x}(0)\in\mathcal{M}_0$ outside a set of Lebesgue measure zero, satisfies
\begin{equation}\label{almost_attractivity_equation}
    \lim_{t\to\infty}d_s(\mathbf{x}(t), \mathbf{x}_d) = 0.
\end{equation}

Consider a solution $\mathbf{x}(t)$ to the closed-loop system \eqref{system_dynamics}-\eqref{n-sphere-control-law} with $\mathbf{x}(0)\in\mathcal{M}_0\setminus\{-\mathbf{x}_d\}$.
There are two possible cases: either $\mathbf{x}(t)\in\mathcal{M}_{\epsilon}\setminus\mathcal{R}$ for all $t\geq 0$ or there exists $t_1\geq 0$ such that $\mathbf{x}(t_1)\in\mathcal{R}$, where the set $\mathcal{R}$ is given by
\begin{equation}\label{collection_of_all_Ri}
    \mathcal{R} = \bigcup_{i\in\mathbb{I}}\mathcal{R}_i,
\end{equation}
and for $i\in\mathbb{I}$, the sets $\mathcal{R}_i$ are defined in \eqref{definition:Ri} and \eqref{definition:Ri_forspecialindex}.
First consider the former case, where $\mathbf{x}(0)\in\mathcal{M}_0\setminus\{-\mathbf{x}_d\}$ and $\mathbf{x}(t)\in\mathcal{M}_{\epsilon}\setminus\mathcal{R}$ for all $t\geq 0$.
Since, according to \eqref{n-sphere-control-law}, $\mathbf{u}(\mathbf{x}) = k_1\mathbf{x}_d$ for all $\mathbf{x}\in\mathcal{M}_{\epsilon}$, it follows that $\dot{d}_s(\mathbf{x}, \mathbf{x}_d) < 0$ for all $\mathbf{x}\in\mathcal{M}_{\epsilon}\setminus\{\mathbf{x}_d, -\mathbf{x}_d\}$, where $d_s(\mathbf{x}, \mathbf{x}_d)$ is defined as in \eqref{definition:spherical_distance_function}. 
Consequently, $\mathbf{x}(t)$ satisfies \eqref{almost_attractivity_equation}.
Now, we proceed to analyze the case where $\mathbf{x}(0)\in\mathcal{M}_0\setminus\{-\mathbf{x}_d\}$ and there exists $t_1\geq 0$ such that $\mathbf{x}(t_1)\in\mathcal{R}$.

It follows from Assumption \ref{assumption:sufficient_separation} that there exists a unique $i\in\mathbb{I}$ such that $\mathbf{x}(t_1)\in\mathcal{R}_i$.
Depending on whether $-\mathbf{x}_d\in\mathcal{D}_{\epsilon}(\mathcal{U}_i)$, there are two possibilities: either $i\in\mathbb{I}_a$ or $i\in\mathbb{I}\setminus\mathbb{I}_a$, where the set $\mathbb{I}_a$ is defined in \eqref{definition:Ia} and contains indices $i$ corresponding to the constraint sets $\mathcal{U}_i$ such that $-\mathbf{x}_d \notin \mathcal{D}_{\epsilon}(\mathcal{U}_i)$. 
Additionally, since $\mathcal{Z}_i\cap\mathcal{R}_i\ne \emptyset$, it follows that either $\mathbf{x}(t_1)\in\mathcal{R}_i\setminus\mathcal{Z}_i$ or $\mathbf{x}(t_1)\in\mathcal{R}_i\cap\mathcal{Z}_i$, where the set $\mathcal{Z}_i$ is defined in \eqref{set_definition_Vi_Zi}. 
In other words, if $\mathbf{x}(t_1)\in\mathcal{R}_i$ for some $i\in\mathbb{I}$, then
there are three mutually exclusive possibilities as follows:
\begin{enumerate}
    \item $\mathbf{x}(t_1)\in\mathcal{R}_i\setminus\mathcal{Z}_i$ for some $i\in\mathbb{I}_a$;
    \item $\mathbf{x}(t_1)\in\mathcal{R}_i\setminus\mathcal{Z}_i$ with $i\in\mathbb{I}\setminus\mathbb{I}_a$;
    \item $\mathbf{x}(t_1)\in\mathcal{R}_i\cap\mathcal{Z}_i$ for some $i\in\mathbb{I}$.
\end{enumerate}
We proceed to analyze the behaviour of the solution $\mathbf{x}(t)$ that satisfies $\mathbf{x}(t_1)\in\mathcal{R}_i\setminus\mathcal{Z}_i$ for some $i\in\mathbb{I}$.
We show that if $i\in\mathbb{I}\setminus\mathbb{I}_a$, then, for sufficiently large $\kappa$, $\mathbf{x}(t)$ reaches $\mathcal{M}_{0}\setminus\mathcal{R}_i$ at some time $t_2  > t_1$.
Additionally, we show that $\mathbf{x}(t)$ remains in $\mathcal{M}_0$ away from $\mathcal{R}_i$ for all $t\geq t_2$.

If $i\in\mathbb{I}_a$, then we show that for sufficiently large $\kappa$, $\mathbf{x}(t)$ enters the set $\mathcal{M}_0\setminus\mathcal{R}_i$ only via the set $\mathcal{P}(\mathbf{x}_d, \mathcal{R}_i)$, where the set $\mathcal{P}(\mathbf{x}_d, \mathcal{R}_i)$ is defined as in \eqref{projection_set_definition} and contains points $\mathbf{x}$ in $\mathcal{R}_i$ that satisfy $d_s(\mathbf{x}_d, \mathbf{x}) = d_s(\mathbf{x}_d, \mathcal{R}_i)$.
Note that since $\mathbf{x}_d$ does not belong to $\mathcal{R}_i$, the set $\mathcal{P}(\mathbf{x}_d, \mathcal{R}_i)$ is a subset of the boundary of $\mathcal{R}_i$.
These statements are formally established in the following lemma:
\begin{lemma}\label{lemma:behaviour_of_x_in_Ri} Consider the closed-loop system \eqref{system_dynamics}--\eqref{n-sphere-control-law} under Assumptions \ref{assumption:non-overlapping-constraints} and \ref{assumption:sufficient_separation}. Suppose a solution $\mathbf{x}(t)$ satisfies $\mathbf{x}(t_1)\in\mathcal{R}_i\setminus\mathcal{Z}_i$ for some $i\in\mathbb{I}$ and for some $t_1\geq 0$. Then, there exist $\kappa_i>0$ such that for all $\kappa>\kappa_i$, the following statements hold: 
\begin{enumerate} 
\item \label{lemma:claim1}If $i\in\mathbb{I}_a$, then there exist $t_2\geq t_1$ and $\gamma > 0$ such that $\mathbf{x}(t_2)\in\mathcal{P}(\mathbf{x}_d, \mathcal{R}_i)$, $\mathbf{x}(t)\in\mathcal{R}_i$ for all $t\in[t_1, t_2]$, and $\mathbf{x}(t)\in\mathcal{M}_0\setminus\mathcal{R}_i$ for all $t\in(t_2, t_2 + \gamma]$.
\item \label{lemma:claim2}If $i\in\mathbb{I}\setminus\mathbb{I}_a$, then there exists $t_2> t_1$ such that $\mathbf{x}(t_2)\in\mathcal{M}_{0}\setminus\mathcal{R}_i$ and $\mathbf{x}(t)\notin\mathcal{R}_i$ for all $t\geq t_2$,\end{enumerate}
where the set $\mathcal{P}(\mathbf{x}_d, \mathcal{R}_i)$ is defined as in \eqref{projection_set_definition}. \end{lemma}
\proof{See Appendix \ref{proof:lemma:behavious_of_x_in_Ri}.}

Since the index set $\mathbb{I}$ is finite, define $\bar{\kappa}:=\max_{i\in\mathbb{I}}\kappa_i$, where the existence of each $\kappa_i$ is guaranteed by Lemma \ref{lemma:behaviour_of_x_in_Ri}.
If $\mathbf{x}(t_1)\in\mathcal{R}_i\setminus\mathcal{Z}_i$ for $i\in\mathbb{I}\setminus\mathbb{I}_a$, then it follows from Claim \ref{lemma:claim2} of Lemma \ref{lemma:behaviour_of_x_in_Ri} that for $\kappa > \bar{\kappa}$, $\mathbf{x}(t_2)\in\mathcal{M}_{0}\setminus\mathcal{R}_i$ for some $t_2 > t_1$ and $\mathbf{x}(t)\notin\mathcal{R}_i$ for all $t\geq t_2$.
Loosely speaking, if the solution $\mathbf{x}(t)$ satisfies $\mathbf{x}(t_1)\in\mathcal{R}_i\setminus\mathcal{Z}_i$ for $i\in\mathbb{I}\setminus\mathbb{I}_a$ and for some $t_1\geq 0$, where $-\mathbf{x}_d\in\mathcal{D}_{\epsilon}(\mathcal{U}_i)$, then for sufficiently large $\kappa$, the control input drives $\mathbf{x}$ out of $\mathcal{R}_i$ and there exists a finite time $t_2 > t_1$ such that the solution $\mathbf{x}(t)$ remains in $\mathcal{M}_0$ away from $\mathcal{R}_i$ for all $t\geq t_2$.

According to Claim \ref{lemma:claim1} of Lemma \ref{lemma:behaviour_of_x_in_Ri}, if the solution $\mathbf{x}(t)$ satisfies $\mathbf{x}(t_1)\in\mathcal{R}_i\setminus\mathcal{Z}_i$ for some $i\in\mathbb{I}_a$ and for some $t_1\geq 0$, then, after $t_1$, for sufficiently large $\kappa$, it enters into $\mathcal{M}_0\setminus\mathcal{R}_i$ only via the set $\mathcal{P}(\mathbf{x}_d, \mathcal{R}_i)$, which is defined using \eqref{projection_set_definition} and 
contains points $\mathbf{x}$ in $\mathcal{R}_i$ that satisfy $d_s(\mathbf{x}_d, \mathbf{x}) = d_s(\mathbf{x}_d, \mathcal{R}_i)$.
In other words, if $\mathbf{x}(t_1)\in\mathcal{R}_i\setminus\mathcal{Z}_i$ for some $i\in\mathbb{I}_a$, and if $\kappa > \bar{\kappa}$, then $\mathbf{x}(t_2)\in\mathcal{P}(\mathbf{x}_d, \mathcal{R}_i)$ for some $t_2\geq t_1$, $\mathbf{x}(t)\in\mathcal{R}_i$ for all $t\in[t_1, t_2]$, and there exists $\gamma > 0$ such that $\mathbf{x}(t)\in\mathcal{M}_0\setminus\mathcal{R}_i$ for all $t\in(t_2, t_2 + \gamma]$.
In this case, it follows that 
\[d_s(\mathbf{x}(t_2), \mathbf{x}_d) \leq d_s(\mathbf{x}(t_1), \mathbf{x}_d).\]
\begin{figure}
    \centering
    \includegraphics[width=0.75\linewidth]{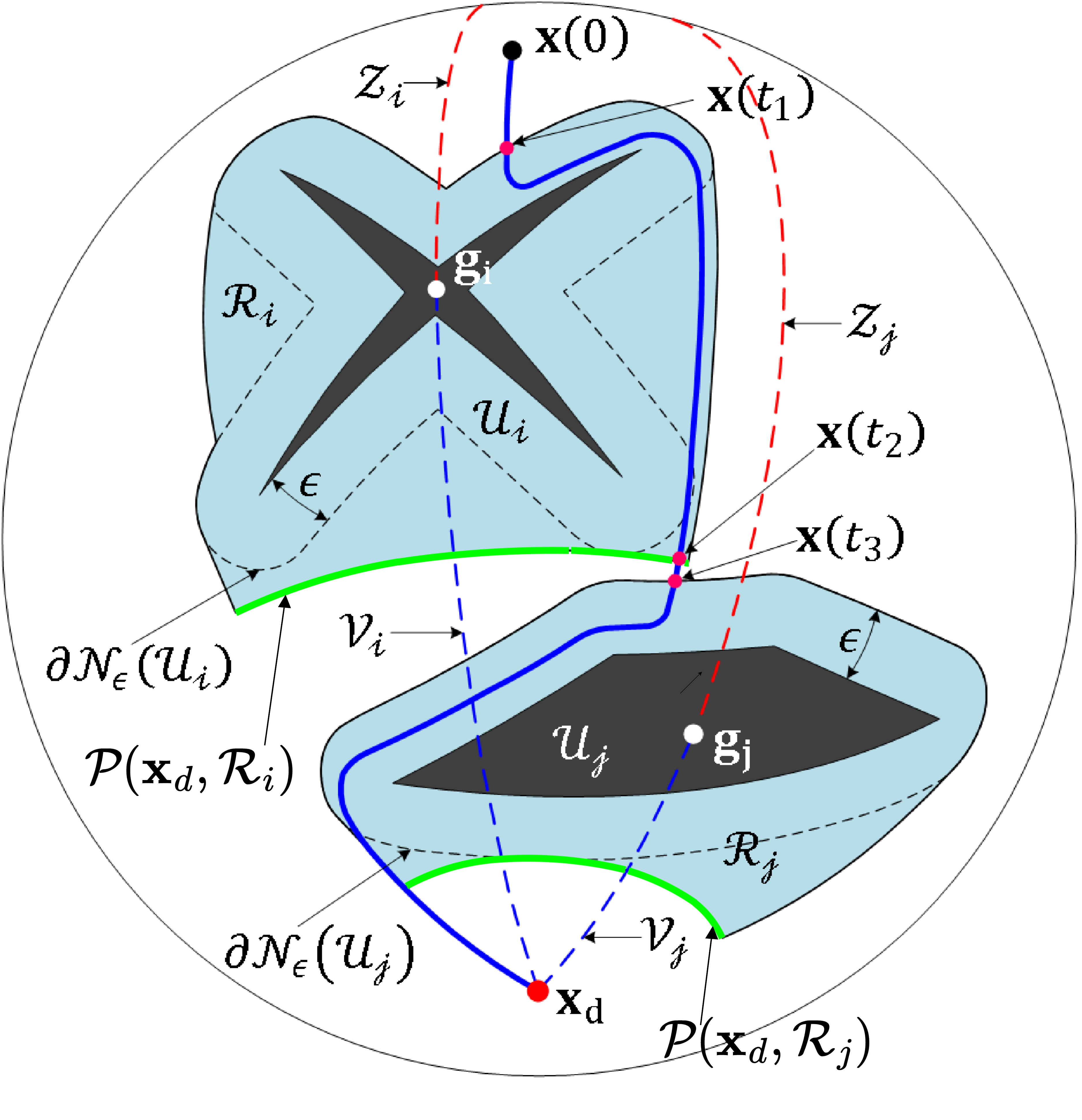}
    \caption{Illustration of an $\mathbf{x}$-trajectory, initialized at $\mathbf{x}(0)$ and converging to $\mathbf{x}_d$.}
    \label{fig:enter-label}
\end{figure}

Suppose there exists $t_3 > t_2$ such that $\mathbf{x}(t_3)\in\mathcal{R}_j$ for some $j\in\mathbb{I}\setminus\{i\}$ and $\mathbf{x}(t)\in\mathcal{M}_{0}\setminus\mathcal{R}$ for all $t\in(t_2, t_3)$.
Consequently, since $\mathbf{u}(\mathbf{x}) = k_1\mathbf{x}_d, k_1 > 0,$ for all $\mathbf{x}\in\mathcal{M}_{0}\setminus\mathcal{R}$, it follows that $\dot{d}_s(\mathbf{x}(t), \mathbf{x}_d) < 0$ for all $t\in [t_2, t_3]$, and therefore
\[d_s(\mathbf{x}(t_3), \mathbf{x}_d) < d_s(\mathbf{x}(t_2), \mathbf{x}_d).\]

Suppose $\mathbf{x}(t)\notin\bigcup_{i\in\mathbb{I}}\left(\mathcal{R}_i\cap\mathcal{Z}_i\right)$ for all $t\geq t_3$.
Then repeated application of Lemma~\ref{lemma:behaviour_of_x_in_Ri} allows us to conclude that, for all $t\geq t_3$, whenever $\mathbf{x}(t)$ belongs to $\mathcal{R}_i\setminus\mathcal{Z}_i$ for some $i\in\mathbb{I}_a$, under the control input \eqref{n-sphere-control-law}, it exits $\mathcal{R}_i$ in finite time and enters the set $\mathcal{M}_0\setminus\mathcal{R}_i$,
and the value of $d_s(\mathbf{x}(t), \mathbf{x}_d)$ at the exit time from $\mathcal{R}_i$ is less than its value at the entry time to $\mathcal{R}_i$.
Moreover, if $\mathbf{x}(t)$ subsequently enters the set $\mathcal{R}_j\setminus\mathcal{Z}_j$, where $j\in\mathbb{I}\setminus\{i\}$, then the value of $d_s(\mathbf{x}(t), \mathbf{x}_d)$ at the instance when $\mathbf{x}(t)$ entered the set $\mathcal{R}_j\setminus\mathcal{Z}_j$ is strictly lower than its value at the instance when it exited the set $\mathcal{R}_i$.
Since $d_s(\mathbf{x}, \mathbf{x}_d)\geq 0$ for all $\mathbf{x}\in\mathbb{S}^n$, a decrease in its value at the consecutive instances at which, under the control input \eqref{n-sphere-control-law}, the solution $\mathbf{x}(t)$ enters the set $\mathcal{R}$ and leaves the set $\mathcal{R}$, combined with the fact that $\dot{d}_s(\mathbf{x}, \mathbf{x}_d) < 0$ for all $\mathbf{x}\in\mathcal{M}_0\setminus\left(\mathcal{R}\cup\{\mathbf{x}_d, -\mathbf{x}_d\}\right)$, implies that $\mathbf{x}(t)$ satisfies \eqref{almost_attractivity_equation}.
Notice that, according to \eqref{set_definition_Vi_Zi}, for each $i\in\mathbb{I}$, the set $\mathcal{R}_i\cap\mathcal{Z}_i$ has zero Lebesgue measure on $\mathbb{S}^n$.
Consequently, if we show that the set of initial conditions in $\mathcal{M}_0$ from which the solutions $\mathbf{x}(t)$ to the closed-loop system \eqref{system_dynamics}-\eqref{n-sphere-control-law} satisfies $\mathbf{x}(t_1)\in\mathcal{R}_i\cap\mathcal{Z}_i$ for some $t_1\geq 0$ and some $i\in\mathbb{I}$ has zero Lebesgue measure on $\mathbb{S}^n$, then one can guarantee that for sufficiently large value of $\kappa$ any solution $\mathbf{x}(t)$ to the closed-loop system \eqref{system_dynamics}-\eqref{n-sphere-control-law}, initialized at any $\mathbf{x}(0)\in\mathcal{M}_0$ outside a set of Lebesgue measure zero, satisfies \eqref{almost_attractivity_equation}, and the proof of Theorem \ref{main_theorem} will be complete.

We proceed to show that the set of initial conditions in $\mathcal{M}_0$ from which the solutions $\mathbf{x}(t)$ to the closed-loop system \eqref{system_dynamics}-\eqref{n-sphere-control-law} satisfy $\mathbf{x}(t_1)\in\mathcal{R}_i\cap\mathcal{Z}_i$ for some $t_1\geq 0$ and some $i\in\mathbb{I}$ has zero Lebesgue measure.
If there exists $t_1\geq 0$ such that $\mathbf{x}(t_1)\in\mathcal{R}_i\setminus\mathcal{Z}_i$ for some $i\in\mathbb{I}$, then it follows from Lemma \ref{lemma:angle_keeps_changing} and Remark \ref{remark:possibilities} that there does not exist $t_2\geq0$ such that $\mathbf{x}(t_2)\in\mathcal{Z}_i$ and $\mathbf{x}(t)\in\mathcal{R}_i$ for all $t\in[t_1, t_2]$.
Consequently, if there exists $t_1\geq 0$ such that $\mathbf{x}(t_1)\in\mathcal{R}_i\cap\mathcal{Z}_i$ for some $i\in\mathbb{I}$, then either $\mathbf{x}(0)\in\mathcal{R}_i\cap\mathcal{Z}_i$ or there exists $s\in[0, t_1]$ such that $\mathbf{x}(s)\in\partial\mathcal{R}_i\cap\mathcal{Z}_i\cap\mathcal{M}_{\epsilon}$ and $\mathbf{x}(0)\in\mathcal{M}_{0}\setminus\mathcal{R}_i$. Since the set $\mathcal{Z}_i$, defined in \eqref{set_definition_Vi_Zi}, has zero Lebesgue measure for every $i\in\mathbb{I}$, it follows that if $\mathbf{x}(0)\in\mathcal{R}_i\cap\mathcal{Z}_i$, then the set of initial conditions in $\mathcal{M}_0$ from which the solutions $\mathbf{x}(t)$ to the closed-loop system \eqref{system_dynamics}-\eqref{n-sphere-control-law} satisfies $\mathbf{x}(t_1)\in\mathcal{R}_i\cap\mathcal{Z}_i$ for some $t_1\geq 0$ and some $i\in\mathbb{I}$ has zero Lebesgue measure.
Therefore, we proceed to analyze the latter case where there exists $s\in[0, t_1]$ such that $\mathbf{x}(s)\in\partial\mathcal{R}_i\cap\mathcal{Z}_i\cap\mathcal{M}_{\epsilon}$ and $\mathbf{x}(0)\in\mathcal{M}_{0}\setminus\mathcal{R}_i$.

According to \eqref{set_definition_Vi_Zi} and \eqref{definition:Ri}, for every $i\in\mathbb{I}$, the intersection set $\partial\mathcal{R}_i\cap\mathcal{Z}_i\cap\mathcal{M}_{\epsilon}$ is a singleton, and the unique element is given by
\begin{equation}\label{definition_of_point_Si}
\partial\mathcal{R}_i\cap\mathcal{Z}_i\cap\mathcal{M}_{\epsilon} = \{\mathbf{s}_i\},
\end{equation}
where if $i\in\mathbb{I}\setminus\mathbb{I}_a$, then $\mathcal{R}_i = \mathcal{S}_i(-\mathbf{x}_d)\setminus\mathcal{U}_i^{\circ}$, as stated in Section \ref{section:stability_analysis} and $\mathcal{S}_i(-\mathbf{x}_d)$ is obtained by replacing $\mathbf{x}_d$ in \eqref{definition:Si} with $-\mathbf{x}_d$.
We show that the set of initial conditions in $\mathcal{M}_0\setminus\mathcal{R}$ from which the solution $\mathbf{x}(t)$ to the closed-loop system \eqref{system_dynamics}-\eqref{n-sphere-control-law}
satisfies $\mathbf{x}(s)= \mathbf{s}_i$ for some $s\geq0$ and for some $i\in\mathbb{I}$ has zero Lebesgue measure, 
where the set $\mathcal{R}$ is defined in \eqref{collection_of_all_Ri}.

According to Lemma \ref{lemma:n-sphere-forward-invariance}, the set $\mathcal{M}_0$ which is a closed subset of $\mathbb{S}^n$, is forward invariant for the closed-loop system \eqref{system_dynamics}-\eqref{n-sphere-control-law}.
Consequently, since $\mathbf{P}(\mathbf{x})\mathbf{u}(\mathbf{x})$ is locally Lipschitz in $\mathbf{x}$ over $\mathcal{M}_0$, as established earlier in Lemma \ref{lemma:lipschitz_continuity}, it follows from \cite[Theorem 3.3]{khalil2002nonlinear} that the solution $\mathbf{x}(t)$ to the closed-loop system \eqref{system_dynamics}-\eqref{n-sphere-control-law}, initialized at any $\mathbf{x}(0)\in\mathcal{M}_0$, is unique and defined for all $t\geq 0$.

Let $\phi(t, \mathbf{x} (0))$ denote the solution $\mathbf{x}(t)$ to the closed-loop system \eqref{system_dynamics}-\eqref{n-sphere-control-law}, starting from the initial condition $\mathbf{x}(0)$.
Since the solution $\phi(t, \mathbf{x}(0))$ is unique for every initial condition $\mathbf{x}(0)\in\mathcal{M}_0$ and defined for all $t\geq 0$, it follows that for any initial conditions $\mathbf{x}_1, \mathbf{x}_2\in\mathcal{M}_0$, if there exist  $t_1\geq 0$ and $t_2\geq 0$ such that $\phi(t_1, \mathbf{x}_1) = \phi(t_2, \mathbf{x}_2)$, then either $\mathbf{x}_2=\phi(\bar{t}, \mathbf{x}_1)$ for some $\bar{t}\in[0, t_1]$  or $\mathbf{x}_1=\phi(\underline{t}, \mathbf{x}_2)$ for some $\underline{t}\in[0, t_2]$.
In other words, if solutions to the closed-loop system \eqref{system_dynamics}-\eqref{n-sphere-control-law} originating from any two distinct initial conditions, $\mathbf{x}_0$ and $\mathbf{x}_1$, in $\mathcal{M}_0$ reach a common point $\mathbf{x}$ in $\mathcal{M}_0$ in a finite time, then one of these solution trajectories must be a subset of the other.
As a result, since the set of points in $\mathcal{M}_0$ that belong to any given solution $\mathbf{x}(t)$ to the closed-loop system \eqref{system_dynamics}-\eqref{n-sphere-control-law} has zero Lebesgue measure, it follows that the set of initial conditions in $\mathcal{M}_0\setminus\mathcal{R}_i$ from which the solution $\mathbf{x}(t)$ satisfies $\mathbf{x}(s) = \mathbf{s}_i$ for some time $s\geq 0$ has zero Lebesgue measure, where for every $i\in\mathbb{I}$, the point $\mathbf{s}_i$ is defined in \eqref{definition_of_point_Si}.
Therefore, the set of initial conditions in $\mathcal{M}_0$ from which the solutions $\mathbf{x}(t)$ to the closed-loop system \eqref{system_dynamics}-\eqref{n-sphere-control-law} satisfy $\mathbf{x}(t_1)\in\mathcal{R}_i\cap\mathcal{Z}_i$ for some $t_1\geq 0$ and some $i\in\mathbb{I}$ has zero Lebesgue measure.

\subsection{Proof of Lemma \ref{lemma:Q_set_is_geodesic}}\label{proof:lemma:Q_set_is_geodesic}
Since $\mathbf{0}\notin\mathcal{L}_s(\mathbf{a}, \mathbf{b})$, it follows that if the vectors $\mathbf{a}$ and $\mathbf{b}$ are collinear, then $\psi(\mathbf{a}) = \psi(\mathbf{b})$ and the results follow directly, where the function $\psi(\cdot)$ is defined in \eqref{psi_function_definition}.
Therefore, we consider the case where the vectors $\mathbf{a}$ and $\mathbf{b}$ are not collinear, \textit{i.e.}, there does not exist $q\in\mathbb{R}$ such that $\mathbf{a} = q\mathbf{b}$.
Consequently, $\theta_{\mathbf{a}, \mathbf{b}}\in(0, \pi)$, where $\theta_{\mathbf{a}, \mathbf{b}} = \arccos\left(
\frac{\mathbf{a}^{\top}\mathbf{b}}
{\|\mathbf{a}\|\|\mathbf{b}\|}
\right)$.

According to \eqref{expression_for_line_segment}, for each $\mathbf{p}_{\mathcal{L}}\in\mathcal{L}_{s}(\mathbf{a}, \mathbf{b})$, there exists a unique $\lambda_{\mathbf{p}}\in[0, 1]$ such that
\begin{equation}\label{point_p_on_line_segment}
    \mathbf{p}_{\mathcal{L}} = (1-\lambda_{\mathbf{p}})\mathbf{a} + \lambda_{\mathbf{p}}\mathbf{b}.
\end{equation}
Using \eqref{psi_function_definition}, $\psi(\mathbf{p}_{\mathcal{L}})$ is evaluated as
\begin{equation}
    \label{expression_for_p_in_a_and_b}\psi(\mathbf{p}_{\mathcal{L}}) = \frac{1-\lambda_{\mathbf{p}}}{\alpha(\lambda_{\mathbf{p}})}\mathbf{a} + \frac{\lambda_{\mathbf{p}}}{\alpha(\lambda_{\mathbf{p}})}\mathbf{b},
\end{equation}
where $\alpha(\lambda_{\mathbf{p}})$ is given by
\[\alpha(\lambda_{\mathbf{p}}) = \sqrt{(1-\lambda_{\mathbf{p}})^2\|\mathbf{a}\|^2 + \lambda_{\mathbf{p}}^2\|\mathbf{b}\|^2 + 2\lambda_{\mathbf{p}}(1-\lambda_{\mathbf{p}})\mathbf{a}^\top\mathbf{b}}.\]
Since $\mathbf{0}\notin\mathcal{L}_s(\mathbf{a}, \mathbf{b})$, $\psi(\mathbf{p}_{\mathcal{L}})$ is well-defined for all $\mathbf{p}_{\mathcal{L}}\in\mathcal{L}_{s}(\mathbf{a}, \mathbf{b})$.
Furthermore, since $\mathbf{p}_{\mathcal{L}}\in\mathcal{L}_s(\mathbf{a}, \mathbf{b})$, it follows from \eqref{Q_set_definition} that $\psi(\mathbf{p}_{\mathcal{L}})\in\mathcal{Q}(\mathbf{a}, \mathbf{b})$.

Now consider the geodesic $\mathcal{G}(\psi(\mathbf{a}), \psi(\mathbf{b}))$, defined as in \eqref{geodesic_expression}.
For every $\mathbf{q}_{\mathcal{G}}\in\mathcal{G}(\psi(\mathbf{a}), \psi(\mathbf{b}))$, there exists a unique $\lambda_{\mathbf{q}}\in[0, 1]$ such that
\begin{equation}\label{expression_q_in_e1_and_e2}
    \mathbf{q}_{\mathcal{G}} = \frac{\sin((1-\lambda_{\mathbf{q}})\theta_{\mathbf{a} ,\mathbf{b}})}{\|\mathbf{a}\|\sin\theta_{\mathbf{a} ,\mathbf{b}}} \mathbf{a} + \frac{\sin(\lambda_{\mathbf{q}}\theta_{\mathbf{a} ,\mathbf{b}})}{\|\mathbf{b}\|\sin \theta_{\mathbf{a} ,\mathbf{b}}}\mathbf{b},
\end{equation}
where $\theta_{\mathbf{a}, \mathbf{b}}\in(0, \pi)$.
Comparing the right-hand sides of \eqref{expression_for_p_in_a_and_b} and \eqref{expression_q_in_e1_and_e2}, equating coefficients of $\mathbf{a}$ and $\mathbf{b}$, one obtains
\begin{equation*}
    \frac{1-\lambda_{\mathbf{p}}}{\alpha(\lambda_{\mathbf{p}})} = \frac{\sin((1-\lambda_{\mathbf{q}})\theta_{\mathbf{a} ,\mathbf{b}})}{\|\mathbf{a}\|\sin\theta_{\mathbf{a} ,\mathbf{b}}}\text{ and } \frac{\lambda_{\mathbf{p}}}{\alpha(\lambda_{\mathbf{p}})} = \frac{\sin(\lambda_{\mathbf{q}}\theta_{\mathbf{a} ,\mathbf{b}})}{\|\mathbf{b}\|\sin \theta_{\mathbf{a} ,\mathbf{b}}}.
\end{equation*}
Therefore, 
\begin{equation}\label{rho_function_definition}
    \lambda_{\mathbf{p}} = \frac{\|\mathbf{a}\|\sin(\lambda_{\mathbf{q}}\theta_{\mathbf{a}, \mathbf{b}})}{\|\mathbf{a}\|\sin(\lambda_{\mathbf{q}}\theta_{\mathbf{a}, \mathbf{b}})  + \|\mathbf{b}\|\sin((1-\lambda_{\mathbf{q}})\theta_{\mathbf{a}, \mathbf{b}})}=:\rho(\lambda_{\mathbf{q}}).
\end{equation}

To show that $\mathcal{Q}(\mathbf{a}, \mathbf{b})=\mathcal{G}(\psi(\mathbf{a}), \psi(\mathbf{b}))$, it is sufficient to show that $\rho(0) = 0$, $\rho(1) = 1$ and $\rho(\lambda_{\mathbf{q}})$ is strictly increasing over $[0, 1]$ for every $\theta_{\mathbf{a}, \mathbf{b}}\in(0, \pi)$.
This will ensure that the mapping $\rho:[0, 1]\to[0, 1]$ is bijective. 
In other words, for every $\lambda_{\mathbf{p}}\in[0, 1]$, there will exist a unique $\lambda_{\mathbf{q}}\in[0, 1]$ such that $\psi\left(\mathbf{p}_{\mathcal{L}}\right) = \mathbf{q}_{\mathcal{G}}$, and for every $\lambda_{\mathbf{q}}\in[0, 1]$, there exists a unique $\lambda_{\mathbf{p}}\in[0, 1]$ such that $\psi\left(\mathbf{p}_{\mathcal{L}}\right) = \mathbf{q}_{\mathcal{G}}$, where $\mathbf{p}_{\mathcal{L}}$ and $\mathbf{q}_{\mathcal{G}}$ are defined in \eqref{point_p_on_line_segment} and \eqref{expression_q_in_e1_and_e2}, respectively.

Using \eqref{rho_function_definition}, it is straightforward to verify that $\rho(0) = 0$ and $\rho(1) = 1$ for every $\theta_{\mathbf{a}, \mathbf{b}}\in(0, \pi)$.
To show that $\rho(\lambda_{\mathbf{q}})$ is strictly increasing over $[0, 1]$ for every $\theta_{\mathbf{a}, \mathbf{b}}\in(0, \pi)$, it is sufficient to show that $\frac{d}{d\lambda_{\mathbf{q}}}\rho(\lambda_{\mathbf{q}}) > 0$ for all $\lambda_{\mathbf{q}}\in[0, 1]$ and for every $\theta_{\mathbf{a}, \mathbf{b}}\in(0, \pi)$.

Differentiating $\rho(\lambda_{\mathbf{q}})$ with respect to $\lambda_{\mathbf{q}}$, one obtains
\begin{equation}
\frac{d}{d{\lambda_{\mathbf{q}}}}\rho(\lambda_{\mathbf{q}}) = \frac{\|\mathbf{a}\|\|\mathbf{b}\|\theta_{\mathbf{a}, \mathbf{b}}\sin(\theta_{\mathbf{a},\mathbf{b}})}{\left(\|\mathbf{a}\|\sin(\lambda_{\mathbf{q}}\theta_{\mathbf{a}, \mathbf{b}}) + \|\mathbf{b}\|\sin((1-\lambda_{\mathbf{q}})\theta_{\mathbf{a}, \mathbf{b}})\right)^2}.
\end{equation}
Since $\theta_{\mathbf{a}, \mathbf{b}}\in(0, \pi)$, it follows that $\frac{d}{d\lambda_{\mathbf{q}}}\rho(\lambda_{\mathbf{q}}) > 0$ for all $\lambda_{\mathbf{q}}\in[0, 1]$.
Therefore, $\rho(\lambda_{\mathbf{q}})$ is strictly increasing over $[0, 1]$ for every $\theta_{\mathbf{a}, \mathbf{b}}\in(0, \pi)$.
This completes the proof of Lemma \ref{lemma:Q_set_is_geodesic}.

\subsection{Proof of Lemma \ref{lemma:behaviour_of_x_in_Ri}}
\label{proof:lemma:behavious_of_x_in_Ri}
\subsubsection{Proof of Claim \ref{lemma:claim1}}

Since $\left(\mathcal{R}_i\setminus\mathcal{Z}_i\right)\subset\mathcal{F}_i$, it follows from Lemma \ref{lemma:angle_keeps_changing} and as described in Remark \ref{remark:possibilities} that 
if $\mathbf{x}(t_1)\in\mathcal{R}_i\setminus\mathcal{Z}_i$ for some $t_1 \geq 0$, then there are three possible cases as mentioned below:

\noindent{\bf Case 1:}\label{case1}
Suppose there exists $s_1 > t_1$ such that $\mathbf{x}(s_1)\in\mathcal{M}_{0}\setminus\mathcal{R}_i$.
Consequently, since $\mathbf{x}(t)$ is continuous for all $t\in[t_1, s_1]$, there exist $t_2\in[t_1, s_1)$ and $\gamma > 0$ such that $\mathbf{x}(t_2)\in\partial\mathcal{R}_i$, $\mathbf{x}(t)\in\mathcal{R}_i$ for all $t\in[t_1, t_2]$ and $\mathbf{x}(t)\in\mathcal{M}_{0}\setminus\mathcal{R}_i$ for all $t\in(t_2, t_2 + \gamma]$.
To show that $\mathbf{x}(t_2)\in\mathcal{P}(\mathbf{x}_d, \mathcal{R}_i)$, we require the following fact:
\begin{fact}\label{fact:only_one_exit}
    For the closed-loop system \eqref{system_dynamics}-\eqref{n-sphere-control-law}, under Assumptions \ref{assumption:non-overlapping-constraints} and \ref{assumption:sufficient_separation}, for $\mathbf{x}\in\mathcal{R}_i$, $\mathbf{P}(\mathbf{x})\mathbf{u}(\mathbf{x})\notin\mathcal{T}_{\mathbf{x}}(\mathcal{R}_i)$ if and only if $\mathbf{x}\in\mathcal{P}(\mathbf{x}_d, \mathcal{R}_i)$, where $i\in\mathbb{I}_a$ and $\mathcal{T}_{\mathbf{x}}(\mathcal{R}_i)$ denotes the tangent cone to $\mathcal{R}_i$ on $\mathbb{S}^n$ at $\mathbf{x}$.
\end{fact}
\proof{Since $\mathcal{R}_i^{\circ}\subset\mathbb{S}^n$, one has $\mathcal{T}_{\mathbf{x}}(\mathcal{R}_i) = \mathbf{T}_{\mathbf{x}}(\mathbb{S}^n)$ for all $\mathbf{x}\in\mathcal{R}_i^{\circ}$, and it follows that $\mathbf{P}(\mathbf{x})\mathbf{u}(\mathbf{x})\in\mathcal{T}_{\mathbf{x}}(\mathcal{R}_i)$ for every $\mathbf{x}\in\mathcal{R}_i^{\circ}$.
Since $\epsilon < \bar{\epsilon}$, as stated in Section \ref{section:feedback_control_design}, one has $\mathbf{x}_d\notin\mathcal{N}_{\epsilon}(\mathcal{U}_i)$ for all $i\in\mathbb{I}_a$.
Therefore, by \eqref{definition:Ri}, one has $\mathbf{x}_d\notin\mathcal{R}_i$ for each $i\in\mathbb{I}_a$.
It follows that for each $i\in\mathbb{I}_a$, $\mathcal{P}(\mathbf{x}_d, \mathcal{R}_i)\subset\partial\mathcal{R}_i$.
Consequently, using \eqref{definition:Ri}, the boundary of $\mathcal{R}_i$ on $\mathbb{S}^n$ can be partitioned as follows:
\begin{equation*}
    \partial\mathcal{R}_i = \partial\mathcal{U}_i   \cup \mathcal{Y}_i \cup \mathcal{P}(\mathbf{x}_d, \mathcal{R}_i),
\end{equation*}
where $\mathcal{Y}_i = \left(\partial\mathcal{R}_i\cap\mathcal{M}_{\epsilon}\right)\setminus\mathcal{P}(\mathbf{x}_d, \mathcal{R}_i)$ and $i\in\mathbb{I}_a$.

Since $\mathcal{M}_0$ is forward invariant for the closed-loop system \eqref{system_dynamics}-\eqref{n-sphere-control-law}, as established in Lemma \ref{lemma:n-sphere-forward-invariance}, it follows that $\mathbf{P}(\mathbf{x})\mathbf{u}(\mathbf{x})\in\mathcal{T}_{\mathbf{x}}(\mathcal{R}_i)$ for all $\mathbf{x}\in\partial\mathcal{U}_i$.
We proceed to show that $\mathbf{P}(\mathbf{x})\mathbf{u}(\mathbf{x})\in\mathcal{T}_{\mathbf{x}}(\mathcal{R}_i)$ for all $\mathbf{x}\in\mathcal{Y}_i$.
By construction in \eqref{definition:Si} and \eqref{definition:Ri}, for every $\mathbf{x}\in\mathcal{Y}_i$, there exists $\mathbf{p}(\mathbf{x})\in\mathcal{R}_i\setminus\{\mathbf{x}\}$ such that the geodesic $\mathcal{G}(\mathbf{x}, \mathbf{p}(\mathbf{x}))$ connecting $\mathbf{x}$ to $\mathbf{p}(\mathbf{x})$ is a subset of the geodesic $\mathcal{G}(\mathbf{x}, \mathbf{x}_d)$ and belongs to $\mathcal{R}_i$.
In other words, for each $\mathbf{x}\in\mathcal{Y}_i$, there exists $\mathbf{p}(\mathbf{x})\in\mathcal{R}_i\setminus\{\mathbf{x}\}$ such that $\mathcal{G}(\mathbf{x}, \mathbf{p}(\mathbf{x}))\subset \mathcal{G}(\mathbf{x}, \mathbf{x}_d)\cap\mathcal{R}_i$.
Additionally, for every $\mathbf{x}\in\mathcal{M}_0\setminus\{\mathbf{x}_d, -\mathbf{x}_d\}$, the vector field $\mathbf{P}(\mathbf{x})\mathbf{x}_d$ steers $\mathbf{x}$ along the geodesic $\mathcal{G}(\mathbf{x}, \mathbf{x}_d)$ towards $\mathbf{x}_d$.
Consequently, for every $\mathbf{x}\in\mathcal{Y}_i$, the vector field $\mathbf{P}(\mathbf{x})\mathbf{x}_d$ cannot drive $\mathbf{x}$ out of $\mathcal{R}_i$, and therefore, $\mathbf{P}(\mathbf{x})\mathbf{x}_d\in\mathcal{T}_{\mathbf{x}}(\mathcal{R}_i)$ for all $\mathbf{x}\in\mathcal{Y}_i$.
In addition, since $\mathcal{Y}_i\subset\mathcal{M}_{\epsilon}$, by \eqref{n-sphere-control-law}, one has $\mathbf{u}(\mathbf{x}) = k_1\mathbf{x}_d$ for all $\mathbf{x}\in\mathcal{Y}_i$.
As a result, $\mathbf{P}(\mathbf{x})\mathbf{u}(\mathbf{x})\in\mathcal{T}_{\mathbf{x}}(\mathcal{R}_i)$ for all $\mathbf{x}\in\mathcal{Y}_i$.
It remains to show that for every $\mathbf{x}\in\mathcal{P}(\mathbf{x}_d, \mathcal{R}_i)$, $\mathbf{P}(\mathbf{x})\mathbf{u}(\mathbf{x})\notin\mathcal{T}_{\mathbf{x}}(\mathcal{R}_i)$.

The set $\mathcal{P}(\mathbf{x}_d, \mathcal{R}_i)$, which is defined using \eqref{projection_set_definition}, contains points $\mathbf{x}$ in $\mathcal{R}_i$ that satisfy $d_s(\mathbf{x}_d, \mathbf{x}) = d_s(\mathbf{x}_d, \mathcal{R}_i)$.
Consequently, since $\mathcal{P}(\mathbf{x}_d, \mathcal{R}_i)\subset\partial\mathcal{R}_i$, it holds that $\mathcal{G}(\mathbf{x}, \mathbf{x}_d)\cap\mathcal{R}_i^{\circ} = \emptyset$ for all $\mathbf{x}\in\mathcal{P}(\mathbf{x}_d, \mathcal{R}_i)$.
Therefore, for each $\mathbf{x}\in\mathcal{P}(\mathbf{x}_d, \mathcal{R}_i)$, $\mathbf{P}(\mathbf{x})\mathbf{x}_d\notin\mathcal{T}_{\mathbf{x}}(\mathcal{R}_i)$.
Furthermore, according to \eqref{definition:Ri}, one has $\mathcal{P}(\mathbf{x}_d, \mathcal{R}_i)\subset\mathcal{M}_{\epsilon}$.
Therefore, it follows from \eqref{n-sphere-control-law} that $\mathbf{u}(\mathbf{x}) = k_1\mathbf{x}_d$ for every $\mathbf{x}\in\mathcal{P}(\mathbf{x}_d, \mathcal{R}_i)$.
Consequently, $\mathbf{P}(\mathbf{x})\mathbf{u}(\mathbf{x})\notin\mathcal{T}_{\mathbf{x}}(\mathcal{R}_i)$ for all $\mathbf{x}\in\mathcal{P}(\mathbf{x}_d, \mathcal{R}_i)$, and the proof is complete.}

Since $\mathbf{x}(t_2)\in\partial\mathcal{R}_i$ and there exists $\gamma > 0$ such that $\mathbf{x}(t)\notin\mathcal{R}_i$ for all $t\in(t_2, t_2 + \gamma]$,
it follows from Fact \ref{fact:only_one_exit} that $\mathbf{x}(t_2)\in\mathcal{P}(\mathbf{x}_d, \mathcal{R}_i)$.
Therefore, Claim \ref{lemma:claim1} of Lemma \ref{lemma:behaviour_of_x_in_Ri} holds in Case 1.

\noindent{\bf Case 2:}\label{case2}
Suppose $\mathbf{x}(t)\in\mathcal{R}_i\setminus(\mathcal{V}_i\cup\mathcal{Z}_i)$ for all $t\geq t_1$ and $\Lim_{t\to\infty}d_s(\mathbf{x}(t), \mathcal{V}_i) = 0$, where $i\in\mathbb{I}_a$ and the set $\mathcal{V}_i$ is defined in \eqref{set_definition_Vi_Zi}.
It follows that for every $\nu > 0$ there exists $s_{\nu} \geq t_1$ such that $\mathbf{x}(t)\in\mathcal{D}_{\nu}(\mathcal{V}_i)\cap\mathcal{R}_i$ for all $t\geq s_{\nu}$, where $\mathcal{D}_{\nu}(\mathcal{V}_i)$ is defined as in \eqref{dilation_on_sphere}.
We show that for each $i\in\mathbb{I}_a$, there exist $\kappa_i > 0$, $\nu > 0$ and $\zeta_{\nu} > 0$ such that, for any $\kappa > \kappa_i$, $\dot{d}_s(\mathbf{x}, \mathbf{x}_d) \leq -\zeta_{\nu} < 0$ for all $\mathbf{x}\in\mathcal{D}_{\nu}(\mathcal{V}_i)\cap\mathcal{R}_i$.
Since $\mathbf{x}_d\notin\mathcal{R}_i$  and $\mathcal{D}_{\nu}(\mathcal{V}_i)\cap\mathcal{R}_i$ is a closed subset of $\mathbb{S}^n$ for any $i\in\mathbb{I}_a$, it follows that there exists $t_2 > t_1$ such that $\mathbf{x}(t_2)\notin\mathcal{D}_{\nu}(\mathcal{V}_i)\cap\mathcal{R}_i$, and therefore, Case 2 does not hold.
To proceed with the proof, we require the following fact:

\begin{fact}\label{fact:for_vi}For the closed-loop system \eqref{system_dynamics}-\eqref{n-sphere-control-law} under Assumptions \ref{assumption:non-overlapping-constraints} and \ref{assumption:sufficient_separation}, for each $i\in\mathbb{I}$, there exist $\kappa_i > 0$ and $\zeta_i > 0$ such that if $\kappa > \kappa_i$, then $\dot{d_s}(\mathbf{x}, \mathbf{x}_d) \leq -\zeta_i < 0$ for all $\mathbf{x}\in\mathcal{V}_i\cap\mathcal{R}_i$.
\end{fact}
\proof{See Appendix \ref{proof:fact:for_vi}.}

According to Fact \ref{fact:for_vi}, for each $i\in\mathbb{I}$, there exist $\kappa_i > 0$ and $\zeta_i > 0$ such that if $\kappa > \kappa_i$, then $\dot{d}_s(\mathbf{x}, \mathbf{x}_d) \leq -\zeta_i$ for all $\mathbf{x}\in\mathcal{V}_i\cap\mathcal{R}_i$, where $\mathcal{V}_i\cap\mathcal{R}_i$ is a closed set on $\mathbb{S}^n$.
Additionally, $d_s(\mathbf{x}, \mathbf{x}_d)$ is smooth with respect to $\mathbf{x}$ over $\mathcal{M}_0$ and $\mathbf{P}(\mathbf{x})\mathbf{u}(\mathbf{x})$ is continuous over $\mathcal{M}_0$.
Therefore, $\dot{d}_s(\mathbf{x}, \mathbf{x}_d)$ is continuous in $\mathbf{x}$ over $\mathcal{M}_0$.
Consequently, there exist $\nu > 0$ and $\zeta_{\nu}\in(0, \zeta_i]$ such that if $\kappa > \kappa_i$, then $\dot{d}_s(\mathbf{x}, \mathbf{x}_d) \leq -\zeta_{\nu}$ for all $\mathbf{x}\in\mathcal{D}_{\nu}(\mathcal{V}_i)\cap\mathcal{R}_i$.
Furthermore, as mentioned earlier, there exists $s_{\nu} \geq t_1$ such that $\mathbf{x}(t)\in\mathcal{D}_{\nu}(\mathcal{V}_i)\cap\mathcal{R}_i$ for all $t\geq s_{\nu}$.
Therefore, $\dot{d}_s(\mathbf{x}(t), \mathbf{x}_d)\leq -\zeta_{\nu}$ for all $t\geq s_{\nu}$ as long as $\mathbf{x}(t)\in\mathcal{D}_{\nu}(\mathcal{V}_i)\cap\mathcal{R}_i$.
However, since $\epsilon < \bar{\epsilon}$, as stated in Section \ref{section:feedback_control_design}, one has $\mathbf{x}_d\notin\mathcal{D}_{\epsilon}(\mathcal{U}_i)$ for any $i\in\mathbb{I}$.
Therefore, it follows from \eqref{definition:Ri} and \eqref{definition:Ri_forspecialindex}, $\mathbf{x}_d\in\mathcal{M}_{0}\setminus\mathcal{R}$, where $\mathcal{R}$ is defined in \eqref{collection_of_all_Ri}.
Moreover, since the set $\mathcal{D}_{\nu}(\mathcal{V}_i)\cap\mathcal{R}_i$ is a closed subset of $\mathbb{S}^n$, it is compact.
Furthermore, since $\dot{d}_s(\mathbf{x}(t), \mathbf{x}_d)\leq -\zeta_{\nu} < 0$ for all $t\geq s_{\nu}$ as long as $\mathbf{x}(t)\in\mathcal{D}_{\nu}(\mathcal{V}_i)\cap\mathcal{R}_i$ and $\mathbf{x}_d\notin\mathcal{R}_i$, there exists $t_2> s_{\nu}$ such that $\mathbf{x}(t_2)\notin\mathcal{D}_{\nu}(\mathcal{V}_i)\cap\mathcal{R}_i$.
Consequently, for any $\kappa > \kappa_i$, where the existence of $\kappa_i > 0$ is established in Fact \ref{fact:for_vi}, Case 2 does not hold.

\noindent{\bf Case 3:} Suppose $\mathbf{x}(t)\in\mathcal{R}_i\cap\mathcal{V}_i$ for all $t\geq t_1$ and for $i\in\mathbb{I}_a$.
By Fact \ref{fact:for_vi}, there exist $\kappa_i > 0$ and $\zeta_i > 0$ such that, for any $\kappa > \kappa_i$, one has $\dot{d}_s(\mathbf{x}, \mathbf{x}_d) \leq -\zeta_i$ for all $\mathbf{x}\in\mathcal{R}_i\cap\mathcal{V}_i$.
Additionally, since $\epsilon < \bar{\epsilon}$, one has $\mathbf{x}_d\notin\mathcal{D}_{\epsilon}(\mathcal{U}_i)$ for any $i\in\mathbb{I}_a$.
Therefore, by \eqref{definition:Ri} and \eqref{definition:Ri_forspecialindex}, $\mathbf{x}_d\notin\mathcal{R}_i$. 
Moreover, since $\mathcal{R}_i\cap\mathcal{V}_i$ is a closed subset of $\mathbb{S}^n$ that does not contain $\mathbf{x}_d$, the value of $d_s(\mathbf{x}, \mathbf{x}_d)$ is lower bounded by a positive value on this set.
Consequently, there exists $t_2> t_1$ such that $\mathbf{x}(t_2)\notin\mathcal{R}_i\cap\mathcal{V}_i$.
As a result, for any $\kappa > \kappa_i$, Case 3 does not hold.
This completes the proof of Claim \ref{lemma:claim1} of Lemma \ref{lemma:behaviour_of_x_in_Ri}.

\subsubsection{Proof of Claim \ref{lemma:claim2}}
Since $\left(\mathcal{R}_i\setminus\mathcal{Z}_i\right)\subset\mathcal{F}_i$, it follows from Lemma \ref{lemma:angle_keeps_changing} and as described in Remark \ref{remark:possibilities} that 
if $\mathbf{x}(t_1)\in\mathcal{R}_i\setminus\mathcal{Z}_i$ for some $t_1 \geq 0$, then there are three possible cases as mentioned below:
\begin{enumerate}
    \item There exists $t_2> t_1$ such that $\mathbf{x}(t_2)\in\mathcal{M}_{\epsilon}\setminus\mathcal{R}_i$,
    \item $\mathbf{x}(t)\in\mathcal{R}_i\setminus(\mathcal{V}_i\cup\mathcal{Z}_i)$ for all $t\geq t_1$ and $\Lim_{t\to\infty}d_s(\mathbf{x}(t), \mathcal{V}_i) = 0$,
    \item $\mathbf{x}(t)\in\mathcal{R}_i\cap\mathcal{V}_i$ for all $t\geq t_1$,
\end{enumerate}
where for $i\in\mathbb{I}\setminus\mathbb{I}_a$, the set $\mathcal{V}_i$ is defined in \eqref{set_definition_Vi_Zi}.
Case 1 directly implies the existence of $t_2 > t_1$ such that $\mathbf{x}(t_2)\in\mathcal{M}_{\epsilon}\setminus\mathcal{R}_i$.
Using arguments similar to those used in the proof of Claim \ref{lemma:claim1} of Lemma \ref{lemma:behaviour_of_x_in_Ri}, one can establish the existence of $\kappa_i > 0$ such that for any $\kappa > \kappa_i$, cases 2 and 3 do not hold.

Next, we show that for $i\in\mathbb{I}\setminus\mathbb{I}_a$, if $\mathbf{x}(t_2)\in\mathcal{M}_{\epsilon}\setminus\mathcal{R}_i$ for some $t_2> t_1$, then $\mathbf{x}(t)\notin\mathcal{R}_i$ for all $t\geq t_2$.
According to \eqref{definition:Ri_forspecialindex}, for $i\in\mathbb{I}\setminus\mathbb{I}_a$, $\mathcal{R}_i = \mathcal{S}_i(-\mathbf{x}_d)\setminus\mathcal{U}_i^{\circ}$, where $\mathcal{S}_i(-\mathbf{x}_d)$ is obtained from \eqref{definition:Si} by replacing $\mathbf{x}_d$ with $-\mathbf{x}_d$.
It follows from \eqref{definition:Si} that $\mathcal{G}(\mathbf{x}, -\mathbf{x}_d)\subset\mathcal{S}(-\mathbf{x}_d)$ for all $\mathbf{x}\in\mathcal{R}_i$.
Consequently, for every $\mathbf{x}\in\mathcal{M}_{\epsilon}\setminus\mathcal{R}_i$, $\mathcal{G}(\mathbf{x}, \mathbf{x}_d)\cap\mathcal{R}_i = \emptyset$.
Moreover, since $i\in\mathbb{I}\setminus\mathbb{I}_a$, by \eqref{definition:Ia}, one has $-\mathbf{x}_d\in\mathcal{D}_{\epsilon}(\mathcal{U}_i)$.
Therefore, by \eqref{n-sphere-control-law}, for any $\mathbf{x}\in\mathcal{M}_{\epsilon}\setminus\mathcal{R}_i$, the control input becomes $\mathbf{u}(\mathbf{x}) = k_1\mathbf{x}_d$, and it steers $\mathbf{x}$ along the geodesic $\mathcal{G}(\mathbf{x}, \mathbf{x}_d)$ towards $\mathbf{x}_d$.
Consequently, since  $\mathcal{G}(\mathbf{x}, \mathbf{x}_d)\cap\mathcal{R}_i = \emptyset$ for all $\mathbf{x}\in\mathcal{M}_{\epsilon}\setminus\mathcal{R}_i$, it follows that if there exists $\mathbf{x}(t_2)\in\mathcal{M}_{\epsilon}\setminus\mathcal{R}_i$, then $\mathbf{x}(t)\notin\mathcal{R}_i$ for all $t\geq t_2$, where $i\in\mathbb{I}\setminus\mathbb{I}_a$.
This completes the proof of Claim \ref{lemma:claim2} of Lemma \ref{lemma:behaviour_of_x_in_Ri}.

\subsection{Proof of Fact \ref{fact:for_vi}}
\label{proof:fact:for_vi}
Since $\epsilon< \bar{\epsilon}$, as stated in Section \ref{section:feedback_control_design}, one has $\mathbf{x}_d\notin\mathcal{N}_{\epsilon}(\mathcal{U})$.
Therefore, by \eqref{definition:Ri}, $\mathbf{x}_d\notin\mathcal{R}_i$ for every $i\in\mathbb{I}_a$.
Additionally, by \eqref{definition:Ri_forspecialindex}, $\mathbf{x}_d\notin\mathcal{R}_i$ for $i\in\mathbb{I}\setminus\mathbb{I}_a$.
Moreover, by \eqref{set_definition_Vi_Zi}, $-\mathbf{x}_d\notin\mathcal{V}_i$ for any $i\in\mathbb{I}$.
Furthermore, since $\mathcal{V}_i\cap\mathcal{R}_i$ is a closed subset of $\mathbb{S}^n$ and $\mathbb{S}^n$ is compact in $\mathbb{R}^{n+1}$, the set $\mathcal{V}_i\cap\mathcal{R}_i$ is compact.
Therefore, since $\mathbf{P}(\mathbf{x})$ is continuous over $\mathcal{V}_i\cap\mathcal{R}_i$, there exists $c_i > 0$ such that \[\underset{\mathbf{x}\in\mathcal{V}_i\cap\mathcal{R}_i}{\max}-\mathbf{x}_d^\top\mathbf{P}(\mathbf{x})\mathbf{x}_d\leq -c_i < 0.\]
Consequently, for each $i\in\mathbb{I}$, to show the existence of $\zeta_i > 0$ such that $\dot{d}_s(\mathbf{x}, \mathbf{x}_d) \leq -\zeta_i$ over $\mathcal{V}_i\cap\mathcal{R}_i$, it is sufficient to show that there exists $\bar{\lambda}_i > 0$ such that the control input satisfies
\begin{equation}\label{sufficient_condition_for_upper_bound}\mathbf{P}(\mathbf{x})\mathbf{u}(\mathbf{x}) = \lambda_i(\mathbf{x})\mathbf{P}(\mathbf{x})\mathbf{x}_d \;\text{for some}\;\lambda_i(\mathbf{x}) \geq \bar{\lambda}_i,
\end{equation}
for all $\mathbf{x}\in\mathcal{V}_i\cap\mathcal{R}_i$.

Using \eqref{set_definition_Vi_Zi}, \eqref{definition:Ri} and \eqref{definition:Ri_forspecialindex}, we partition the set $\mathcal{V}_i\cap\mathcal{R}_i$ as follows:
\begin{equation}\label{partition1}\mathcal{V}_i\cap\mathcal{R}_i = \left(\mathcal{V}_i\cap\mathcal{R}_i\cap\mathcal{M}_{\epsilon}\right)\cup\left(\mathcal{V}_i\cap\mathcal{N}_{\epsilon}(\mathcal{U}_i)\right),\end{equation}
where the set $\mathcal{M}_{\epsilon}$ is obtained by replacing $p$ with $\epsilon$ in \eqref{eroded_free_space_definition}.
First, consider the case where $\mathbf{x}\in\mathcal{V}_i\cap\mathcal{R}_i\cap\mathcal{M}_{\epsilon}$.
Since $\mathbf{x}\in\mathcal{M}_{\epsilon}$, according to \eqref{n-sphere-control-law}, $\mathbf{u}(\mathbf{x}) = k_1\mathbf{x}_d$ for all $\mathbf{x}\in\mathcal{V}_i\cap\mathcal{R}_i\cap\mathcal{M}_{\epsilon}$ for any $i\in\mathbb{I}$.
Therefore, for each $i\in\mathbb{I}$ and for every $\mathbf{x}\in\mathcal{V}_i\cap\mathcal{R}_i\cap\mathcal{M}_{\epsilon}$, $\mathbf{P}(\mathbf{x})\mathbf{u}(\mathbf{x}) = \lambda_i(\mathbf{x})\mathbf{P}(\mathbf{x})\mathbf{x}_d$ with $\lambda_i(\mathbf{x}) = k_1$.

Now, we consider the case where $\mathbf{x}\in\mathcal{V}_i\cap\mathcal{N}_{\epsilon}(\mathcal{U}_i)$ for some $i\in\mathbb{I}$.
According to \eqref{set_definition_Vi_Zi}, $\mathcal{V}_i=\mathcal{G}(\mathbf{x}_d, \mathbf{g}_i)\cup\mathcal{G}(\mathbf{x}_d, -\mathbf{g}_i)$.
It follows that the vectors $\mathbf{P}(\mathbf{x})\mathbf{x}_d$ and $\mathbf{P}(\mathbf{x})\mathbf{g}_i$ are linearly dependent for all $\mathbf{x}\in\mathcal{V}_i\cap\mathcal{N}_{\epsilon}(\mathcal{U}_i)$.
Using \eqref{set_definition_Vi_Zi}, we partition the set $\mathcal{V}_i\cap\mathcal{N}_{\epsilon}(\mathcal{U}_i)$ as follows:
\begin{equation}\label{partition2}
    \mathcal{V}_i\cap\mathcal{N}_{\epsilon}(\mathcal{U}_i) = \mathcal{N}_{\mathcal{V}_i}^+\cup\mathcal{N}_{\mathcal{V}_i}^-,
\end{equation}
where the sets $\mathcal{N}_{\mathcal{V}_i}^+$ and $\mathcal{N}_{\mathcal{V}_i}^-$ are defined as
\begin{equation*}
\begin{aligned}
    \mathcal{N}_{\mathcal{V}_i}^+ &=\mathcal{G}(\mathbf{x}_d, \mathbf{g}_i)\cap\mathcal{N}_{\epsilon}(\mathcal{U}_i),\\
    \mathcal{N}_{\mathcal{V}_i}^- &=\mathcal{G}(\mathbf{x}_d, -\mathbf{g}_i)\cap\mathcal{N}_{\epsilon}(\mathcal{U}_i).
    \end{aligned}
\end{equation*}

First, we analyze the case where $\mathbf{x}\in\mathcal{N}_{\mathcal{V}_i}^+$.
Since $\epsilon< \bar{\epsilon}$, as stated in Section \ref{section:feedback_control_design}, one has $\mathbf{x}_d\notin\mathcal{N}_{\epsilon}(\mathcal{U}_i)$.
Moreover, since $\mathbf{g}_i\in\mathcal{U}_i^{\circ}$, one has $\mathbf{g}_i\notin\mathcal{N}_{\epsilon}(\mathcal{U}_i)$.
Furthermore, we know that $\mathcal{G}(\mathbf{x}_d, \mathbf{g}_i)\cap\{-\mathbf{x}_d, -\mathbf{g}_i\} = \emptyset$.
Consequently, $\mathcal{N}_{\mathcal{V}_i}^+\cap\{\mathbf{x}_d, -\mathbf{x}_d, \mathbf{g}_i, -\mathbf{g}_i\} = \emptyset$.
Therefore, one can show that for every $\mathbf{x}\in\mathcal{N}_{\mathcal{V}_i}^+$,  there exists $q_i^+(\mathbf{x}) > 0$ such that $\mathbf{P}(\mathbf{x})\mathbf{g}_i = -q_i^+(\mathbf{x})\mathbf{P}(\mathbf{x})\mathbf{x}_d$.
Therefore, by \eqref{n-sphere-control-law} and \eqref{individual_control_input_vector_design}, for every $\mathbf{x}\in\mathcal{N}_{\mathcal{V}_i}^+$, $\mathbf{P}(\mathbf{x})\mathbf{u}(\mathbf{x}) = \lambda_i^+(\mathbf{x})\mathbf{P}(\mathbf{x})\mathbf{x}_d,$
where $\lambda_i^+(\mathbf{x})$ is given by
\begin{equation}\label{plus_lambda_function}
    \lambda_i^+(\mathbf{x}) = \frac{k_1}{\epsilon\kappa}\left(\kappa d_s(\mathbf{x}, \mathcal{U}_i) + q_i^+(\mathbf{x})(\epsilon - d_s(\mathbf{x}, \mathcal{U}_i))\right),
\end{equation}
with $\kappa > 0$ and $d_s(\mathbf{x}, \mathcal{U}_i)\in[0, \epsilon]$.

One can represent $q_i^+(\mathbf{x})$ as $q_i^+(\mathbf{x}) = -\frac{\mathbf{x}_d^\top\mathbf{P}(\mathbf{x})\mathbf{g}_i}
{\|\mathbf{P}(\mathbf{x})\mathbf{x}_d\|^2}$, where, by construction, $q_i^+(\mathbf{x}) > 0$ over $\mathcal{N}_{\mathcal{V}_i}^+$.
Since $\mathbf{x}_d, -\mathbf{x}_d\notin\mathcal{N}_{\mathcal{V}_i}^+$, it holds that $\mathbf{P}(\mathbf{x})\mathbf{x}_d\ne \mathbf{0}$ for any $\mathbf{x}\in\mathcal{N}_{\mathcal{V}_i}^+$.
Therefore, since $\mathbf{P}(\mathbf{x})$ is continuous over $\mathcal{N}_{\mathcal{V}_i}^+$, it follows that $q_i^+(\mathbf{x})$ is continuous and well-defined over $\mathcal{N}_{\mathcal{V}_i}^+$.
Furthermore, since $\mathcal{N}_{\mathcal{V}_i}^+$ is a closed subset of $\mathbb{S}^n$, it is compact.
Consequently, since $q_i^+(\mathbf{x}) > 0$ over $\mathcal{N}_{\mathcal{V}_i}^+$, there exists $\underline{q}_i^+ > 0$ such that $\underline{q}_i^+ = \underset{\mathbf{x}\in\mathcal{N}_{\mathcal{V}_i}^+}{\min}q_i^+(\mathbf{x})$. 
It is straightforward to verify that for each $i\in\mathbb{I}$, for all $\mathbf{x}\in\mathcal{N}_{\mathcal{V}_i}^+$, $\lambda_i^+(\mathbf{x}) \geq \min\left\{k_1, \frac{k_1\underline{q}_i^+}{\kappa}\right\}$, where $\lambda_i^+(\mathbf{x})$ is defined in \eqref{plus_lambda_function}.

We proceed to analyze the case where $\mathbf{x}\in\mathcal{N}_{\mathcal{V}_i}^-$ for some $i\in\mathbb{I}$.
One can show that for every $\mathbf{x}\in\mathcal{N}_{\mathcal{V}_i}^-$, there exists $q_i^-(\mathbf{x})\geq0$ such that $\mathbf{P}(\mathbf{x})\mathbf{g}_i = q_i^-(\mathbf{x})\mathbf{P}(\mathbf{x})\mathbf{x}_d$.
Therefore, for every $\mathbf{x}\in\mathcal{N}_{\mathcal{V}_i}^-$, by \eqref{n-sphere-control-law} and \eqref{individual_control_input_vector_design}, $\mathbf{P}(\mathbf{x})\mathbf{u}(\mathbf{x}) = \lambda_i^-(\mathbf{x})\mathbf{P}(\mathbf{x})\mathbf{x}_d,$
where $\lambda_i^-(\mathbf{x})$ is given by
\begin{equation}\label{minus_lambda_function}
    \lambda_i^-(\mathbf{x}) = \frac{k_1}{\epsilon\kappa}\left(\kappa d_s(\mathbf{x}, \mathcal{U}_i)-q_i^-(\mathbf{x})(\epsilon - d_s(\mathbf{x}, \mathcal{U}_i))\right),
\end{equation}
and $q_i^-(\mathbf{x}) \geq 0$ over $\mathcal{N}_{\mathcal{V}_i}^-$.
Since $\frac{k_1}{\epsilon\kappa} > 0$, a sufficient condition for $\lambda_i^-(\mathbf{x}) > 0$ over $\mathcal{N}_{\mathcal{V}_i}^-$ is given by
\begin{equation}\label{expression_intermediate}
    \kappa > q_i^-(\mathbf{x})\left(\frac{\epsilon - d_s(\mathbf{x}, \mathcal{U}_i)}{d_s(\mathbf{x}, \mathcal{U}_i)}\right).
\end{equation}
We proceed to obtain the upper bound on the right-hand side of \eqref{expression_intermediate} over $\mathcal{N}_{\mathcal{V}_i}^-$.

Since $\mathbf{g}_i\in\sigma(\mathcal{U}_i)$, as stated in Section \ref{section:feedback_control_design}, one has $-\mathbf{g}_i\notin\mathcal{U}_i$.
Furthermore, as mentioned earlier $\mathbf{x}_d\notin\mathcal{R}_i$.
Therefore, one can show the existence of $\mu_i > 0$ such that $d_s(\mathcal{G}(\mathbf{x}_d, -\mathbf{g}_i), \mathcal{U}_i) > \mu_i$.
Consequently, since $f(p) = \frac{\epsilon - p}{p}$ is strictly decreasing on $(0, \infty)$, and $d_s(\mathbf{x}, \mathcal{U}_i)> \mu_i$ for all $\mathbf{x}\in\mathcal{N}_{\mathcal{V}_i}^-$,
it follows that
\begin{equation}\label{one_more_inequality}
    \frac{\epsilon - d_s(\mathbf{x}, \mathcal{U}_i)}{d_s(\mathbf{x}, \mathcal{U}_i)} < \frac{\epsilon - \mu_i}{\mu_i},
\end{equation}
for all $\mathbf{x}\in\mathcal{N}_{\mathcal{V}_i}^-$.

We proceed to show the existence of an upper bound on the value of $q_i^-(\mathbf{x})$ over $\mathcal{N}_{\mathcal{V}_i}^-$.
One can represent ${q}_i^-(\mathbf{x})$ as ${q}_i^-(\mathbf{x}) = \frac{\mathbf{x}_d^\top\mathbf{P}(\mathbf{x})\mathbf{g}_i}{\|\mathbf{P}(\mathbf{x})\mathbf{x}_d\|^2}$, where, by construction, $q_i^-(\mathbf{x}) \geq 0$ over $\mathcal{N}_{\mathcal{V}_i}^-$.
Since $\mathbf{x}_d, -\mathbf{x}_d\notin\mathcal{N}_{\mathcal{V}_i}^-$, it holds that $\mathbf{P}(\mathbf{x})\mathbf{x}_d\ne \mathbf{0}$ for any $\mathbf{x}\in\mathcal{N}_{\mathcal{V}_i}^-$.
Therefore, since $\mathbf{P}(\mathbf{x})$ is continuous over $\mathcal{N}_{\mathcal{V}_i}^-$, it follows that $q_i^-(\mathbf{x})$ is continuous and well-defined over $\mathcal{N}_{\mathcal{V}_i}^-$.
Furthermore, since $\mathcal{N}_{\mathcal{V}_i}^-$ is a closed subset of $\mathbb{S}^n$, it is compact.
Consequently, there exists $\bar{q}_i^- \geq 0$ such that $\bar{q}_i^- = \underset{\mathbf{x}\in\mathcal{N}_{\mathcal{V}_i}^-}{\max}q_i^-(\mathbf{x})$.
Therefore, by \eqref{one_more_inequality}, it follows that
\begin{equation*}
    q_i^-(\mathbf{x})\left(\frac{\epsilon - d_s(\mathbf{x}, \mathcal{U}_i)}{d_s(\mathbf{x}, \mathcal{U}_i)}\right) < \bar{\kappa}_i,
\end{equation*}
for all $\mathbf{x}\in\mathcal{N}_{\mathcal{V}_i}^-$, where $\bar{\kappa}_i = \frac{\bar{q}_i^-(\epsilon - \mu_i)}{\mu_i}$.
Therefore, for each $i\in\mathbb{I}$ and for any $\kappa > \bar{\kappa}_i$,  $\lambda_i^-(\mathbf{x}) \geq  \frac{k_1\mu_i(\kappa - \bar{\kappa}_i)}{\epsilon\kappa}$ for all $\mathbf{x}\in\mathcal{N}_{\mathcal{V}_i}^-$, where $\lambda_i^-(\mathbf{x})$ is defined in \eqref{minus_lambda_function}.
Combining this with the bounds obtained on the remaining two subsets $\mathcal{V}_i\cap\mathcal{R}_i\cap\mathcal{M}_{\epsilon}$ and $\mathcal{N}_{\mathcal{V}_i}^+$, and using \eqref{partition1} and \eqref{partition2}, it follows that for each $i\in\mathbb{I}$ and for every $\mathbf{x}\in\mathcal{V}_i\cap\mathcal{R}_i$, condition \eqref{sufficient_condition_for_upper_bound} is satisfied with $\bar{\lambda}_i = \min\left\{k_1, \frac{k_1\underline{q}_i^+}{\kappa},\frac{k_1\mu_i(\kappa - \bar{\kappa}_i)}{\epsilon\kappa} \right\}$.
This completes the proof of Fact \ref{fact:for_vi}. 

\bibliographystyle{IEEEtran}
\bibliography{reference}

\end{document}